\UseRawInputEncoding
\documentclass[12pt]{amsart}
\usepackage{amssymb}
\usepackage[all]{xy}
\usepackage[toc,page,title,titletoc,header]{appendix}
\usepackage{xcolor}

\begin{document}
\theoremstyle{plain}
\newtheorem{thm}{Theorem}[section]
\newtheorem*{thm1}{Theorem 1}
\newtheorem*{thm1.1}{Theorem 1.1}
\newtheorem*{thmM}{Main Theorem}
\newtheorem*{thmA}{Theorem A}
\newtheorem*{thm2}{Theorem 2}
\newtheorem{lemma}[thm]{Lemma}
\newtheorem{lem}[thm]{Lemma}
\newtheorem{cor}[thm]{Corollary}
\newtheorem{pro}[thm]{Proposition}
\newtheorem{propose}[thm]{Proposition}
\newtheorem{variant}[thm]{Variant}
\theoremstyle{definition}
\newtheorem{notations}[thm]{Notations}
\newtheorem{rem}[thm]{Remark}
\newtheorem{rmk}[thm]{Remark}
\newtheorem{rmks}[thm]{Remarks}
\newtheorem{defi}[thm]{Definition}
\newtheorem{exe}[thm]{Example}
\newtheorem{claim}[thm]{Claim}
\newtheorem{ass}[thm]{Assumption}
\newtheorem{prodefi}[thm]{Proposition-Definition}
\newtheorem{que}[thm]{Question}
\newtheorem{con}[thm]{Conjecture}
\newtheorem*{assa}{Assumption A}
\newtheorem*{algstate}{Algebraic form of Theorem \ref{thmstattrainv}}

\newtheorem*{dmlcon}{Dynamical Mordell-Lang Conjecture}
\newtheorem*{condml}{Dynamical Mordell-Lang Conjecture}
\newtheorem*{pdd}{P(d)}
\newtheorem*{bfd}{BF(d)}

\newtheorem*{probreal}{Realization problems}
\numberwithin{equation}{section}
\newcounter{elno}                
\def\points{\list
{\hss\llap{\upshape{(\roman{elno})}}}{\usecounter{elno}}}
\let\endpoints=\endlist
\newcommand{\SH}{\rm SH}
\newcommand{\Tan}{\rm Tan}
\newcommand{\res}{\rm res}
\newcommand{\Om}{\Omega}
\newcommand{\om}{\omega}
\newcommand{\la}{\lambda}
\newcommand{\mc}{\mathcal}
\newcommand{\mb}{\mathbb}
\newcommand{\surj}{\twoheadrightarrow}
\newcommand{\inj}{\hookrightarrow}
\newcommand{\zar}{{\rm zar}}
\newcommand{\Exc}{{\rm Exc}}
\newcommand{\an}{{\rm an}}
\newcommand{\red}{{\rm red}}
\newcommand{\codim}{{\rm codim}}
\newcommand{\Supp}{{\rm Supp}}
\newcommand{\rank}{{\rm rank}}
\newcommand{\Ker}{{\rm Ker \ }}
\newcommand{\Pic}{{\rm Pic}}
\newcommand{\Der}{{\rm Der}}
\newcommand{\Div}{{\rm Div}}
\newcommand{\Hom}{{\rm Hom}}
\newcommand{\im}{{\rm im}}
\newcommand{\Spec}{{\rm Spec \,}}
\newcommand{\Nef}{{\rm Nef \,}}
\newcommand{\Frac}{{\rm Frac \,}}
\newcommand{\Sing}{{\rm Sing}}
\newcommand{\sing}{{\rm sing}}
\newcommand{\reg}{{\rm reg}}
\newcommand{\Char}{{\rm char}}
\newcommand{\Tr}{{\rm Tr}}
\newcommand{\ord}{{\rm ord}}
\newcommand{\id}{{\rm id}}
\newcommand{\NE}{{\rm NE}}
\newcommand{\Gal}{{\rm Gal}}
\newcommand{\Min}{{\rm Min \ }}
\newcommand{\Max}{{\rm Max \ }}
\newcommand{\Alb}{{\rm Alb}\,}
\newcommand{\GL}{{\rm GL}\,}        
\newcommand{\PGL}{{\rm PGL}\,}
\newcommand{\Bir}{{\rm Bir}}
\newcommand{\Aut}{{\rm Aut}}
\newcommand{\End}{{\rm End}}
\newcommand{\Per}{{\rm Per}\,}
\newcommand{\ie}{{\it i.e.\/},\ }
\newcommand{\niso}{\not\cong}
\newcommand{\nin}{\not\in}
\newcommand{\soplus}[1]{\stackrel{#1}{\oplus}}
\newcommand{\by}[1]{\stackrel{#1}{\rightarrow}}
\newcommand{\longby}[1]{\stackrel{#1}{\longrightarrow}}
\newcommand{\vlongby}[1]{\stackrel{#1}{\mbox{\large{$\longrightarrow$}}}}
\newcommand{\ldownarrow}{\mbox{\Large{\Large{$\downarrow$}}}}
\newcommand{\lsearrow}{\mbox{\Large{$\searrow$}}}
\renewcommand{\d}{\stackrel{\mbox{\scriptsize{$\bullet$}}}{}}
\newcommand{\dlog}{{\rm dlog}\,}    
\newcommand{\longto}{\longrightarrow}
\newcommand{\vlongto}{\mbox{{\Large{$\longto$}}}}
\newcommand{\limdir}[1]{{\displaystyle{\mathop{\rm lim}_{\buildrel\longrightarrow\over{#1}}}}\,}
\newcommand{\liminv}[1]{{\displaystyle{\mathop{\rm lim}_{\buildrel\longleftarrow\over{#1}}}}\,}
\newcommand{\norm}[1]{\mbox{$\parallel{#1}\parallel$}}
\newcommand{\boxtensor}{{\Box\kern-9.03pt\raise1.42pt\hbox{$\times$}}}
\newcommand{\into}{\hookrightarrow}
\newcommand{\image}{{\rm image}\,}
\newcommand{\Lie}{{\rm Lie}\,}      
\newcommand{\CM}{\rm CM}
\newcommand{\sext}{\mbox{${\mathcal E}xt\,$}}  
\newcommand{\shom}{\mbox{${\mathcal H}om\,$}}  
\newcommand{\coker}{{\rm coker}\,}  
\newcommand{\sm}{{\rm sm}}
\newcommand{\pgcd}{\text{pgcd}}
\newcommand{\trd}{\text{tr.d.}}
\newcommand{\tensor}{\otimes}
\newcommand{\hotimes}{\hat{\otimes}}

\renewcommand{\iff}{\mbox{ $\Longleftrightarrow$ }}
\newcommand{\supp}{{\rm supp}\,}
\newcommand{\ext}[1]{\stackrel{#1}{\wedge}}
\newcommand{\onto}{\mbox{$\,\>>>\hspace{-.5cm}\to\hspace{.15cm}$}}
\newcommand{\propsubset}
{\mbox{$\textstyle{
\subseteq_{\kern-5pt\raise-1pt\hbox{\mbox{\tiny{$/$}}}}}$}}
\newcommand{\sA}{{\mathcal A}}
\newcommand{\sB}{{\mathcal B}}
\newcommand{\sC}{{\mathcal C}}
\newcommand{\sD}{{\mathcal D}}
\newcommand{\sE}{{\mathcal E}}
\newcommand{\sF}{{\mathcal F}}
\newcommand{\sG}{{\mathcal G}}
\newcommand{\sH}{{\mathcal H}}
\newcommand{\sI}{{\mathcal I}}
\newcommand{\sJ}{{\mathcal J}}
\newcommand{\sK}{{\mathcal K}}
\newcommand{\sL}{{\mathcal L}}
\newcommand{\sM}{{\mathcal M}}
\newcommand{\sN}{{\mathcal N}}
\newcommand{\sO}{{\mathcal O}}
\newcommand{\sP}{{\mathcal P}}
\newcommand{\sQ}{{\mathcal Q}}
\newcommand{\sR}{{\mathcal R}}
\newcommand{\sS}{{\mathcal S}}
\newcommand{\sT}{{\mathcal T}}
\newcommand{\sU}{{\mathcal U}}
\newcommand{\sV}{{\mathcal V}}
\newcommand{\sW}{{\mathcal W}}
\newcommand{\sX}{{\mathcal X}}
\newcommand{\sY}{{\mathcal Y}}
\newcommand{\sZ}{{\mathcal Z}}
\newcommand{\A}{{\mathbb A}}
\newcommand{\B}{{\mathbb B}}
\newcommand{\C}{{\mathbb C}}
\newcommand{\D}{{\mathbb D}}
\newcommand{\E}{{\mathbb E}}
\newcommand{\F}{{\mathbb F}}
\newcommand{\G}{{\mathbb G}}
\newcommand{\HH}{{\mathbb H}}
\newcommand{\I}{{\mathbb I}}
\newcommand{\J}{{\mathbb J}}
\newcommand{\M}{{\mathbb M}}
\newcommand{\N}{{\mathbb N}}
\renewcommand{\P}{{\mathbb P}}
\newcommand{\Q}{{\mathbb Q}}
\newcommand{\R}{{\mathbb R}}
\newcommand{\T}{{\mathbb T}}
\newcommand{\U}{{\mathbb U}}
\newcommand{\V}{{\mathbb V}}
\newcommand{\W}{{\mathbb W}}
\newcommand{\X}{{\mathbb X}}
\newcommand{\Y}{{\mathbb Y}}
\newcommand{\Z}{{\mathbb Z}}
\newcommand{\bk}{{\mathbf{k}}}
\newcommand{\bbk}{{\overline{\mathbf{k}}}}
\newcommand{\Fix}{\mathrm{Fix}}

\title[]{The existence of Zariski dense orbits for endomorphisms of projective surfaces \\(with an appendix in collaboration with Thomas Tucker)}

\author{Junyi Xie}


\address{Univ Rennes, CNRS, IRMAR - UMR 6625, F-35000 Rennes, France}

\email{junyi.xie@univ-rennes1.fr}

\thanks{The author is partially supported by project ``Fatou'' ANR-17-CE40-0002-01 and  PEPS CNRS}

\date{\today}

\bibliographystyle{plain}

\maketitle

\begin{abstract}
In this paper we prove the following theorem.
Let $f$ be a dominant endomorphism of a smooth projective surface over an algebraically closed field of characteristic $0$. If there is no nonconstant invariant rational function under $f$,
then there exists a closed point whose orbit under $f$ is Zariski dense.
This result gives us a positive answer to the Zariski dense orbit conjecture for endomorphisms of  smooth projective surfaces.

We define a new canonical topology on varieties over an algebraically closed field which has finite transcendence degree over $\mathbb{Q}$.  
We call it the adelic topology.
This topology is stronger than the Zariski topology and an irreducible variety is still irreducible in this topology.
Using the adelic topology, we propose an adelic version of the Zariski dense orbit conjecture, which is stronger than the original one and quantifies how many such orbits there are. We also prove this adelic version for endomorphisms of  smooth projective surfaces, for endomorphisms of abelian varieties, and split polynomial maps.
This yields new proofs of the original conjecture in the latter two cases.

In Appendix A, we study endomorphisms of $k$-affinoid spaces. We show that for certain endomorphisms $f$ on a $k$-affinoid space $X$, the attractor $Y$ of $f$ is a Zariski closed subset and 
the dynamics of $f$  semi-conjugates to its restriction to $Y.$  A special case of this result is used in the proof of the main theorem.

In Appendix B, written in collaboration with Thomas Tucker, we prove the Zariski dense orbit conjecture for endomorphisms of $(\P^1)^N.$
\end{abstract}

\tableofcontents

\section{Introduction}
Denote by $\bk$ an algebraically closed field of characteristic $0$.  Let $X$ be an irreducible variety over $\bk$ and let $f:X\dashrightarrow X$ be a dominant rational self-map.

A rational function $H\in \bk(X)$ is \emph{$f$-invariant} if $H\circ f=H$. 
Denote by $\bk(X)^f$ the field of $f$-invariant rational functions of $X$.
We have $\bk\subseteq \bk(X)^f\subseteq \bk(X).$


Give a point $x\in X(\bk)$, we say that its orbit under $f$ is well defined if for every $n\geq 0$, $f$ is well defined at $f^n(x)$.

\subsection{Zariski dense orbit conjecture}
The Zariski dense orbit conjecture was proposed by Medvedev and Scanlon \cite[Conjecture 5.10]{Medvdevv1} and by Amerik, Bogomolov and Rovinsky \cite{E.Amerik2011}.
\begin{con}\label{conexistszdo}(\emph{Zariski dense orbit conjecture}) Let $X$ be an irreducible variety over $\bk$ and let $f:X\dashrightarrow X$ be a dominant rational self-map.
If $\bk(X)^f=\bk$, then there exists  $x\in X(\bk)$ whose orbit is well-defined and is Zariski dense in $X$.
\end{con}

A dominant endomorphism $f$ on a projective variety $X$ is said to be \emph{polarized} if there exists an ample line bundle $L$ on $X$ satisfying $f^*L=L^{\otimes d}$ for some integer $d>1$.
Conjecture \ref{conexistszdo} strengthens the following conjecture of Zhang \cite{zhang}.
\begin{con}\label{conexistszdozhang}Let $X$ be an irreducible projective variety and let $f:X\to X$ be a polarized endomorphism defined over $\bk$. 
Then there exists a point $x\in X(\bk)$ whose orbit $O_f(x):=\{f^n(x)|\,\,n\geq 0\}$ is Zariski dense in $X$.
\end{con}
Conjecture \ref{conexistszdo} implies Conjecture \ref{conexistszdozhang}.
\begin{rem}The Zariski dense orbit conjecture is mainly interesting when $\bk$ is countable (e.g. $\overline{\Q}$); a difficulty is that $X(\bk)$ is covered by a countable union of proper subvarieties. 
One of the most important case of this conjecture is when $f$ is an endomorphism of $\P^N_{\overline{\Q}}$. 
\end{rem}

As stated, Conjecture \ref{conexistszdo} is slightly unbalanced. The nonexistence of nonconstant invariant 
rational function is birationally invariant, but it is not obvious that the existence of Zariski dense orbits is birationally invariant without assuming the dynamical Mordell-Lang Conjecture.
Thus, we propose to reformulate Conjecture \ref{conexistszdo} in the following strong form\footnote{This strong form is inspired by \cite[Conjecture 7.2]{Benedetto2019}. They also proposed a strong form of Conjecture \ref{conexistszdo}. However, their strong and the original form are equivalent for every pair $(X,f)$.}.
\begin{con}\label{conexistszdostr}Let $X$ be an irreducible variety over $\bk$ and let $f:X\dashrightarrow X$ be a dominant rational self-map. If $\bk(X)^f=\bk,$
then for every Zariski dense open subset $U$ of $X$, there exists a point $x\in X(\bk)$ whose orbit $O_{f}(x)$ under $f$ is well-defined, contained in $U$ and Zariski dense in $X$.
\end{con}
For a given pair $(X,f)$, Conjecture \ref{conexistszdostr} implies Conjecture \ref{conexistszdo}. Conversely 
Conjecture \ref{conexistszdo} for \emph{every} birational model of a pair $(X,f)$ implies Conjecture \ref{conexistszdostr}. See Section \ref{sectionszdp} for more discussion.

\medskip

\begin{defi}
We say that $(X,f)$ satisfies the \emph{ZD-property}, if 
either $\bk(X)^f\neq \bk$ or there exists a point $x\in X(\bk)$ whose orbit $O_{f}(x)$ under $f$ is well-defined and Zariski dense in $X$.
\end{defi}

\begin{defi}
We say that $(X,f)$ satisfies the \emph{strong ZD-property}, if
either  $\bk(X)^f\neq \bk$  or for every Zariski dense open subset $U$ of $X$, there exists a point $x\in X(\bk)$ whose orbit $O_{f}(x)$ under $f$ is well-defined, contained in $U$ and  Zariski dense in $X$.
\end{defi}

A pair $(X,f)$ satisfies the ZD-property (resp. strong ZD-property) if and only if Conjecture \ref{conexistszdo} (resp. Conjecture \ref{conexistszdostr}) holds for it.
The strong ZD-property implies the ZD-property.

\medskip

The following result is directly implied by our Theorem \ref{thmzaridenseorbitsurfendoadelic}.

\begin{cor}\label{corzaridenseorbitsurfendo}
If $f: X\to X$ is a dominant endomorphism of an irreducible smooth projective surface over $\bk$, then $(X,f)$ satisfies the strong ZD-property.
\end{cor}

\subsubsection{Historical note}
When $\bk$ is uncountable, Conjecture \ref{conexistszdo} was proved by Amerik and Campana \cite{Amerik2008}. If $\bk$ is countable, it 
has only been proved in a few special cases.
By Medvedev and Scanlon in \cite{Medvdev}, Conjecture \ref{conexistszdo} was proved when $f:=(f_1(x_1),\cdots,f_N(x_N))$ is a split polynomial endomorphism of $\mathbb{A}^N_{\bk}$.
By Bell, Ghioca and Tucker in \cite{Bell2016}, Conjecture \ref{conexistszdo} was proved when $f:=(f_1(x_1),\cdots,f_N(x_N))$ is a split endomorphism of $(\mathbb{P}^1_{\overline{\Q}})^N$
such that the $f_i$'s are endomorphisms of $\P^1_{\overline{\Q}}$ which are not post-critically finite. 
In \cite{E.Amerik2011}, Amerik, Bogomolov and Rovinsky proved Conjecture \ref{conexistszdo} under the assumption that $\bk=\overline{\Q}$ and $f$ has a fixed point $o$ which is smooth and such that the eigenvalues of $df|_o$ are nonzero and multiplicatively independent. 
The author proved Conjecture \ref{conexistszdo} for birational self-maps of surfaces with dynamical degree greater than $1$ in \cite{Xie2015}, and  for all polynomial endomorphisms $f$ of $\mathbb{A}^2$ in \cite{Xie2017}.

In \cite{Ghioca2018a}, Ghioca and the author shows that if Conjecture \ref{conexistszdostr} holds for a rational self-map $g: X\dashrightarrow X$ of a projective variety then Conjecture \ref{conexistszdo}  holds for 
all skew-linear rational self-maps $f: X\times \A^N \dashrightarrow X\times \A^N$ of type $(x,y)\mapsto (f(x), A(x)y)$ where $A(x)\in M_{N\times N}(\bk(X))$.\footnote{
The original statement of \cite[Theorem 1.4]{Ghioca2018a} said that if $(X,g)$ satisfies the strong ZD-property, then the pair $(X\times\A^N,f)$ satisfies the strong ZD-property.
But its proof only showed the slightly weaker conclusion as we stated here.}

We mention that in \cite{Amerik}, Amerik proved that there exists a nonpreperiodic algebraic point when  $f$ is of infinite order. In \cite{Bell}, Bell, Ghioca and Tucker proved that if $f$ is an automorphism without nonconstant invariant rational function, then there exists a
subvariety of codimension $2$ whose orbit under $f$ is Zariski dense. See \cite{E.Amerik2011,Fakhruddin2014,Bell2017,Bell2017a,Ghioca2017a,Ghioca2018b,Ghioca2019} for more results.

\subsection{Adelic topology}
Assume that the transcendence degree of $\bk$ over $\Q$ is finite.
We define a new canonical topology on $X(\bk)$ and call it the \emph{adelic topology}.
This topology is defined by considering all embeddings of $\bk$ in $\C$ and $\C_p$ for every prime $p$.
See Section \ref{sectionadelictop} for the precise definition. 
The terminology  ``adelic topology" usually means the topology on $X(\bf{A})$ coming from the topology on the adele ring $\bf{A}$ (see \cite[Chapter 1.2]{Weil1995} and \cite{Lorscheid2016}). This topology is different from what we define in this article. For example, the ``adelic topology" on $X(\bf{A})$ in the usual sense is Hausdorff, but the adelic topology on $X(\bk)$ we define is not Hausdorff in general.

The adelic topology has the following basic properties.
\begin{points}
\item It is stronger than the Zariski topology.
\item It is ${\mathsf{T}}_1$ i.e. for every distinct points $x, y\in X(\bk)$ there are adelic open subsets $U,V$ of $X(\bk)$ such that 
$x\in U, y\not\in U$ and $y\in V, x\not\in V.$
\item Morphisms between algebraic varieties over $\bk$ are continuous for the adelic topology.
\item Flat morphisms are open w.r.t. the adelic topology.
\item The irreducible components of $X(\bk)$ in the Zariski topology are the irreducible components of $X(\bk)$ in the  adelic topology.
\item Let $K$ be any subfield of $\bk$ which is finitely generated over $\Q$ and such that $X$ is defined over $K$ and $\overline{K}=\bk$. Then the action 
$$\Gal(\bk/K)\times X(\bk)\to X(\bk)$$ is continuous w.r.t. the adelic topology.
\end{points}
When $X$ is irreducible, (v) implies that the intersection of finitely many nonempty adelic open subsets of $X(\bk)$ is nonempty.
So, if $\dim X\geq 1$, the adelic topology is not Hausdorff.

In general, the adelic topology is strictly stronger than the Zariski topology. For example, on $\P^1(\bk)$ there exists an adelic open subset $U$ such that both $U$ and $\P^1(\bk)\setminus U$ are infinite.

\rem In \cite[Section 2.1]{Yuan2013}, Yuan and Zhang introduced an
analytic space $\sX^{\an}$, which is the union of the Berkovich
spaces over all places.  As with the adelic topology, the space
$\sX^{an}$ is canonical and is defined by considering all valuations.
On the other hand, the topology of $\sX^{an}$ is Hausdorff and does not induce a topology on $X(\bk)$, whereas the
adelic topology is on $X(\bk)$ and is not Hausdorff in general.  \endrem

\subsubsection{Adelic general points}\label{subsubade}
We have the notions of ``general" and ``very general" points to speak of  ``almost all" points in $X(\bk).$
We say  ``a property $P$ holds for a general point" if there exists a Zariski dense open subset $U$ of $X$, such that $P$ holds for all points in $U(\bk);$
and ``a property $P$ holds for a very general point" if there exists a sequence of  Zariski dense open subsets $U_n, n\geq 0$ of $X$, such that $P$ holds for all points in $\cap_{n\geq 0}U_n(\bk).$
These two notions are stable under  finite intersections. When $\bk$ is countable, there may be no $\bk$-points satisfying some given property $P$ even if $P$ holds for a very general point.
Using the adelic topology, we get a third natural notion: We say that ``property $P$ holds for an \emph{adelic general point}" if there exists a nonempty adelic open subset $U$ of $X(\bk)$, such that $P$ holds for all points in $U.$

\subsection{Adelic version of the Zariski dense orbit conjecture}
A phenomenon is that, if we have one Zariski dense orbit, we expect many such orbits. However, the Zariski topology is too weak to describe  such phenomena.
\begin{exe}Let $f: \P^1_{\bk}\to \P^1_{\bk}$ be an endomorphism of degree at least $2$. There are many Zariski dense orbits. On the other hand, the set of periodic points  of $f$ is Zariski dense in $\P^1_{\bk}$. 
However, we can show that the orbit of an adelic general point $x\in X(\bk)$ is Zariski dense in $\P^1(\bk).$
\end{exe}

In order to describe  this phenomenon, we propose an adelic version of the Zariski dense orbit conjecture.
\begin{con}\label{conexistszdoade}
Assume that the transcendence degree of $\bk$ over $\Q$ is finite.
Let $X$ be an irreducible variety over $\bk$ and let $f:X\dashrightarrow X$ be a dominant rational self-map. If $\bk(X)^f= \bk$,
then there exists a nonempty adelic open subset $A\subseteq X(\bk)$ such that the orbit of every point $x\in A$ is 
well-defined and Zariski dense in $X$.
\end{con}
In Section \ref{sectionadelictop}, we will show that this conjecture
has good behavior under birational conjugacy and is stronger
than the original version of the Zariski dense orbit conjecture.

\begin{defi}
We say that a pair $(X,f)$ satisfies the \emph{SAZD-property}, if 
there exists a nonempty  adelic open  subset $A$ of $X(\bk)$ such that for every  
point $x\in A$, its orbit $O_{f}(x)$ under $f$ is well-defined and Zariski dense in $X$.

We say that a pair $(X,f)$ satisfies the \emph{AZD-property}, if 
either $\bk(X)^f\neq \bk$ or it satisfies the SAZD-property.
\end{defi}

%
%

\rem
The SAZD-property implies the AZD-property.
\endrem

\subsubsection{Applications of the adelic topology}

Most known results on the Zariski dense orbit conjecture can be generalized to the adelic version. In Section \ref{sectionadelictop}, we generalize results which  are needed in the proof of 
Theorem \ref{thmzaridenseorbitsurfendoadelic}.

As first application of the adelic topology, we give quick proofs of Conjecture \ref{conexistszdoade} for split polynomial endomorphisms on $(\A^1)^N$ and for endomorphisms of abelian varieties.

\begin{thm}\label{thmspendopoly}Assume that the transcendence degree of $\bk$ over $\Q$ is finite. Let $f: \A^N\to \A^N, N\geq 1$ be a dominant endomorphism over $\bk$ taking form $(x_1,\dots,x_N)\mapsto (f_1(x_1),\dots, f_N(x_N)).$
Then $(X,f)$ satisfies the AZD-property.
Moreover, if $\deg f_i\geq 2,$ for every $i=1,\dots,N$, then $(X,f)$ satisfies the SAZD-property.
\end{thm}

\begin{thm}\label{thmendoabadelic}Assume that the transcendence degree of $\bk$ over $\Q$ is finite.
Let $A$ be an abelian variety over $\bk.$ Let $f:A\to A$ be a dominant endomorphism.
Then $(A,f)$ satisfies the AZD-property.
\end{thm}
These results generalize and give new proof of \cite[Theorem 7.16]{Medvdev} and \cite[Theorem 1.2]{Ghioca2017a}.
In their original proof,  Medvedev and Scanlon used their deep 
 result on the classification of all invariant subvarieties of split polynomial endomorphisms on $(\A^1)^N$, which is proved using Model theory and polynomial decomposition theory.
The original proof of \cite[Theorem 1.2]{Ghioca2017a} relies on the Mordell-Lang conjecture, due to Faltings. 
Our proof is different in nature, more conceptual, and does not rely on any classification argument, or Falting's theorem.

\subsection{Endomorphisms of surfaces}
In this paper, we prove Conjecture \ref{conexistszdoade} for endomorphisms of smooth projective surfaces. 
\begin{thm}\label{thmzaridenseorbitsurfendoadelic}
Assume that the transcendence degree of $\bk$ over $\Q$ is finite.
Let $X$ be an irreducible smooth projective surface over $\bk$.
Let $f:X\to X$ be a dominant endomorphism.
Then $(X,f)$ satisfies the AZD-property.
\end{thm}

For the study of Conjecture \ref{conexistszdostr} and Conjecture \ref{conexistszdo}, the hypothesis on $\bk$ does not cause any problem. 
There exists an algebraically closed subfield $K$ of $\bk$ whose transcendence degree over $\Q$ is finite and  such that $(X,f)$ is defined over
$K$ i.e. there exists a pair $(X_K,f_K)$ such that $(X,f)$ is its base change by $\bk$. 
Corollary \ref{coradelicstrong} shows that if $(X_K,f_K)$ satisfies the ADZ-property, then $(X,f)$ satisfies the strong ZD-property.
Then Theorem \ref{thmzaridenseorbitsurfendoadelic} implies Corollary \ref{corzaridenseorbitsurfendo}.



\subsubsection{Strategy}
For simplicity, we explain a strategy of  a direct proof of Conjecture \ref{conexistszdo} for endomorphisms of smooth projective surfaces.
However, except for the systematic use of the adelic topology and some technical difficulties, Theorem \ref{thmzaridenseorbitsurfendoadelic} follows the same idea.


\medskip

We first explain our strategy for an endomorphism $f$ of  $\P^2$ of degree at least $2.$
For simplicity, we assume that $\bk=\overline{\Q}.$
In this case, there is no nonconstant rational function invariant under $f.$ So we only need to show that there exists a closed point which has a Zariski dense orbit.
The idea is to combine the $p$-adic local dynamics near a certain periodic point with a constraint, 
obtained by some global information, on the field of definition of any invariant curve.

By \cite{E.Amerik2011}, if there exists a fixed point $o$ of $f^m,$ for some $m\geq 1$ such that the two eigenvalues $\la_1,\la_2$ of $df^m|_o$ are multiplicatively independent, then there exists a closed point which has a Zariski dense orbit. So we may assume that such points do not exist.

At first, we study the $f$-invariant curves.
A curve $C$ is $f$-invariant, if $f(C)=C.$
There is a number field $K$ such that $f$ and all fixed points of $f$ are defined over $K.$
We show that there exists a positive integer $N$ depending on $f$, such that every irreducible $f$-invariant curve $C$ is defined over a field $K_C$ satisfying $[K_C:K]| N^n$ for some $n\geq 0.$
Moreover, we show that the number of invariant branches of $C$ at any fixed point of $f$ is bounded from above by some integer $B>N.$

Next we want to find a fixed point $o$ of $f^m$ for some  $m\geq 1$ and a field embedding $\tau:\overline{\Q} \hookrightarrow \C_p$ 
such that 
\begin{points}
\item $df^m|_o$ is invertible;
\item $|\tau(\la_1)|$ and $|\tau(\la_2)|\leq 1$ where
 $\la_1,\la_2\in \overline{\Q}$ are the two eigenvalues of $df^m|_o$;
 \item  $|\tau(\la_1)||\tau(\la_2)|<1.$ 
\end{points}
By studying multipliers of endomorphisms on curves and 
assuming that there are no Zariski dense orbits of closed point, we show that
the existence of such a point is ensured by the existence of a repelling periodic point.
The latter is ensured by \cite[Theorem 3.1, iv)]{Guedj2005} (see \cite[Theorem 1.1]{Dinh2015}  and \cite[Theorem 2]{Briend1999} also).
After replacing $f$ by $f^m$, we assume that $o$ is a fixed point of $f$.
Using $\tau$, we view $\overline{\Q}$ as a subfield of $\C_p.$
We assume that $\la_1,\la_2$ are contained in $K.$ Denote by $K_p$ the closure of $K$ in $\C_p.$
 
If $|\la_1|=1$ and $|\la_2|<1$, we show that there exists an $f$-invariant $p$-adic neighborhood $U$ of $o$ in $\P^2(K_p)$ which is isomorphic to a polydisc $(K_p^{\circ})^2$.
Moreover, after shrinking $U$, there exist an analytic curve $Y\subseteq U$ which is invariant by $f$, and an analytic morphism $\psi: U\to Y$ such that 
$\psi|_Y=\id$, $\psi\circ f=f|_Y\circ \psi,$
and $\cap_{n\geq 0}f^n(U)=Y$.  Appendix A contains a more general result for endomorphisms of affinoid spaces.
For an endomorphism of $\P^2$ of degree at least two, the periodic points are isolated, hence $f|_Y$ is not of finite order. In this case, we then show that there is a point in $\P^2(\overline{\Q})\cap U$ whose orbit is Zariski dense.

If both $|\la_1|$ and $|\la_2|$ are $<1$, then, since $\la_1,\la_2$ are not multiplicatively independent, there exists $m_1,m_2\geq 1$ such that 
$\la_1^{m_1}=\la_2^{m_2}.$ After replacing $f$ by a suitable iterate, we assume $(m_1,m_2)=1.$ We show that there exist a birational morphism $\pi: X'\to \P^2$ which is a composition of blowups of $K$-points, an irreducible component $E$ in $\pi^{-1}(o)$ such that the induced rational self-map $f'$ on $X'$ is regular along $E$ and fix $E$, and
a fixed point $o'\in E(K)$ such that the two eigenvalues of  $df'|_{o'}$ are $1$ and $\mu$ with $|\mu|<1$ and the eigenvectors for $1$ are in the tangent space of $E.$
 Let $M$ be a finite field extension of $K$ such that $[M:K]$ is prime to $B!.$
 Denote by $M_p$ the closure of $M$ in $\C_p.$ The argument in the previous paragraph shows that 
 there exists a $p$-adic neighborhood $U$ of $o$ in $X'(M_p)$ which is isomorphic to a polydisc $(M_p^{\circ})^2$ and invariant by $f'$ satisfying 
 $\cap_{n\geq 0}(f')^n(U)=Y:=U\cap E$
and an analytic morphism $\psi: U\to Y$ such that 
$\psi|_Y=\id$ and $\psi\circ f'=f'|_Y\circ \psi.$
Moreover the construction of $U$ and $\psi$ shows that they are defined over $K_p.$
If $f'|_Y\neq \id$, we may conclude the proof by the argument in the previous paragraph.
If $f'|_Y= \id$, such an argument is not sufficient. Here we need the constraint on definition fields of invariant curves.
We show that every irreducible periodic curve $C$ passing through $U$ is indeed invariant by $f'$. 
So it is defined over a field $K_C$ such that $[K_C:K]| N^n$ for some $n\geq 0.$
It follows that $C\cap U$ is defined over the closure $(K_C)_p$ of $K_C$ in $\C_p.$
We show that $C\cap U$ is a disjoint union of $\psi^{-1}(x_i), i=1,\dots,s$ where $s\leq B$ and $x_i\in Y=U\cap E$.
Then there exists a field extension $H_p$ of $K_p$ satisfying $[H_p: (K_C)_p]| B!$ such that all 
 $x_i$ are defined over $H_p$. It follows that there exists $n\geq 0$ such that $[H_p: K_p]| (B!)^n.$
 Since $[M_p:K_p]$ is prime to $B!$, $M_p\cap H_p=K_p.$ Then $x_i\in X'(K_p)\cap Y, i=1,\dots,s.$
 Since $X'(K_p)\cap Y$ is not dense in $Y,$ we deduce that there exists a point $x\in X'(\overline{\Q})\cap \psi^{-1}(Y\setminus X'(K_p))$ which has a Zariski dense orbit for $f'$.
 Then $\pi(x)$ has a Zariski dense orbit for $f$.

In the general case, by  \cite[Theorem 1.3]{Bell}, we may assume that $f$ is not an automorphism.
Using the classification of surface and the works of Fujimoto, Nakayama,  Matsuzawa, Sano, and Shibata on endomorphisms of surfaces, we may reduce to a case that can either be treated by the same 
argument for $\P^2$ or preserves a fibration to a curve.
In the latter case, we conclude the proof using this fibration.

\subsection{Endomorphisms of $(\P^1)^N$}
In Appendix B,  we prove
Conjecture \ref{conexistszdoade} for endomorphism of $(\P^1)^N$ which generalizes Theorem \ref{thmspendopoly}.
As a consequence, this settles the Zariski dense orbit conjecture for endomorphisms of $(\P^1)^N$.
\begin{thm}\label{thmzaridenseorbitsurfendop1n}
Let $f:(\P^1)^N\to (\P^1)^N, N\geq 1$ be a dominant endomorphism over $\bk$.
Then $((\P^1)^N, f)$ satisfies the strong ZD-property.
\end{thm}

This result generalizes \cite[Theorem 7.16]{Medvdev} and \cite[Theorem 14.3.4.2]{Bell2016}.

%
%
%

Now we explain the strategy of the proof of Theorem \ref{thmzaridenseorbitsurfendop1n}. 
We only need to prove its adelic version stated as Theorem \ref{thmzaridenseorbitsurfendop1nadelic}.
It is easy to show that after replacing $f$ by a positive iterate, we may assume that $f$ takes form 
$f:=(f_1(x_1),\cdots,f_N(x_N)).$
As in of \cite{Medvdev}, we first need a description of invariant subvarieties of endomorphisms of $(\P^1)^N$ (see Proposition \ref{proinvsubvspl}). Basically this description shows that when all factors $f_i$ are not of some special type,  all invariant subvarieties come from some invariant curves of $(f_i,f_j):\P^1\times \P^1\to \P^1\times \P^1$ for some $i\neq j.$
Such a description was already obtained by Medvedev and Scanlon in \cite{Medvdev} using Model theory. Here we give a new and elementary proof.
Using this description and some basic properties of the adelic topology, we reduce the problem to the case $N=2.$  Then we may conclude the proof by Theorem \ref{thmzaridenseorbitsurfendoadelic}.

\subsection{Further problems and potential applications}\label{subsectionreal}
Inspired by the Zariski dense orbit conjecture, we introduce a new class of problems.
\begin{probreal}
Let $X$ be a variety over an algebraically closed field $k$,
Let $S$ be an algebraic structure on $X$ (e.g. $\emptyset$, a rational self-map, a line bundle, ect).
Let $P(k,X,S)$ be an algebraic property for points in $X(k)$ which depends on $k,X,S$. 
What is the minimal transcendence degree $R_P(k,X,S)$ of a field extension $K$ of $k$ such that $P(K,X_K,S_K)$
holds for at least one point $x\in X(K)$, if such $K$ exists. Here $X_K,S_K$ are the base change of $X,S$ by $K.$
\end{probreal}

Roughly speaking, we ask whether  a property which generally holds can be realized for at least one algebraic point.
In many cases, this problem is easy when $\bk$ is uncountable, but it could be very difficult over $\overline{\Q}.$ 
This resonates with the notion of essential dimension of Zinovy Reichstein \cite{Buhler1997}.
  
\medskip

Because the Zariski dense orbit conjecture holds over an uncountable field \cite{Amerik2008}, it can be rewritten as a realization problem by letting $k$ be an algebraically closed field of characteristic $0$, $S$ be a dominant rational self-map $f: X\dashrightarrow X$ satisfying $\bk(X)^f=\bk,$
and $P(k,X,S=\{f\})$ be the property ``the $f$-orbit of $x$ is Zariski dense in $X$"  for a point $x\in X(k)$.

One may ask a similar question in positive characteristic. In this case $R_P(k,X,f)$ depends on the pair $(X,f)$. 
Assume that $k=\overline{\F_p}$ for some prime $p.$  In this case $R_P(k,X,f)$ can not be zero since every orbit is finite. 
For example, when $X=\P^1$ and $f: z\mapsto z^2$, we have $R_P(k,X,f)=1;$
when $X=\P^2$ and $f: [x:y:z]\mapsto [x^p:y^p:z^p],$ we can show that $R_P(k,X,f)=2.$

\subsubsection{Examples of realization problems}
Changing $P$, we get many interesting questions. 

\begin{que}Let $X$ be a variety over an algebraically closed field $k$.  Assume that the property $P$ for a point  $x\in X(k)$ is ``$x$ is not contained in a rational curve." 
We may ask whether $R_P(k,X,\emptyset)=0$, when $X$ is not uniruled.
\end{que}
When $k$ is $\overline{\Q}$ or the algebraic closure of a finite field, and $X$ is a non-uniruled K3-surface over $k$,
this question was asked by Bogomolov in 1981 (see \cite{Bogomolov2005}).

\begin{que}\label{queper}Let $X$ be an irreducible variety over $k=\overline{\F_p}$ for some prime $p.$
Let $S=G$ be a finitely generated subgroup of the group $\Bir(X)$ of birational transformations of $X$.
Assume that the property $P$ for $x\in X(k)$ is ``there exists a finite index subgroup $H$ of $G$ such that for every $g\in H$, $g$ is well defined at $x$." 
We ask whether $R_P(k,X,G)=0$.
\end{que}
Note that, if Question \ref{queper} has a positive answer for a group $G$, then there exists $x\in X(k)$ and a finite index subgroup $H'$ of $G$ such that for every $g\in H'$ $g(x)=x$
(because $k=\overline{\F_p}$).
\begin{points}
\item[(1)] When $G\simeq\Z$, Question \ref{queper} has a positive answer by the twisted Lang-Weil estimates of Hrushovski \cite{hu7,fa}.
\item[(2)] When $G$ has Property (T), Cantat and the author showed that Question \ref{queper} has a positive answer \cite{Cantata} (see also \cite{Cantat2019,Cornulier,Lonjou2020} for Property (FW)).
\item[(3)] When $G\simeq\Z^2$, this question is still open.
\end{points}
The $p$-adic method introduced in \cite{Cantata} suggests that any progress on Question \ref{queper} will be useful in the study of $\Bir(X).$

\begin{que}Let $B$ be an irreducible variety over $\overline{\Q}$ and let $\pi:\sX\to B$ be a projective scheme over $B$. Let $f_B: \sX\to \sX$ be an endomorphism over $B$ and let $x_B\in \sX(B)$ be a section.
\begin{points}
\item[(1)]
Assume that $x_B$ is not preperiodic in $\sX$. Let $S$ be the triple $(\pi: \sX\to B, f_B, x_B)$.
Consider the property $P$ for $b\in B(\overline{\Q})$ given by ``in the special fiber $X_b$ of $X$ at $b$, the specialization $x_b$ of $x_B$ is not preperiodic for the specialization $f_b$ of $f_B$."
We ask whether $R_P(\overline{\Q},B, S)=0$.
\item[(2)]
Let $\eta$ be the generic point of $B.$
Assume that the generic fiber $X_{\eta}$ does not have any non-constant $f_{\eta}$-invariant rational function. Let $S$ be the pair $(\pi: \sX\to B, f_B)$.
Consider the property $P$ for $b\in B(\overline{\Q})$ is given by ``the special fiber $X_b$ does not have any non-constant $f_{b}$-invariant rational function."
We ask whether $R_P(\overline{\Q},B,S)=0$.
\end{points}
\end{que}
A special case of (1) was used in \cite{Ghioca2018} to reduce some instance of the dynamical Mordell-Lang conjecture to $k=\overline{\Q}$ and
(2) can be used to reduce the Zariski dense orbit conjecture to $k=\overline{\Q}$.
More generally, solving a realization problem helps to reduce a problem over algebraically closed fields to some specific field.

\subsubsection{Adelic version of realization problems} 
There is also a variant of the realization problem.
Assume that  $\bk$ is an algebraically closed field of finite transcendence degree over $\mathbb{Q}.$
Let $X$ be an irreducible variety over $\bk$. 
Let $S$ be an algebraic structure on $X$.
Let $P(\bk,X,S)$ be an algebraic property for points in $X(\bk)$. Then we may ask the question $Q_{P}(\bk,X,S)$: does an adelic general point in $X(\bk)$ satisfy $P(\bk,X,S)$?

\medskip

For two algebraic properties $P_1,P_2$,
in general,  $R_{P_1}(\bk,X,S)=R_{P_2}(\bk,X,S)=0$ does not imply $R_{P_1\text{ and } P_2}(\bk,X,S)=0.$
On the other hand, if the answer of both $Q_{P_1}(\bk,X,S)$ and $Q_{P_2}(\bk,X,S)$ are positive, the answer of $Q_{P_1\text{ and } P_2}(\bk, X,S)$ is also positive. 
This property allows us to split a problem into problems that are easier to tackle. 
Using the notion of "adelic general point", many realization problems have a  strengthened adelic versions.
An example is the Zariski dense orbit conjecture.
As shown in the proof of Theorem \ref{thmspendopoly}, \ref{thmendoabadelic} and  \ref{thmzaridenseorbitsurfendop1n},
though the adelic version of a realization problem is stronger than the original one, it is easier to handle in many cases.

\subsection{Organization of the paper}
Section \ref{sectionszdp} contains some basic facts about the Zariski dense orbit conjecture. 
We show that Conjecture \ref{conexistszdostr} implies Conjecture \ref{conexistszdo}. 
and discuss the relation between the Zariski dense orbit conjecture and the dynamical Mordell-Lang conjecture.

\medskip

Section \ref{sectionadelictop} introduces the adelic topology and proves some of its basic properties. Using this topology, we propose the adelic version of the Zariski dense orbit conjecture.
We show that Conjecture \ref{conexistszdoade} implies Conjecture \ref{conexistszdostr}. In particular, Theorem \ref{thmzaridenseorbitsurfendoadelic} implies Corollary \ref{corzaridenseorbitsurfendo}. 
We generalize former results on the Zariski dense orbit conjecture to their adelic version. 

\medskip

In Section \ref{sectionappadelic}, we give some applications of the adelic topology,
thereby proving 
Theorem \ref{thmspendopoly} and Theorem \ref{thmendoabadelic}.

\medskip

Section \ref{secgen} concerns endomorphisms of projective surfaces. We prove a constraint on the definition field of their invariant curves.

\medskip 

In Section \ref{seclocdy}, we first study the multipliers of periodic points and the dynamics near a fixed point;
then we focus on amplified endomorphisms. In particular, we prove Theorem \ref{thmzaridenseorbitsurfendoadelic} for endomorphisms of $\P^2.$

\medskip

In Section \ref{secproof}, we prove Theorem \ref{thmzaridenseorbitsurfendoadelic} in the general case.

\medskip
In Appendix A, we study the endomorphisms on the $\bk$-affinoid spaces. We show that for certain endomorphism $f$ on a $\bk$-affinoid space $X$, the attractor $Y$ of $f$ is a Zariski closed subset and 
the dynamics of $f$ semi-conjugates to the its restriction to $Y.$  

\medskip
In Appendix B, written in collaboration with Thomas Tucker, we prove Theorem \ref{thmzaridenseorbitsurfendop1n}.

\subsection{Notation and Terminology}
\begin{points}
\item[$\bullet$] A \emph{variety} is a reduced separated scheme of finite type over a field.
\item[$\bullet$] For a variety $X$ (resp. a rational self-map $f: X\dashrightarrow Y$) over a field $k$ and a subfield $K$ of $k$, we say that $X$ (resp. $f$) is \emph{defined over $K$} if there is a variety $X_K$ (resp. a rational map $f_K$) over $K$ such that $X$ (resp. $f$) is the base change by $k$ of $X$ (resp. $f$).
\item[$\bullet$] For a prime $p$, $\C_p$ is the completion of $\overline{\Q_p}$, $\C_p^{\circ}=\{x\in \C_p|\,\, |x|\leq 1\}$ is the valuation ring of $\C_p$ and $\C_p^{\circ\circ}=\{x\in \C_p|\,\, |x|<1\}$ is the maximal ideal of $\C_p^{\circ}.$
\item[$\bullet$]When $p=\infty$, $\C_p$ is the field $\C$  of complex numbers. 
\item[$\bullet$] For a dominant rational map $f: X\dashrightarrow Y$ between irreducible varieties over a field $k$, denote by $d_f:=[k(X):f^*k(Y)]$ the topological degree of $f.$  When $k$ is algebraically closed of characteristic $0$, $d_f=|f^{-1}(y)|$ where $y$ is a general point in $Y(k).$
\item[$\bullet$]For a rational map $f: X\dashrightarrow Y$ between varieties. Denote by $I(f)$ the indeterminacy locus of $f$. For a subvariety $V$ of $X$, if $I(f)$ does not contain any irreducible component of $V$, we write $f(V)$ for the strict transform $\overline{f(V\setminus I(f))}$ of $V$.
\item[$\bullet$] For a dominant rational self-map $f: X\dashrightarrow X$ between varieties. A subvariety $V$ of $X$ is said to be \emph{$f$-invariant} if 
$I(f)$ does not contain any irreducible component of $V$ and 
$f(V)\subseteq V.$ 
\item[$\bullet$] For a projective variety $X$, $N^1(X)$ is the $\R$-N\'eron-Severi group of $X$. 
\item[$\bullet$] For an irreducible variety $X$, $\kappa(X)$ is the Kodaira dimension of $X.$
\item[$\bullet$] For a field extension $k/K$, $\trd_Kk$ is the transcendence degree of $k/K.$
\end{points}

\subsection*{Acknowledgement}
I would like to thank Yang Cao, Laurent Moret-Bailly and Miao Niu for useful discussions.
I thank Serge Cantat, Charles Favre, Dragos Ghioca, Deqi Zhang and Shou-Wu Zhang for their comments of the first version of the paper.
I thank Yohsuke Matsuzawa and a referee, who told me Lemma \ref{lemindamplify}.
I especially thank Thomas Tucker. The discussion with him improved very much the paper. 
I'm grateful to the referees, who find a problem in my previous definition of the adelic topology, suggest me to state an algebraic version of  Theorem \ref{thmstattrainv}. \ref{lemindamplify}.

\section{The Zariski dense orbit conjecture}\label{sectionszdp}

Let $X$ be an irreducible variety over $\bk$ and $f:X\dashrightarrow X$ be a dominant rational self-map.

\begin{lem}\label{leminvratfunite}Let $X'$ be an irreducible variety over $\bk$, $f': X'\dashrightarrow X'$ be a rational self-map and $\pi: X'\dashrightarrow X$ be a generically finite dominant rational map satisfying $f\circ \pi=\pi\circ f',$ then we have the following properties.
\begin{points}
\item If there exists $m\geq 1$, and $H\in \bk(X)^{f^m}\setminus \bk$, then there exists $G\in \bk(X)^f\setminus \bk$.
\item  There exists $H'\in \bk(X')^{f'}\setminus \bk$, if and only if there exists $H\in \bk(X)^f\setminus \bk$.
\end{points}
\end{lem}

\proof[Proof of Lemma \ref{leminvratfunite}]
(i). There exists a point $b\in X(\bk)$ such that $f,\dots,f^{m-1}$ are well-defined at $b$ and $H, f^*H,\dots, (f^{m-1})^*H$ are regular at $b.$
Set $a:=H(b)\in \bk$ and 
$P:=\prod_{i=0}^{m-1}((f^i)^*H-a)\in \bk(X).$ We have $f^*P=P.$ Since $H$ is not constant, for every $i=0,\dots,m$, $(f^i)^*H-a\neq 0.$ 
So we have $P\neq 0.$ Since $P(b)=0$, $P$ is not constant, which concludes the proof.

(ii). Assume that there exists $H'\in \bk(X')^{f'}\setminus \bk$.
Set $m:=[\bk(X'):\bk(X)].$
Then $ \bk(X')$ is a $m$-dimensional $\bk(X)$ vector space. 
Denote by $T^m+\sum_{i=1}^{m}(-1)^iP_iT^{m-i}$ the characteristic polynomial of the $\bk(X)$-linear operator 
$\bk(X')\to \bk(X'): g\mapsto H'g.$
We have $P_i\in \bk(X)$ and $f^*(P_i)=P_i, i=1,\dots,m.$
If $P_i\in \bk$ for every $i=1,\dots, m$, then $H'\in \bk,$ which is a contradiction. Then there exists $i=1,\dots,m$, such that $P_i\in \bk(X)^f\setminus \bk.$
If $H\in \bk(X)^f\setminus \bk$,
then $H':=H\circ \pi\in \bk(X')^{f'}\setminus \bk$.
\endproof

The following proposition shows that the strong ZD-property is invariant under birational conjugation and iterations. 

\begin{pro}\label{proszdundbiriter}The following statements are equivalents:
\begin{points}
\item $(X,f)$ satisfies the strong ZD-property;
\item  there exists $m\geq 1$, such that $(X,f^m)$ satisfies the strong ZD-property;
\item there exists a pair $(Y,g)$ which is birational to the pair $(X,f)$, and $(Y,g)$ satisfies the strong ZD-property.
\end{points}
\end{pro}

\begin{rem}Lemma \ref{leminvratfunite} implies also that the ZD-property is invariant under iterations. However, a priori it is not clear that whether the ZD-property is invariant under birational conjugation.
\end{rem}

\proof[Proof of Proposition \ref{proszdundbiriter}]
It is clear that (i) implies (ii), (i) implies (iii) and (iii) implies (i).
We only need to prove that  (ii) implies (i).
Assume that $(X,f^m)$ satisfies the strong ZD-property.
By Lemma \ref{leminvratfunite}, we assume that $\bk(X)^{f^m}=\bk.$
Let $U$ be a Zariski dense open subset of $X$. 
We may assume that $U\cap (I(f)\cup I(f^2)\dots \cup I(f^{m-1}))=\emptyset.$ 
Set 
$V:=U\cap(\cap_{i=1}^{m-1} f|_U^{-i}(U)).$
For every point $x\in V(\bk)$, $x,f(x),\dots,f^{m-1}(x)$ are well-defined and contained in $U$.

Since $(X,f^m)$ satisfies the strong ZD-property, 
there exists a point $x\in X(\bk)$ whose orbit $O_{f^m}(x)$ under $f^m$ is contained in $V$ and is Zariski dense in $X$.
It follows that the orbit $O_{f}(x)$ under $f$ is contained in $U$ and is Zariski dense in $X$, which concludes the proof.
\endproof

The following definition was introduced in \cite{Xie2014} when $X$ is a surface.
\begin{defi}We say that a pair $(X,f)$ satisfies the \emph{DML-property}, if for every subvariety $V$ of $X$ and every $x\in X(\bk)$ whose orbit is well-defined, 
the set $\{n\geq 0|\,\, f^n(x)\in V\}$ is a finite union of arithmetic progressions\footnote{An arithmetic progression is a set of the form $\{an + b|\,\, n\in \N\}$ with $a,b \in \N$ possibly with $a = 0$.}.
\end{defi}
The dynamical Mordell-Lang conjecture was proposed by Ghioca and Tucker \cite{Ghioca2009}.
The following is a slight generalization for rational self-maps.
\begin{con}\label{condml}
 Let $X$ be a variety defined over $\bk$, $f: X\dashrightarrow X$ be a rational self-map. Then $(X,f)$ satisfies the DML-property.
\end{con}

\begin{pro}\label{prodmlzdimpszd}Assume that $(X,f)$ has both the ZD-property and the DML-property, then it has the strong ZD-property.
\end{pro}
\proof[Proof of Proposition \ref{prodmlzdimpszd}]We may assume that $\bk(X)^f=\bk$.
By the ZD-property, there exists $x\in X(\bk)$ whose orbit is well-defined and Zariski dense in $X$.

Let $U$ be any Zariski dense open subset of $X$. Set $Z:=X\setminus U.$ Since $(X,f)$ satisfies the DML-property,
the set $\{n\geq 0|\,\, f^n(x)\in Z\}$ is a finite union of arithmetic progressions. If it is infinite, there exist $a,b \in \N, a\neq 0$ such that 
$\{an + b|\,\, n\in \N\}\subseteq \{n\geq 0|\,\, f^n(x)\in Z\}.$
Then
$\{f^n(x)|\,\, n\geq b\}\subseteq Z\cup \cdots \cup f^{a-1}(Z).$
So the orbit $O_f(x)$ is not Zariski dense, which contradicts our assumption. 
So there exists $N>0$ such that $f^n(x)\not\in Z$ for every $n\geq N.$
Then the orbit of $f^N(x)$ is well-defined, contained in $U$ and Zariski dense in $X$.
\endproof

\section{The adelic topology}\label{sectionadelictop}
Assume $\trd_{\Q}\bk<\infty.$ Let $X$ be a variety over $\bk$.
Let $K_0$ be a subfield of $\bk$ such that
\begin{points}
\item[(1)] $K_0$ is finitely generated over $\Q$;
\item[(2)] $\overline{K_0}=\bk$;
\item[(3)]  $X$ is defined over $K_0$ i.e. there exists a variety $X_{K_0}$ defined over $K_0$ such that 
$X=X_{K_0}\times_{\Spec K_0}\Spec \bk.$
\end{points}
In this section, we define the adelic topology on $X(\bk).$

\subsection{The adelic topology}

For every algebraic extension $K$ over $K_0$, define $X_K=X_{K_0}\times_{\Spec K_0}\Spec K.$ We may canonically
identify $X_K(\bk)$ with $X(\bk)$.

\medskip
\subsubsection{The basic adelic subsets}
Let $K$ be any finite field extension of $K_0.$
Denote by $\sI_K$ the set of  field embeddings $\tau:K\hookrightarrow \C_{p_{\tau}}$ for some prime $p_{\tau}$  or $p_{\tau}=\infty$. 
Denote by $\sI_K^{f}$ the set of $\tau\in \sI_K$ for which $p_{\tau}$ is a prime.
We say two embeddings $\tau,\tau'\in \sI_K$ are equivalent if the absolute values $|\tau(\cdot)|, |\tau'(\cdot)|$ are the same functions on $K$.
Denote by $\sM_K$ (resp. $\sM_K^f$) the set of equivalent classes in $\sI_K^f.$
For every $\tau\in \sI_K$, denote by $\sI_{\tau}$ the set of field embeddings $\overline{\tau}:\bk=\overline{K}\hookrightarrow \C_{p_{\tau}}$ satisfying $\overline{\tau}|_{K}=\tau.$
Every $\overline{\tau}\in \sI_{\tau}$ induces an embedding $\phi_{\overline{\tau}}: X(\bk)\hookrightarrow X_K(\C_{p_{\tau}}).$
For all $\overline{\tau}\in \sI_{\tau}$, the images $\phi_{\overline{\tau}}(X(\bk))\subseteq X_K(\C_{p_{\tau}})$ are the same.
On $X_K(\C_{p_{\tau}})$, we have the $p_{\tau}$-adic topology when $p_{\tau}$ is a prime.
When $p_{\tau}=\infty$, we have the complex topology on $X_K(\C_{p_{\tau}})$ and we also call it the $p_{\tau}$-adic topology in this case.
The images $\phi_{\overline{\tau}}(X(\bk))\subseteq X_K(\C_{p_{\tau}})$ is dense in the $p_{\tau}$-adic topology.

\medskip

For $\tau\in \sI_K$ and a $p_{\tau}$-adic open subset $U$ of $X_K(\C_{p_i})$, define
$$X_K(\tau, U):=\cup_{\overline{\tau}\in \sI_{\tau}}\phi_{\overline{\tau}}^{-1}(U)\subseteq X(\bk).$$
More generally, let $\tau_i: K\hookrightarrow \C_{p_i}, i=1,\dots,m$ be (not necessarily distinct) elements of $\sI_K$,
and $U_i$ be a $p_i$-adic open subset of $X_K(\C_{p_i})$.
Define
$$X_K((\tau_i, U_i), i=1,\dots,m):=\cap_{i=1}^{m}X_K(\tau_i, U_i)\subseteq X(\bk).$$
Such a subset of $X(\bk)$ is called a \emph{basic adelic subset over $K$}.
A subset of $X$ is called a \emph{basic adelic subset}, if it is a basic adelic subset over \emph{some} finite field extension $K$ of $K_0.$

\begin{rem}\label{reminterzari}
For every Zariski open subset $U$ of $X$ defined over $K$, the subset $U(\bk)\subseteq X(\bk)$ is 
a basic adelic subset over $K$ and we have 
 $$X_K((\tau_i, U_i), i=1,\dots,m)\cap U=X_K((\tau_i, U_i\cap U_K(\C_{p_{\tau_i}})), i=1,\dots,m).$$
\end{rem}

\begin{rem}\label{remsurject}Let $Y$ be a variety over $K$ and let
$\pi:Y\to X$ be a  morphism over $K$. Set $\pi_{\tau_i}:=\pi\otimes_{K}\C_{p_{\tau_i}}: Y_K(\C_{p_{\tau_i}})\to X_K(\C_{p_{\tau_i}})$ the morphism induced by $\pi.$
Then 
$$\pi^{-1}(X_K((\tau_i, U_i), i=1,\dots,m))=Y_K((\tau_i, \pi_{\tau_i}^{-1}(U_i)), i=1,\dots,m).$$
\end{rem}

\begin{rem}\label{reminter}
We have $$X_{K}((\tau_i, U_i), i=1,\dots,m)\cap X_K((\tau_i', U_i'), i=1,\dots,m')$$
$$=X_K((\tau_i, U_i), i=1,\dots,m, (\tau_j', U_j'), j=1,\dots,m').$$
\end{rem}

\begin{rem}\label{remfieldextad}Let $K'$ be a finite extension of $K$. 
For every field embedding $\tau\in \sI_K$, the set $\sI_{\tau}^{K'}:=\{\tau'\in \sI_{K'}|\,\, \tau'|_{K'}=\tau\}$ is finite.
 Then we have $$X_K((\tau_i, U_i), i=1,\dots,m)=\cap_{i=1}^mX_K((\tau_i, U_i))=\cap_{i=1}^m(\cup_{\tau_i'\in \sI^{K'}_{\tau_i}}X_K((\tau_i', U_i)))$$
 $$=\cup_{(\tau_1',\dots,\tau_m')\in \prod_{i=1}^m\sI^{K'}_{\tau_i}}X_{K'}((\tau_i', U_i), i=1,\dots,m).$$
 
 In particular, if for $i=1,\dots,m,$ $\tau_{i}': K'\hookrightarrow \C_{p_i}$ is an extension of $\tau_i$,
then  
$$X_{K'}((\tau_i', U_i), i=1,\dots,m)\subseteq X_K((\tau_i, U_i), i=1,\dots,m).$$
\end{rem}

\begin{rem}\label{remgalaction}For $\sigma\in \Gal(\bk/K_0)$, we have 
$$\sigma(X_K((\tau_i, U_i), i=1,\dots,m))=X_{\sigma(K)}((\tau_i\circ \sigma^{-1}, U_i), i=1,\dots,m).$$

In particular, for a basic adelic subset $A$ of $X(\bk)$ over $K$ and $\sigma\in \Gal(\bk/K)$, $\sigma(A)=A.$
\end{rem}

\medskip

\begin{exe}\label{exeaoneadelic}Assume that $\bk=\overline{\Q},$ $X:=\A^1.$ Then $X$ is defined over $\Q$. Let $\tau: \Q\hookrightarrow \C$ be the unique embedding. Let $U_1:=\{x\in \C=\A^1(\C)|\,\, |x|<1\}$ and $U_2:=\{x\in \C=\A^1(\C)|\,\, |x|>1\}$ be two disjoint open subsets of $\A^1(\C).$ Then  $$\{1/2,1/3,\dots\}\subseteq \A^1_{\Q}(\tau, U_1)\setminus \A^1_{\Q}(\tau, U_2);$$
$$\{2,3,\dots\}\subseteq \A^1_{\Q}(\tau, U_2)\setminus \A^1_{\Q}(\tau, U_1)$$ and $$\{n\pm \sqrt{n^2-1}, n=2,3,\dots\}\subseteq \A^1_{\Q}(\tau, U_1)\cap \A^1_{\Q}(\tau, U_2)=\A^1_{\Q}((\tau, U_1),(\tau,U_2)).$$

We note that $\A^1_{\Q}(\tau, U_1)$ is a basic adelic subset of  $\A^1(\overline{\Q})$.  Both  $\A^1_{\Q}(\tau, U_1)$ and $\A^1(\overline{\Q})\setminus(\A^1_{\Q}(\tau, U_1))$ are infinite. 
\end{exe}

In Example \ref{exeaoneadelic}, $\A^1_{\Q}(\tau, U_1)\cap \A^1_{\Q}(\tau, U_2)\neq\emptyset$, even when $U_1\cap U_2=\emptyset$. This phenomenon motivates the following property.

\begin{pro}\label{prostrappro}
If $X$ is irreducible and the $U_i, i=1,\dots,m$ are nonempty, then the basic adelic subset $X_K((\tau_i, U_i), i=1,\dots,m)$ is nonempty.
\end{pro}

\proof[Proof of Proposition \ref{prostrappro}]
Since $K$ is finitely generated over $\Q$, $M:=K\cap \overline{\Q}$ is a number field. For every nontrivial absolute value $|\cdot|$ on $M$, there exists a finite extension $M'$ such that there are two distinct absolute values $|\cdot|_1, |\cdot|_2$ on $M'$ which extend $|\cdot|.$
Then there exists a finite field extension $K'$ of $K$, such that there are extensions $\tau_i': K'\hookrightarrow \C_{p_i}$ of $\tau_i, i=1,\dots,m,$ such that
the absolute values $|\tau_i'(\cdot)|, i=1,\dots,m$ on $K'$ are distinct.

By Remark \ref{remfieldextad}, we have $X_{K'}((\tau_i', U_i), i=1,\dots,m)\subseteq X_K((\tau_i, U_i), i=1,\dots,m).$
Replacing $K$ by $K'$ and $\tau_i$ by $\tau_i', i=1,\dots,m$, we assume that $|\cdot |_i:=|\tau_i(\cdot)|, i=1,\dots,m$ on $K$ are distinct. 
Let $K_{p_i}$ be the closure of $\tau_i(K)$ in $\C_{p_i}.$

\begin{lem}\label{lemapproaone}\cite[Page 35, Theorem 1]{Lang1994}The image of the diagonal embedding $$x\mapsto (\tau_1(x),\dots, \tau_m(x)): K\hookrightarrow \prod_{i=1}^m K_{p_i}$$ is dense. 
\end{lem}

By Remark \ref{reminterzari}, we may assume that $X$ is smooth and affine.
Set $d:=\dim X.$ There exists a finite morphism $\pi: X_K\to \A_K^d$. 
(We still denote by $\pi$ the induced morphism $X\to \A_{\bk}^d.$)
There exists a Zariski dense open subset $V$ of $\A_K^d$ such that $\pi|_{\pi^{-1}(V)}:\pi^{-1}(V)\to V$
is an \'etale covering.

Let $\psi_i: \A_K^d(K)\hookrightarrow  \A_K^d(K_{p_i})$ be the morphism $\psi_i: (x_1,\dots, x_d)\mapsto (\tau_i(x_1),\dots, \tau_i(x_d)).$
and $\psi: \A_K^d(K)\hookrightarrow \prod_{i=1}^m \A_K^d(K_{p_i})$ be the diagonal embedding  $\psi: x\mapsto (\psi_1(x),\dots,\psi_m(x)).$
Lemma \ref{lemapproaone} shows that the image of $\psi$ is dense. So the image of 
$\psi|_{V(K)}: V(K)\hookrightarrow \prod_{i=1}^m V(K_{p_i})$ is dense.
Then the image of the diagonal map $\phi_{tot}: V(K)\hookrightarrow \prod_{\tau\in \sM_K} V(K_{p_{\tau}})$ is dense.

Now we replace $X_K$ by $\pi^{-1}(V)$.  By Remark \ref{remsurject}, we may replace $X_K$ by a Galois \'etale cover over $V$ which dominates $X_K$. Then we may assume that $\pi$ is Galois.
By \cite[Proposition 3.3.1]{Serre2008}, there exists a thin set $A\subseteq V(K)$ in the sense of Serre,  such that for every point $x\in V(K)\setminus A$, the fiber $\pi^{-1}(x)$ is integral. 

Set $W_i:=(\pi\otimes_{\bk}\C_{p_i})(U_i)\subseteq V(\C_{p_i}), i=1,\dots, m.$ They are open subsets.  After replacing $K$ by a finite extension, we may assume that 
$V(K_{p_i})\cap W_i\neq\emptyset$ for every $i=1,\dots,m$


\begin{lem}\label{lempsivsad}The image $\psi(V(K)\setminus A)$ is dense in $\prod_{i=1}^m V(K_{p_i}).$
\end{lem}
By Lemma \ref{lempsivsad}, which we prove below, there exists a point $x\in V(K)\setminus A$ such that for $i=1,\dots,m,$ $\psi_i(x)\in V(K_{p_i})\cap W_i.$
Since $\pi^{-1}(x)$ is integral, we have $\pi^{-1}(x)=\Spec(L)$ for some finite field extension $L$ of $K.$
An inclusion $L\subseteq \bk=\overline{K}$ gives a morphism $\Spec(\bk)\to \pi^{-1}(x)\subseteq X.$
This defines a point $z\in X(\bk).$

For $i=1,\dots,m,$
there exists $y_i\in \pi^{-1}(x)(\C_{p_i})\cap U_i;$ it gives a morphism $\tau_{y_i}: L\hookrightarrow \C_{p_i}$ which extends $\tau_i.$
We extend $\tau_{y_i}$ to a morphism $\overline{\tau_i}: \overline{K}\hookrightarrow \C_{p_i}.$ 
Let $\phi_{\overline{\tau_i}}: X(\overline{K})\hookrightarrow X(\C_{p_i})$ be the embedding induced by $\overline{\tau_i}$.
Then we have $\phi_{\overline{\tau_i}}(z)\in U_i, i=1,\dots,m.$ So we have $z\in X_K((\tau_i, U_i), i=1,\dots,m)$.
\endproof

\proof[Proof of Lemma \ref{lempsivsad}]
When $K$ is a number field, the following lemma is \cite[Theorem 3.5.3]{Serre2008}.
\begin{lem}\label{lemthinsetthin}
Let $V$ be a geometrically irreducible variety over $K$ and let $A$ be a thin subset of $V(K).$
For every finite subset $S_0\subseteq \sM_K^f$,
there exists a finite set $S\subseteq \sM_K^f\setminus S_0$ such that  
the image of $A$ in $\prod_{\tau \in S}V(K_{p_{\tau}})$ under the diagonal map $\phi_S$ is not dense.
\end{lem}

It follows that the image of the diagonal map $\phi_{f}:A\hookrightarrow \prod_{\tau \in \sM_K^f}V(K_{p_{\tau}})$ is nowhere dense. 
Then the image of $\phi_{tot}:A\hookrightarrow \prod_{\tau \in \sM_K}V(K_{p_{\tau}})$ is nowhere dense. 
Since $\phi_{tot}(V(K))$ is dense in
$\prod_{\tau \in \sM_K}V(K_{p_{\tau}}),$ $\phi_{tot}(V(K)\setminus A)$ is dense in
$\prod_{\tau \in \sM_K}V(K_{p_{\tau}}).$ This concludes the proof.
\endproof

\proof[Proof of Lemma \ref{lemthinsetthin}]
Let $L$ be a finite Galois extension of $K.$
There exists a subring $R$ of $K$ which is finitely generated over $\Z$ such that $K=\Frac R$.
Set $W_K:=\Spec R$ and let $W_L$ be the normalization of $W_K$ in $L.$ After shrinking $W_K$, we may assume that $W_K$ is regular.

\begin{lem}\label{lemexsp}
For every $N\geq 0$, and every $g\in R\setminus \{0\}$,
there exists a prime $p\geq N$ and an embedding $\tau: K\hookrightarrow \C_p$ in $\sI_K^f$ such that 
the absolute value $|\tau(\cdot)|$ on $K$ splits completely in $L$, and $\tau(R_{g})\subseteq \C_p^{\circ}.$
\end{lem}

In the proof of \cite[Theorem 3.5.3]{Serre2008}, Serre used the fact that,
when $K$ is a number field, there are infinitely many places of $K$
that split completely in $L.$
Once we replace this fact by Lemma \ref{lemexsp}, the same proof works.
\endproof

\proof[Proof of Lemma \ref{lemexsp}]
After replacing $R$ by $R_{g}$, we may assume that $g=1.$
Pick any $y\in W_K(\overline{\Q})$, such that the morphism $W_L\to W_K$ is \'etale at $y.$ 
Denote by $o$ the image of $y$. Then the residue field $\kappa(o)$ is a number field and $o$ induces an embedding $\iota: \Spec O_{\kappa(o),S}\hookrightarrow W_K$, where $S$ is a finite set of places of $\kappa(o)$ containing all infinite places. 

Denote by $q_1,\dots, q_s$ the pre-image of $o$ is $W_L.$ Since the extension $L$ over $K$ is Galois, the extensions $\kappa(q_i)$ over $\kappa(o), i=1,\dots,s$ are isomorphic to each other.

For every $N>0$, there exists a closed point $x\in \Spec O_{\kappa(o),S}$ which is completely splitting in the extension $\kappa(q_1)$ over $\kappa(o)$ and whose residue field $\kappa(x)$ is of characteristic $p>N$.  Then it is completely splitting in the extension $\kappa(q_i)$ over $\kappa(o)$ for every $i=1,\dots,s.$  
Let $m$ be the maximal ideal of $R$ corresponding to $\iota(x)$. The pre-image of $\iota(x)$ in $W_L$ has exactly  $[L:K]$ points. Now we only need to show that there exists an embedding $\tau: R\hookrightarrow \C_p^{\circ},$ such that $m=\tau^{-1}(\C_p^{\circ\circ}).$

For every $P\in R\setminus \{0\},$ denote by $Z_P$ the set $\{z\in W_K(\C_p)|\,\, P(z)=0\}.$ It is a nowhere dense closed subset of $W_K(\C_p).$ Observe that the topology of $W_K(\C_p)$ can be defined by a complete metric. Since $R\setminus \{0\}$ is countable, by Baire category theorem, $W_K(\C_p)\setminus (\cup_{R\setminus \{0\}}Z_P)$ is dense in $W_K(\C_p).$
Let $P_1,\dots, P_n$ be a set of generator of $m.$ Then $B:=\{z\in W_K(\C_p)|\,\, |P_i(z)|<1\}$ is an open subset of $W_K(\C_p).$
Pick any inclusion $\overline{\Q}\subseteq \C_p$. Using this inclusion, we may view $y$ as a point in $B\subseteq W_K(\C_p).$ So $B$ is not empty. Pick any point $z\in B\setminus ( (\cup_{R\setminus \{0\}}Z_P)).$ Then the inclusion $\tau_z: R\hookrightarrow \C_p$ defined by $z$ is what we need. 
\endproof

\subsubsection{General adelic subsets}
\begin{defi}A \emph{general adelic subset of $X(\bk)$ (over $K$)} is a subset taking form $\pi(B)$ where $\pi: Y\to X$ is a flat morphism (defined over $K$) and  $B$  is a basic  adelic subset of $Y$ (over $K$).
\end{defi}
In particular, every basic adelic subset over $K$ is a general adelic subset.

%
%
%
%
%
%
%
%
%

\begin{pro}\label{progenadelicsub}The following properties for general adelic subsets hold:
\begin{points}
\item Let $f: Y\to X$ be a morphism defined over $K$. Let $A\subseteq X(\bk)$ be a general adelic subset over $K$, then $f^{-1}(A)\subseteq Y(\bk)$ is a general adelic subset over $K$.
\item For every general adelic subset $A\subseteq X(\bk)$ over $K$, there is a flat morphism $\pi: Y\to X$ defined over $K_0$ and a basic adelic subset $B$ over $K$, such that 
$A=\pi(B).$ 
\item For every $\sigma\in \Gal(\bk/K_0)$ and every general adelic subset $A\subseteq X(\bk)$ over $K$, $\sigma(A)$ is a general adelic subset of $X(\bk)$ over $\sigma(K).$ In particular, if $\sigma\in \Gal(\bk/K)$, $\sigma(A)=A.$
\item For general adelic subsets $A_1,A_2$ of $X(\bk)$ over $K$, $A_1\cup A_2$ is a general adelic subset over $K.$
\item Let $A$ be a general adelic subset of $X(\bk)$ over $K.$ Let $K'$ be a finite extension of $K$. 
Then $A$ is a general adelic subset of $X(\bk)$ over $K'.$
\item Let $f: Y\to X$ be a flat morphism defined over $K$. Let $A\subseteq Y(\bk)$ be a general adelic subset over $K$, then $f(A)\subseteq X(\bk)$ is a general adelic subset over $K$.
\item For general adelic subsets $A_1, A_2\subseteq X(\bk)$ over $K$, their intersection $A_1\cap A_2$ is a general adelic subset over $K$.
\item Assume that $X$ is irreducible. Let $A_1,A_2$ be two nonempty general adelic subsets of $X(\bk).$ Then $A_1\cap A_2\neq \emptyset.$
\end{points}
\end{pro}

\proof[Proof of Proposition \ref{progenadelicsub}]

(i). 
Write $A=\pi(B)$, where  $\pi: Z\to X$ is a flat morphism defined over $K$ and $B$ is a basic adelic subset of $Z(\bk)$ over $K$.
Consider the fiber product $Z_Y:=Y\times_X Z$, and the projections $\pi_Y: Z_Y\to Y$ and $\pi_Z: Z_Y\to Z.$ We note that $\pi_Z$ and $\pi_Y$ are defined over $K$, $\pi_Y$ is flat and
$f\circ \pi_Y=\pi\circ \pi_Z.$
Then 
$$f^{-1}(A)=f^{-1}(\pi(B))=\{y\in Y(\bk)|\,\, \text{there exists } b\in B \text{ such that } f(y)=\pi(b)\}$$
$$=\pi_Y(\{(y,b)\in Y\times_X Z=Z_Y|\,\, b\in  B\})=\pi_Y(\pi_Z^{-1}(B)).$$
By Remark \ref{remsurject},  $\pi_Z^{-1}(B)$ is a basic adelic subset over $K.$
This concludes the proof.

\medskip

(ii). Write $A=\pi_1(C)$ where $\pi_1: Z\to X$ is a flat morphism over $K$ and $C$ is a basic adelic subset of $Z(\bk)$ over $K$.
Let $L$ be the normal closure of $K$ over $K_0.$
For every $\sigma\in \Gal(L/K_0)$, consider the Galois conjugate of $\pi_1$ by $\sigma$, $$\pi_1^{\sigma}: Z^{\sigma}\to X^{\sigma}=X.$$
Then the fiber product $Y:=\prod_X^{\sigma\in \Gal(L/K_0)}Z^{\sigma}$ and the natural morphism $\pi:\prod_XZ^{\sigma}\to X$ are defined over $K_0$. Moreover $\pi$ is flat.
Consider the projection $$\psi: Y:=\prod_X^{\sigma\in \Gal(K/K_0)}Z^{\sigma}\to Z^{\id}=Z.$$ By (i), since $\phi$ is defined over $K,$
$B:=\psi^{-1}(C)$ is a basic adelic subset of $Y(\bk)$ over $K.$
Since $\pi_1\circ \psi=\pi$, $A=\pi(\psi^{-1}(C))=\pi(B).$

\medskip

(iii). 
By (ii), write $A=\pi(B)$, where  $\pi: Y\to X$ is a flat morphism defined over $K_0$ and $B$ is a basic adelic subset of $Y(\bk)$ over $K$.
Then we have  $\sigma(A)=\sigma(\pi(B))=\pi(\sigma(B)).$ By Remark \ref{remgalaction}, $\sigma(B)$ is a basic adelic subset of $Y(\bk)$ over $\sigma(K)$.
If $\sigma\in \Gal(\bk/K)$, we get $\sigma(A)=\pi(\sigma(B))=\pi(B)=A.$

\medskip

(iv). Write $A_i=\pi_i(B_i), i=1,2,$ where $\pi_i: Y_i\to X$ is a flat morphism defined over $K$ and $B_i$ is a basic adelic subset of $Y_i(\bk)$ over $K$.
Set $Y:=Y_1 \sqcup Y_2$ and $\pi:=\pi_1\sqcup \pi_2: Y=Y_1\sqcup Y_2\to X.$ Then $\pi$ is flat. Since $A_1\sqcup A_2$ is a basic adelic subset of $(Y_1\sqcup Y_2)(\bk)$ over $K$ and $A_1\cup A_2=\pi(A_1\sqcup A_2)$, $A_1\cup A_2$ is a general adelic subset of $X(\bk)$ over $K.$

\medskip

(v).
Write $A=\pi(B)$, where  $\pi: Y\to X$ is a flat morphism defined over $K$ and $B$ is a general adelic subset of $Y(\bk)$ over $K$.
By Remark \ref{remfieldextad}, $B=\cup_{i=1}^m B_i$ where $B_i, i=1,\dots,m$ are basic adelic subsets of $Y(\bk)$ over $K'$.
Set $A_i:=\pi(B_i), i=1,\dots,m$. Then $A_i, i=1,\dots,m$ are general adelic subsets of $Y(\bk)$ over $K'$.
Then $$A=\pi(B)=\pi(\cup_{i=1}^m B_i)=\cup_{i=1}^m\pi(B_i)=\cup_{i=1}^mA_i.$$
By (iv), $A$ is a general adelic subset of $X(\bk)$ over $K'$.

\medskip

(vi). This follows from the definition of the the general adelic subsets.

\medskip

(vii). Write $A_i=\pi_i(B_{i}), i=1,2$, where  $\pi: Y_i\to X$ is a flat morphism defined over $K$ and $B_{i}$ is a basic adelic subset of $Y_i(\bk)$ over $K$.

Consider the fiber product $Y:= Y_1\times_X Y_2$. The natural morphism $\pi: Y\to X$ is flat and defined over $K.$
For every $i=1,2$, denote by $\psi_i: Y\to Y_i$ the $i$-th projection. They are flat and defined over $K$. We have $\pi=\pi_i\circ\psi_i$ for $i=1,2.$
Then $$A_1\cap A_2=(\pi_1(B_1))\cap (\pi_2(B_2))=\pi(B_1\times_{X}B_2)=\pi(\psi_1^{-1}(B_1)\cap \psi_2^{-1}(B_2))$$
By (i), every $\psi_i^{-1}(B_{i}), i=1,2$ is a basic adelic subset over $K$. So $A_1\cap A_2$ is a general adelic subset over $K$.

\medskip

(viii). 
By  the proof of (vi), we may assume that $A_i=\pi(B_i), i=1,2$, where  $\pi: Z\to X$ is a flat morphism defined over $K$ and $B_i, i=1,2$ are nonempty basic adelic subsets of $Z(\bk)$.
Denote by $Z_1,\dots,Z_s$ the irreducible components of $Z.$

If there is $j=1,\dots, s$ such that both $Z_j\cap B_1$ and $Z_j\cap B_2$ are nonempty, by Proposition \ref{prostrappro}, $B_1\cap B_2\neq \emptyset.$ Since $\pi(B_1\cap B_2)\subseteq \pi(B_1)\cap \pi(B_2)=A_1\cap  A_2,$ $A_1\cap A_2$ is nonempty.

Now we may assume that for every $j=1,\dots,s$, one of $Z_j\cap B_1$, $Z_j\cap B_2$ is empty.
Since both $B_1,B_2$ are nonempty, we may assume that  $B_1\cap Z_1\neq\emptyset$ and $B_2\cap Z_2\neq\emptyset.$
For $i=1,2$, pick a nonempty Zariski open subset $U_i$ of $Z_i$ such that $U_i\cap Z_j=\emptyset$ for $j\neq i.$

Since every $U_i, i=1,2$ is a basic adelic subset, $U_i\cap B_i$ is a nonempty basic adelic subset.
After replacing $Z_i, i=1,2$ by $U_i, i=1,2$, $B_i, i=1,2$ by $B_i\cap U_i, i=1,2$ and $Z$ by $U_1\sqcup U_2$, we may assume that $Z=Z_1\sqcup Z_2$, $B_1\subseteq Z_1$ and $B_2\subseteq Z_2$.
Set $\pi_i:=\pi|_{Z_i}: Z_i\to X, i=1,2.$ Then $\pi_i, i=1,2$ are flat and $A_i=\pi_i(B_i), i=1,2.$

Consider the fiber product $W:=Z_1\times_X Z_2$. The natural morphism $\psi: W\to X$ and the projections $\psi_i: W\to Z_i, i=1,2$ are flat. 
Pick $U$ a nonempty, irreducible, Zariski open subset of $W$.  Then every $\psi_i|_U, i=1,2$ is flat. So $\psi_i(U)$ is open and nonempty in $Z_i.$
By Proposition \ref{prostrappro}, $\psi_i(U)\cap B_i\neq \emptyset.$ By Remark \ref{remsurject}, $\psi_i^{-1}(B_i)\cap U, i=1,2$ are nonempty basic adelic subsets of $U$. By Proposition \ref{prostrappro} again,
$\psi_1^{-1}(B_1)\cap \psi_2^{-1}(B_2)\cap U\neq\emptyset.$
Since $$\psi(\psi_1^{-1}(B_1)\cap \psi_2^{-1}(B_2)\cap U)\subseteq \psi(\psi_1^{-1}(B_1)\cap \psi_2^{-1}(B_2))$$
$$\subseteq \psi(\psi_1^{-1}(B_1))\cap \psi(\psi_2^{-1}(B_2))=\pi_1(B_1)\cap \pi_2(B_2)=A_1\cap A_2,$$
$A_1\cap A_2$ is nonempty.
\endproof

\subsubsection{The adelic topology}
By (vii) of Proposition \ref{progenadelicsub}, the set of general adelic subsets is closed under finite intersection. So we may define the adelic topology on $X(\bk)$ as follows:
\begin{prodefi}
A subset $S\subseteq X(\bk)$ is an \emph{adelic open subset}, if it is a union of general adelic subsets. 
Adelic open subsets form a topology, called the \emph{adelic topology}. General adelic subsets forms a basis for this topology.
\end{prodefi}
By (v) of  Proposition \ref{progenadelicsub},  the adelic topology does not depend on the choice of the base field $K_0$.

\begin{rem}
The adelic topology can be described as the minimal topology defined on $X(\bk)$ for every 
variety $X$ over $\bk$, which contains all basic adelic subsets and makes all flat morphisms open.
\end{rem}

\begin{pro}\label{proadelic}The adelic topology has the following basic properties.
\begin{points}
\item It is stronger than the Zariski topology.
\item It is ${\mathsf{T}}_1$ i.e. for every pair of distinct points $x, y\in X(\bk)$ there are adelic open subsets $U,V$ of $X(\bk)$ such that 
$x\in U, y\not\in U$ and $y\in V, x\not\in V.$
\item Morphisms between algebraic varieties over $\bk$ are continuous.
\item Flat morphisms are open w.r.t. the adelic topology.
\item The irreducible components of $X(\bk)$ in the Zariski topology are the irreducible components of $X(\bk)$ in the  adelic topology.
\item Let $K$ be any subfield of $\bk$ which is finitely generated over $\Q$ such that $X$ is defined over $K$ and $\overline{K}=\bk$. Then the action 
$$\rho:\Gal(\bk/K)\times X(\bk)\to X(\bk)$$ sending $(\sigma,x)$ to $\sigma(x)$ is continuous w.r.t. the adelic topology.
Here the topology of $\Gal(\bk/K)$ is the topology of profinite groups.
\end{points}
\end{pro}

\begin{rem}
Example \ref{exeaoneadelic} shows that the adelic topology on $\A^1(\overline{\Q})$ is strictly stronger than the Zariski topology.  
\end{rem}

\begin{rem}
When $X$ is irreducible, (v) of Proposition \ref{proadelic} implies that the intersection of finitely many nonempty adelic open subsets of $X(\bk)$ is nonempty.
This also shows that, in general, the adelic topology is not Hausdorff.
\end{rem}

\proof[Proof of Proposition \ref{proadelic}]
The properties (i)-(v) easily follows from Proposition \ref{progenadelicsub}.
We only need to prove (vi). Let $U$ be a general adelic subset of $X(\bk)$ over a finite Galois extension $L$ of $K.$
We only need to show that $\rho^{-1}(U)$ is open in $\Gal(\bk/K)\times X(\bk)$.
Let $(\sigma,x)$ be a point in $\rho^{-1}(U)$. 
Since $L$ is Galois over $K$, $\sigma(L)=L$. By (iii) of Proposition \ref{progenadelicsub}, $\sigma^{-1}(U)$ is a general adelic subset of $X(\bk)$ over $L$ and $x\in \sigma^{-1}(U).$ 
By (iii) of Proposition \ref{progenadelicsub} again, for every $\tau\in \Gal(\bk/L)$, $\tau(\sigma^{-1}(U))=\sigma^{-1}(U).$
Then we have $(\sigma,x)\in (\Gal(\bk/L)\sigma)\times (\sigma^{-1}(U))$ and $(\Gal(\bk/L)\sigma)\times (\sigma^{-1}(U))\subseteq \rho^{-1}(U).$
This concludes the proof. 
\endproof

\begin{exe}
Assume that $\bk=\overline{\Q},$ $X:=\A^1.$ 
The diagonal embedding $\A^1(\bk)\hookrightarrow \prod_{\tau\in \sI_{\bk}} \A^1(\C_{p_{\tau}})$ induces a topology on $\A^1(\bk)$, which is called the product topology.
These two topologies are different on $\A^1(\overline{\Q})$. Here, the product topology is Hausdorff, but the adelic topology is not.
\end{exe}

\begin{rem}
We defined the adelic topology by using all valuations of $\bk$ given by any embeddings in $\C_p$ where $p$ is a prime or $\infty.$ 
If we reduce or enlarge the set of valuations and consider, 
for example, only the archimedean ones or the nonarchimedean ones, or all possible valuations which are non trivial on $\overline{\Q}$ etc, we may get another topology which satisfies similar properties.
However, the current range is sufficient for this paper.
\end{rem}

\begin{rem}\label{remhowuse}
In this paper, we usually only use the adelic topology in the following way.
Assume that $X$ is irreducible.
Let $S$ be an algebraic structrue on $X$ as in Section \ref{subsectionreal}.
Let $P_1,\dots , P_n$ be finitely many algebraic properties for points in $X$, which depends on $X,S$. 

There exists a subfield $K$ of $\bk$ which is finitely generated over $\Q$ such that $\overline{K}=\bk$ and $X,S$ are defined over $K.$
Assume we know that, for every $i=1,\dots,n$, there exists an embedding $\tau_i:\bk\hookrightarrow \C_{p_i}$ and a nonempty open subset $U_i$ of $X(\C_{p_i})$ such that 
$P_i$ is satisfied for all points in $U_i;$  we get automatically that all points in the nonempty basic adelic subset $X_K((\tau_1|_K,U_1),\dots, (\tau_n|_K,U_n))$ satisfy  $P_1,\dots,P_n.$
\end{rem}

\subsection{Invariant polydisk and the dynamical Mordell-Lang conjecture}
Assume that $X$ is irreducible of dimension $d$.
Let $K$ be a finitely generated field extension over $\Q$ such that
$\overline{K}=\bk$, and $X$ is defined over $K$.
Let $f: X\dashrightarrow X$ be a dominant rational self-map of $X$ defined over $K$.

There exists a variety $X_K$ over $K$ such that $X= X_K \times_{\Spec(K)} \Spec(\bk)$ and an endomorphism $f_K: X_K\dashrightarrow X_K$ such that 
$f=f_K\times_{\Spec (K)} \id.$

In this section, we use the strategy in \cite{Amerik} and \cite{Bell2010} to 
show that the dynamical  Mordell-Lang conjecture holds for an adelic general point in $X(\overline{K})$.


\begin{pro}\label{proinvpoly}There exists $m\geq 1$, a prime $p\geq 3$, an embedding $\iota:K\hookrightarrow \C_p$, and an open subset $V\subseteq X_K(\C_p)$, such that 
\begin{points}
\item $V$ is analytically isomorphic to $(\C_p^{\circ})^d$;
\item  $V$ is $f^m$-invariant, and $f^m|_V=\id \mod p$;
\item for every $x\in V$, its $f$-orbit is well defined.
\end{points}
Moreover, there exists an analytic action $\Phi: \C_p^{\circ}\times V\to V$ of
$(\C_p^{\circ},+)$ on $V$ such that for every $n\in \Z_{\geq 0}$,  $\Phi(n,\cdot)=f^{mn}|_{V}(\cdot).$
In particular, for every $x\in X_K(\iota,V)$, the Zariski closure of the orbit $O_{f^m}(x)$ in $X$ is irreducible. 
\end{pro}

\proof[Proof of Proposition \ref{proinvpoly}]
There exists a subring $R$ of $K$ that is finitely generated over $\Z$
with the property that $\Frac R=K$.

Pick a model $\pi : X_R \to \Spec(R)$ which is projective over $\Spec(R)$ and whose generic fiber is $X_K$. 
Then $f$ extends to a rational self-map $f_R: X_R\dashrightarrow X_R.$ Denote by $B_R$ the union of indeterminacy locus of $f_R$, the non-\'etale locus of $f_R$, and the non-smooth locus of $X_R$.

\begin{lem}\label{lem:good-model-lem}
There exists a nonempty, affine, open subset $U$ of $\Spec(R)$ such that
\begin{enumerate}
\item $U$ is of finite type over $\Spec(\Z)$;
\item for every point $y \in U$, the fiber $X_y$ is geometrically irreducible  and  
$\dim_{K(y)} X_y = \dim_K X_K$, where $K(y)$ is the residue field at $y$;
\item for every  $y \in U$, the fiber $X_y$ is not contained in
$B_{R}$.
\end{enumerate}
\end{lem}

\proof[Proof of Lemma \ref{lem:good-model-lem}]
To prove the lemma, we shall use the following fact:
For any integral affine scheme $\Spec(A)$ of finite type over $\Spec(\Z)$ and any nonempty open
subset $V_1$ of $\Spec(A)$, there exists an affine open subset $V_2$ of $V_1$ which is of finite type over $\Spec(\Z)$.
Indeed, we may pick any nonzero element $g\in I$ where $I$ is the ideal of $A$ that defines the closed subset $\Spec(A)\setminus V$ and set 
$U :=\Spec(A)\setminus \{g = 0\}$. Then $U = \Spec(A[1/ g])$ is of finite type over $\Spec(\Z)$.

Since $X_K$ is geometrically irreducible, \cite[Proposition 9.7.8]{EGA4} gives an affine open subset $V$ of $\Spec(R)$
 such that $X_y$ is geometrically irreducible  for every $y \in V$. We may suppose that $V$
  is of finite type over $\Spec(\Z)$. By generic flatness (see \cite[Thm. 6.9.1]{EGA4}), 
  we may change $V$ in a smaller subset and suppose
  that the restriction of $\pi$  to $\pi^{-1}(V)$ is flat. Then, the fiber $X_y$ is geometrically irreducible  and of dimension
  $\dim_{K(y)} X_y = \dim_K X_K$ for every point $y \in V$.

Denote by $B_{K}$ the union of the indeterminacy locus, the non-\'etale locus of $f$ in $X_K$, the singular locus of $X_K$ and $Z_K$.
Observe that $B_{K}$ is exactly the generic fiber of $\pi_{\vert B_{R}} \colon B_{R}\to  \Spec(R)$.
By generic flatness, there exists a nonempty, affine, open subset $U$ of $V$ such that the restriction of $\pi$ to every irreducible component of
 $\pi_{\vert B_{R}}^{-1} (U)$ is flat. 
Then for $y \in U$, the fiber $X_y$ is not contained in $B_{R}$.
Then, we shrink $U$ to suppose that $U$ is of finite type over $\Spec(\Z)$. Since
\[
\dim_{K(y)}(B_{R}\cap X_y) = \dim_K(B_{K}) < \dim_K X_K = \dim_{K(y)} X_y
\]
 for every  $y \in V$, the fiber $X_y$ is not contained in $B_{R}$. 
 \endproof

By Lemma~\ref{lem:good-model-lem}, we may replace $\Spec(R)$ by $U$ and assume that
\begin{itemize}
\item for every $y \in \Spec(R)$, the fiber $X_y$ is geometrically irreducible ;

\item for every $y \in \Spec(R)$, the fiber $X_y$ is not contained in $B_{R}$.
\end{itemize}

Recall the following Lemma (\cite{Lech1953,Bell2006} and \cite[Proposition 2.5.3.1]{Bell2016}).
\begin{lem}\label{lem:Lech-Bell}
Let $L$ be a finitely generated extension of $\Q$ and $B$ be a finite subset of $L$. The set of primes
$p$ for which there exists an embedding of $L$ into $\Q_p$ that maps $B$ into $\Z_p$ has positive density\footnote{By {\bf{positive density}}, we mean that the proportion of primes $p$ among the first $N$ primes that satisfy the statement
is bounded from below by a positive number if $N$ is large enough. }
among the set of all primes. 
\end{lem}
Since $R$ is integral and finitely generated over $\Z$, by Lemma~\ref{lem:Lech-Bell} 
 there exists infinitely many primes $p \geq 3$ such that $R$ can be embedded into $\Z_p\subseteq \C_p^{\circ}$. This induces an embedding
$\Spec(\Z_p) \to \Spec(R)$. Set $X_{\C_p} := X_R \times_{\Spec(R)}\Spec(\C_p)$, and $f_{\C_p}:=f_R\times_{\Spec(R)} \id.$
All fibers $X_y$, for $y\in \Spec(R)$, are geometrically irreducible  and of dimension $d$;
hence, the special fiber $X_{\overline{\F_p}}$ of $X_{\C_p^{\circ}} \to  \Spec(\C_p^{\circ})$ is geometrically irreducible .
Denote by $B_{\C_p}$ the union of indeterminacy locus, the non\'etale of $f_{\C_p}$, the singular locus of $X_{\C_p}$.
 Since $B_{\C_p} \subset B_{R} \cap  X_{\overline{\F_p}}$, the fiber $X_{\overline{\F_p}}$ is not contained in $B_{\C_p}$. 
We note that $X_{\overline{\F_p}}$ and  $f|_{X_{\overline{\F_p}}\setminus B_{\C_p}}$ are indeed defined over $\F_p$.

Apply \cite[Corollary 2]{Amerik} to the rational map $f|_{X_{\overline{\F_p}}\setminus B_{\C_p}}: X_{\overline{\F_p}}\setminus B_{\C_p}\dashrightarrow X_{\overline{\F_p}}\setminus B_{\C_p}$
 there exists a periodic point $x\in X_{\overline{\F_p}}(\overline{\F_p})\setminus B_{\C_p}$ whose orbit under $f_{\C_p}$ is contained in $X_{\overline{\F_p}}(\overline{\F_p})\setminus B_{\C_p}.$
 Observe that $x$ is a regular closed point in $X_{\C_p^{\circ}}.$ There exists $m\geq 1$ such that 
 $f|_{X_{\overline{\F_p}}}^m(x)=x$ and 
 $(df|_{X_{\overline{\F_p}}}^m)_x=\id$.
Let $U$ be the open subset of $X(\C_p)$ consisting of points whose specialization is $x.$ Then we have $U\simeq (\C_p^{\circ\circ})^{d}.$
Then we have $f^m(U)=U$ and the orbit of the points in $U$ are well-defined.
 The restriction of $f^m$ on $U$ is an analytic automorphism taking form 
 $f^m|_U: (x_1,\dots,x_d)\mapsto (F_1,\dots, F_d)$
 where $F_n=\sum_{I}a_I^nx^I, n=1,\dots d$ are analytic functions on $U$ with
 $a_I^n\in \C_p^{\circ\circ}.$  
 Moreover, since $f,X$ are defined over $\Z_p$, in a suitable coordinate, all $a_l^n$ are contained in a finite extension $K_p$ over $\Z_p.$
 Since  $(df|_{X_{\overline{\F_p}}}^m)_x=\id$, we  have $f^m|_U=\id \mod \C_p^{\circ\circ}.$
 Then there exists $l\in \Q^{+}$, such that $f^m|_U=\id \mod p^{2/l}.$
  Set $V:=\{(x_1,\dots,x_d)\in U|\,\, |x_i|\leq p^{-1/l}\}\simeq (\C^{\circ}_p)^d.$
 Then $V$ is invariant by $f^m$ and  
 $f^m|_{V}=\id \mod p^{1/l}.$
 After replacing $m$ by some multiple of $m$, we get
$f^m|_{V}=\id \mod p.$

By \cite[Theorem 1]{Poonen2014}, 
there exists an analytic action $\Phi: \C_p^{\circ}\times V\to V$ of
$(\C_p^{\circ},+)$ on $V$ such that for every $n\geq \Z_{\geq 0}$,  $\Phi(n,\cdot)=f^{mn}|_{V}(\cdot).$ 
For every $x\in V$,  denote by $\Phi_x:\C_p^{\circ}\to V$ the analytic morphism $t\to \Phi(t,x).$
Denote $Z$ the Zariski closure of $O_{f^m}(x)$ in $X.$ Let $Z_1,\dots,Z_s$ be the irreducible component of $Z$. 
We have $\C_p^{\circ}=\cup_{i=1}^s \Phi_x^{-1}(Z_i).$ Since $\C_p^{\circ}$ is irreducible in the Zariski topology, we have $\Phi_x^{-1}(Z_i)=\C_p$ for some $i=1,\dots,s.$
It follows that $O_{f^m}(x)\subseteq \Phi_x(\C_p^{\circ})\subseteq Z_1.$ It follows that $Z=Z_1$ is irreducible.
 We showed that for every $x\in V\subseteq X(\C_p)$, the Zariski closure of the orbit $O_{f^m}(x)$ in $X(\C_p)$ is irreducible. 
By Remark \ref{remhowuse}, for every $x\in X_K(\iota,V)$, the Zariski closure of the orbit $O_{f^m}(x)$ in $X$ is irreducible.
This concludes the proof.
\endproof

\begin{pro}\label{prodmlad}There exists a prime $p\geq 3$, an embedding $\iota:K\hookrightarrow \C_p$, and an open subset $V\simeq (\C_p^{\circ})^d$ of $X_K(\C_p)$,  such that for every proper subvariety 
$Z$ of $X$ and point $x\in X_K(\iota,V)$, the orbit of $x$ is well-defined and the set $\{n\geq 0|\,\, f^n(x)\in Z\}$ is a finite union of arithmetic progressions.
In particular,  if the orbit of $x$ is Zariski dense in $X$, then $\{n\geq 0|\,\, f^n(x)\in Z\}$ is finite.
\end{pro}
 
 \proof[Proof of Proposition \ref{prodmlad}]
 By Proposition \ref{proinvpoly},
 there exists $m\geq 1$, a prime $p\geq 3$ and an embedding $\iota:K\hookrightarrow \C_p$, such that on $X_K(\C_p)$ there exists an open subset $V\simeq (\C_p^{\circ})^d$,  such that $V$ is invariant by $f^m$, the orbit of the points in $V$ are well-defined  and there exists an analytic action $\Phi: \C_p^{\circ}\times V\to V$ of
$(\C_p^{\circ},+)$ on $V$ such that for every $n\in \Z_{\geq 0}$,  $\Phi(n,\cdot)=f^{mn}|_{V}(\cdot).$

Let $x$ be a point in $X_K(\iota,V)$, there exists $\overline{\iota}\in \sI_{\iota}$ such that $\phi_{\overline{\iota}}(x)\in V.$
Using $\overline{\iota}$ to view $\bk$ as a subfield of $\C_p$ and identify $x$ with $\phi_{\overline{\iota}}(x)\in V.$
Set $Z_j:=(f|_{V}^j)^{-1}(Z), j=0,\dots,m-1.$ Set $g:=f^m|_V$, we only need to show that for every $j=0,\dots,m-1$, the set 
$S_j:=\{n\geq 0|\,\, g^n(x)\in Z_j\}$ is a finite union of arithmetic progressions.
Observe that $S_j:=\{n\geq 0|\,\, g^n(x)\in Z_j\}\subseteq T_j:=\{t\in \C_p^{\circ}|\,\,\Phi(t,x)\in Z_j \}, j=0,\dots,m-1,$ 
and  $T_j$ is a Zariski closed subset of the disk $\C_p^{\circ}$.
 If $S_j$ is infinite, then $T_j$ is Zariski dense in $\C_p^{\circ}$. It follows that $S_j=\Z_{\geq 0}.$
So for every $j=0,\dots,m-1$, $S_j$ is either finite or $\Z_{\geq 0}.$

Assume that the orbit $O_f(x)$ of $x$ is Zariski dense in $X$. If $\{n\geq 0|\,\, f^n(x)\in Z\}$ is not finite, there exists $a\geq 0, b\geq 1$ such that 
$f^{a+bn}(x)\in Z$ for $n\geq 0.$ It follows that 
$$O_f(x)\subseteq \{x,\dots, f^{a-1}(x)\}\cup Z\cup\dots\cup f^{b-1}(Z)$$ which is not Zariski dense. This contradicts to our assumption.
Then we conclude the proof.
 \endproof

\subsection{The adelic version of the Zariski dense orbit conjecture}
Let $X$ be an irreducible variety defined over $\bk$
and let $f:X\dashrightarrow X$ be a dominant rational self-map.

\begin{rem}\label{remcurve}
Proposition \ref{proinvpoly} shows that when $X$ is a curve, $(X,f)$ satisfies the SAZD-property if and only if $f$ is not of finite order. In particular, $(X,f)$ satisfies the AZD-property
\end{rem}

\begin{pro}\label{proadelicpairbasicp}The following statements are equivalents:
\begin{points}
\item $(X,f)$ satisfies the AZD-property (resp. SAZD-property);
\item  there exists $m\geq 1$, such that $(X,f^m)$ satisfies the AZD-property (resp. SAZD-property);
\item there exists a pair $(Y,g)$ which is birational to the pair $(X,f)$, and $(Y,f)$ satisfies the AZD-property (resp. SAZD-property).
\end{points}
\end{pro}

\proof[Proof of Proposition \ref{proadelicpairbasicp}]
We only prove it for AZD-property. For SAZD-property, the proof is similar.
It is clear that (i) implies (ii) and (iii).
Lemma \ref{leminvratfunite} and Proposition \ref{prodmlad} show that (ii) implies (i).
Let us show that (iii) implies (i).
Let $\pi: Y\dashrightarrow X$ be a birational map satisfying $\pi\circ g=f\circ \pi.$
If there exists $H\in \bk(Y)^g\setminus \bk$, then $(\pi^{-1})^*H\in \bk(X)^f\setminus \bk$.

Now assume that there exists a nonempty adelic open subset $A$ of
$Y(\bk)$ such that for every point $x\in A$, 
$O_{g}(x)$ is well-defined and Zariski dense in $Y$.
Let $U$ be a dense Zariski open subset of $Y$ such that $\pi|_U: U\to X$ is an
isomorphism to its image.  Applying Proposition \ref{prodmlad} to
$g|_U: U\dashrightarrow U$, we see that there exists a nonempty adelic
open subset $A_1$ of $U(\bk)$ such that for every point $x\in A_1$, its
orbit under $g|_U$ is well-defined. Then for every $x\in A\cap A_1$,
the orbit of $x$ under $g$ is well-defined, contained in $U$ and is
Zariski dense in $Y$.  By Proposition \ref{prostrappro}, $A\cap A_1$ is a
nonempty adelic open subset of $Y(\bk)$ and $\pi(A\cap A_1)$ is a
nonempty adelic open subset of $X(\bk)$.  Then for every
$x\in \pi(A\cap A_1)$, the orbit of $x$ under $f$ is well-defined,
contained in $\pi(U)$, and Zariski dense in $X$, which concludes the
proof.  \endproof

\begin{lem}\label{lemdesfiniteadelic}Let $X'$ be an irreducible variety over $\bk$. Let $f':X'\dashrightarrow X'$ be a rational dominant self-map.
Let $\pi:X'\dashrightarrow X$ be a generically finite dominant rational map such that $\pi\circ f'=f\circ \pi.$
Then $(X',f')$ satisfies the AZD-property (resp. SAZD-property) if and only if $(X,f)$ satisfies the AZD-property (resp. SAZD-property).
\end{lem}
\proof[Proof of Lemma \ref{lemdesfiniteadelic}]
We only prove it for AZD-property. For SAZD-property, the proof is similar.
After shrinking $X'$ we may assume that $\pi$ is well-defined, quasi-finite and \'etale.

We first assume that $(X',f')$ satisfies the AZD-property.
If $\bk(X')^{f'}\neq \bk$, we conclude the proof by (ii) of Lemma \ref{leminvratfunite}.
Now we may assume that there exists a nonempty adelic open subset $A$ of $X'(\bk)$ such that 
the orbit of every 
point $x\in A$ is well-defined, contained in $\pi^{-1}(X\setminus I(f))$, and Zariski dense in $X'.$ 
Then for every $x$ in the nonempty adelic open subset $\pi(A)\subseteq X(\bk)$, the orbit  $O_{f}(x)$ is Zariski dense in $X.$ 

\medskip

Assume $(X,f)$ satisfies the AZD-property.
By by (ii) of Lemma \ref{leminvratfunite}, we may assume that $\bk(X)^f=\bk.$
Then there exists a nonempty adelic open subset $A$ of $X(\bk)$ such that the orbit of every point
$x\in A$ is well-defined, contained in $(X\setminus \overline{I(f')})\cap (X\setminus I(f))$, and Zariski dense in $X.$ 
Then for every $x$ in the nonempty adelic open subset $\pi^{-1}(A)\subseteq X'(\bk)$, the orbit  $O_{f'}(p)$ is Zariski dense in $X'.$ 
\endproof

The following result shows that the adelic version of the Zariski dense orbit conjecture implies the strong form of  the Zariski dense orbit conjecture.
\begin{cor}\label{coradelicstrong}
Let $\bk'$ be an algebraically closed field extension over $\bk$. 
If the pair $(X,f)$ satisfies the AZD-property, then $(X_{\bk'},f_{\bk'})$ satisfies the strong ZD-property.
\end{cor}
\proof[Proof of Corollary \ref{coradelicstrong}]
Let $U'$ be a nonempty Zariski open of $X_{\bk'}$. 
Set $V:=\cup_{\sigma\in \Gal(\bk'/\bk)}\sigma(U').$
Then there exists a nonempty Zariski open $U$ of $X$ such that $V=U\times_{\Spec\bk}\Spec \bk'.$
Denote by $\phi: X(\bk)\hookrightarrow X_{\bk'}(\bk')$ the natural embedding. Observe that for every $y\in  X(\bk)$, $\phi(y)$ is invariant under the action of $\Gal(\bk'/\bk)$.

If $\bk(X)^f\neq \bk$, then $\bk'(X_{\bk'})^{f_{\bk'}}\neq \bk'.$ We assume that $\bk(X)^f= \bk.$ By Proposition \ref{proadelicpairbasicp} for the pair $(U,f|_U)$, there exists a nonempty adelic open subset $A$ of $U(\bk)$ for every $x\in A$, the orbit of $x$ under $f$ is well-defined, contained in $U$, and Zariski dense in $X$. 
Then $\phi(x)\in V$, its orbit under $f_{\bk'}$ is well-defined, contained in $V$ and is Zariski dense in $X_{\bk'}$. 
For every $n\geq 0$, there exists $\sigma\in \Gal(\bk'/\bk)$ such that  $\phi(f_{\bk}^n(x))=f_{\bk'}^n(\phi(x))\in \sigma(U').$ It follows that 
$\phi(f_{\bk}^n(x))=\sigma^{-1}(\phi(f_{\bk}^n(x)))\in U',$ which concludes the proof. 
\endproof

\subsection{Invariant curves}
In this section, we assume that $X$ is a surface.

\begin{pro}\label{proinfinvcurve}
If the pair $(X,f)$ does not satisfy the SAZD-property, then
 there exists $m\geq 1,$ such that  there exist infinitely many irreducible curves
$C$ of $X$ satisfying $f^m(C)\subseteq C.$ 
\end{pro}

\proof
Let $K$ be a subfield of $\bk$ which is finitely generated over $\Q$, satisfying $\overline{K}=\bk$ and such that $X$ and $f$ are defined over $K.$


Apply Proposition \ref{proinvpoly}; in particular $V\simeq (\C_p^{\circ})^d$ is $f^m$-invariant.
Observe that for every point $x\in V\setminus \Fix(f^m)$, the orbit of $x$ is well-defined and infinite. We may assume that $V\setminus \Fix(f^m)\neq\emptyset.$

Let  $B$ be any proper Zariski closed subset of $X$ containing $\Fix(f)$.
Since $(X,f)$ does not satisfy the SAZD-property, there exists $z\in X_K(\iota,V\setminus B)$, whose orbit is not Zariski dense.
 There exists $\overline{\iota}\in \sI_{\iota}$ such that $\phi_{\overline{\iota}}(z)\in V.$ 
 By Proposition \ref{proinvpoly}, the Zariski closure $Z_z$ of $O_{f^m}(z)$ in $X$ is irreducible. Then we have $f^m(Z_z)=Z_z$ and $Z_z\not\subseteq B.$
 Since $z$ is not preperiodic, $\dim Z_z=1.$  
 It follows that for every proper Zariski closed subset of $X$, there exists an irreducible and $f^m$-invariant curve $C$ of $X$ which is not contained in $B.$ This concludes the proof.
\endproof

The following result generalizes \cite[Theorem 1.3]{Xie2015} and \cite[Theorem 1.3.]{Bell} in the adelic setting.
It is a direct consequence of Proposition \ref{proinfinvcurve} and \cite[Theorem B]{Invarianthypersurfaces}.
\begin{cor}\label{corbirsurzd}Assume that $f$ is a birational self-map on the surface $X$, then the pair $(X,f)$ satisfies the AZD-property.
\end{cor}
%

\subsection{Skew-linear self-maps}
In this section, we prove the following adelic version of \cite[Theorem 1.4]{Ghioca2018a}.
\begin{thm}\label{thmzardenseratadelic}
Let $g:X\dashrightarrow X$ be a dominant rational map defined over $\bk$ and $N\geq 1$. Let $f:X\times\A^N\dashrightarrow X\times\A^N$ be defined by $(x,y)\mapsto (g(x),A(x)y)$ where  $A\in \GL_{N}(\bk(X))$. Let $h:X\times\P^{N-1}\dashrightarrow X\times\P^{N-1}$ be defined by $(x,y)\mapsto (g(x),C(x)y)$ where  $C\in \PGL_{N}(\bk(X)).$
If the pair $(X,g)$ satisfies the AZD-property, then the pairs $(X\times\A^N,f)$ and $(X\times\P^N,h)$ also satisfy the AZD-property.
\end{thm}

\proof
Assume that the pair $(X,g)$ satisfies the AZD-property.

We first prove that $(X\times\A^N,f)$ satisfies the AZD-property.
Denote by $\pi:X\times\A^N\to X$ the first projection. 
Let $\sB$ be the set of points $x\in X$ such that $f$ is not locally an isomorphism along the fiber $\pi^{-1}(x)$. Then $\sB$ is a proper closed subset of $X$. 

Since $\pi^*$ maps $\bk(X)^g\setminus \bk$ to $\bk(X\times \A^N)^f\setminus \bk$, we may assume that $\bk(X)^g=\bk.$
Then there exists a nonempty adelic open subset $A_0$ of $(X\setminus \sB)(\bk)$ such that for all point $x\in A_0$, the orbit of $x$ under $g|_{X\setminus \sB}$ is well-defined and Zariski dense in $X$.

\smallskip

Let $I$ be the set of all $f$-invariant subvarieties $V\subseteq X\times \A^N$ for which every irreducible component of $V$ dominates $X$ under the projection $\pi$.  
Then \cite[Theorem 2.1]{Ghioca2018a} yields that (perhaps, at the expense of  replacing $f$ by a suitable iterate) 
 there exists an irreducible variety $Y$ endowed with a dominant rational self-map 
$g':Y\dashrightarrow  Y$
and a generically finite dominant rational map $\tau: Y\dashrightarrow X$ satisfying $\tau\circ g'=g\circ \tau$ 
 such that there exists 
a birational map $h_1$ on $Y\times \A^N=Y\times_{X} X\times\A^N$ of the form  $(x,y)\mapsto (x,T(x)y)$ where $T(x)\in \GL_N(\bk(Y))$ such that  for any subvariety $V\in I$, we have 
$h_1^{-1}((\tau\times\id)^{\#} (V))=Y\times \alpha(V)\subseteq Y\times \A^N,$  where $\alpha(V)$ is a subvariety of $\A^{N}$ and $(\tau\times\id)^{\#} (V)$ is the strict transform of $V$ by the rational map $\tau\times\id$.  Let $f':Y\times \A^N\dashrightarrow Y\times \A^N$ be the rational self-map defined by 
$$f':=g'\times_{(X,g)} f:(x,y)\mapsto (g'(x), A(\tau(x))y).$$   
We have $(\tau\times \id)\circ f'=f\circ (\tau\times \id)$.
Set $F:=h_1^{-1}\circ f'\circ h_1: Y\times \A^N\to Y\times \A^N.$
Then $F$ is the map $(x,y)\mapsto (g'(x),B(x)y)$ where $B(x):=T^{-1}(g'(x))A(\tau(x))T(x)$.
Set $\rho:=(\tau\times \id)\circ h$. Then $\rho\circ F=f\circ \rho$. For any $V\in I$,  $\rho^{\#}(V)$ is $F$-invariant and it has the form $Y\times \alpha(V).$
Denote by $G^V$ the subgroup of $\GL_N(\bk)$ consisting for $\sigma\in \GL_N(\bk)$ satisfying $\sigma\alpha(V)=\alpha(V).$ Set $G:=\cap_{V\in I}G^V.$

After replacing $Y$ by some smaller open subset, we may assume that $\tau$ and $\rho$ are regular morphisms.  Furthermore, we may assume that $\tau,\rho$ are quasi-finite and \'etale. Then $\tau^{-1}(A_0)$ is a nonempty adelic open subset of $Y.$
Let $p: Y\times \A^N\to Y$ 
be the first projection. 
Let $\sB'$ be the set of points $x\in Y$ such that $F$ is not locally an   isomorphism along the fiber $p^{-1}(x)$. Then $\sB'$ is a proper closed subset of $Y$. 
Since $Y\times \alpha(V)$ is invariant by $F$ for every $V\in I$, we get $B(x)\in G$ for $x\in Y\setminus \sB'.$
By Lemma \ref{lemdesfiniteadelic}, we only need to show that $(Y\times \A^N, F)$ satisfies the AZD-property.

By Proposition \ref{prodmlad}, there exists a nonempty adelic open subset $A_1$ of $Y\setminus \sB'$ such that for every point $x\in A_1$, its orbit under $g'|_{Y\setminus \sB'}$ is well-defined.
Since $\tau$ is \'etale, $\tau(A_1)$ is a nonempty adelic  open subset of $X(\bk).$ For every point $x\in \tau^{-1}(A_0)\cap A_1$, the orbit of $x$ is well-defined, contained in $Y\setminus \sB'$ and Zariski dense in $Y.$

By \cite[Theorem 2]{Rosenlicht1956},
 either there exists $\phi\in \bk(\A^N)\setminus \bk$ such that $\phi\circ \sigma=\phi$ for all $\sigma\in G$ or there exists a nonempty $G$-invariant Zariski open subset  $U_G\subseteq \A^N$ such that for every $y\in U_G$, $G\cdot y$ is Zariski dense in $\A^N.$
 
 First assume that the later holds. Then $(A_1\cap \tau^{-1}(A_0))\times U_G(\bk)$ is a nonempty adelic open subset of $Y\times \A^N.$
 For every $q\in (A_1\cap \tau^{-1}(A))\times U_G(\bk)$, the orbit of $q$ is well-defined and contained in $p^{-1}(Y\setminus \sB').$
 We need to show that that the orbit $O_F(q)$ is Zariski dense in  $Y\times \A^N.$ Denote by $Z$ the Zariski closure of $O_F(q)$.  Since $O_{g'}(p(q))$ is Zariski dense in $Y$, then $Z$ has at least one irreducible component which dominates $Y$. Let $W$ be the union of all irreducible components of $Z$ which dominate $Y$; then $\overline{\rho(W)}\in I$ and $\overline{\rho(W)}\neq Y\times \A^N$. There exists $m\geq 0$ such that  $F^m(q)\in \rho^{\#}(\overline{\rho(W)})=Y\times \alpha(\overline{\rho(W)})$ and so, $F^n(q)\in Y\times \alpha(\overline{\rho(W)})$ for all $n\geq m.$ 
 Write $q=(x,y)\in (A_1\cap \rho^{-1}(A))\times U_G(\bk),$ we have 
 $F^m(q)=(g'^m(x), B(g'^{m-1}(x))\dots B(x)y).$
 It follows that $B(g'^{m-1}(x))\dots B(x)y\in \alpha(\overline{\rho(W)}).$
 Since $B(z)\in G\subseteq G^{\overline{\rho(W)}}$ for $z\in Y\setminus \sB',$ we get 
 $y\in \alpha(\overline{\rho(W)}).$ It follows that $G\cdot y\subseteq \alpha(\overline{\rho(W)})$, which is not Zariski dense in $\A^N$.
Then we get a contradiction.

Now we assume that there exists $\phi\in \bk(\A^N)\setminus \bk$ such that $\phi\circ \sigma=\phi$ for all $\sigma\in G_\alpha$.
Let $\chi$ be the rational function on $Y\times \A^N$ defined by $(x,y)\mapsto \phi(y)$. Then $\chi$ is nonconstant and is invariant by $F$.

\medskip

Now we prove that $(X\times\P^{N-1},h)$ satisfies the AZD-property. 
Denote by $\pi:X\times\P^{N-1}\to X$ the first projection.
Let $\sB$ be the set of points $x\in X$ such that $h$ is not a locally isomorphism along the fiber $\pi^{-1}(x)$. Then $\sB$ is a proper closed subset of $X$. 


Since $\pi^*$ maps $\bk(X)^g\setminus \bk$ to $\bk(X\times \P^{N-1})^h\setminus \bk$, we may assume that $\bk(X)^g=\bk.$
Then there exists a nonempty adelic open subset $A_0$ of $(X\setminus \sB)(\bk)$ such that for every $x\in A_0$, the orbit of $x$ under $g|_{X\setminus \sB}$ is well-defined and is Zariski dense in $X$.

Let $J$ be the set of all $h$-invariant subvarieties in $X\times \P^{N-1}$ for which every irreducible component of $V$ dominates $X$ under the projection $\pi$.  
The following result is an analogue of \cite[Theorem 2.1]{Ghioca2018a} in this setting. 
\begin{lem}\label{leminvpn}At the expense of  replacing $h$ by a suitable iterate,
 there exists an irreducible variety $Y$ endowed with a dominant rational self-map 
$g':Y\dashrightarrow  Y$
and a generically finite dominant rational map $\tau: Y\dashrightarrow X$ satisfying $\tau\circ g'=g\circ \tau$ 
 such that there exists 
a birational map $\beta$ on $Y\times \P^{N-1}=Y\times_{X} X\times\P^{N-1}$ of the form  $(x,y)\mapsto (x,T(x)y)$ where $T(x)\in \PGL_N(\bk(Y))$ such that  for any subvariety $V\in J$, we have 
$$\beta^{-1}((\tau\times_X\id)^{\#} (V))=Y\times \gamma(V)\subseteq Y\times \P^{N-1},$$  where $\gamma(V)$ is a subvariety of $\P^{N-1}$.    
\end{lem}

After replacing \cite[Theorem 2.1]{Ghioca2018a} by Lemma \ref{leminvpn}, the proof above for $(X\times\A^N,f)$ yields directly the proof for the pair  $(X\times \P^{N-1},h)$. 
\endproof

\proof[Proof of Lemma \ref{leminvpn}]
Since $H^{1}_{\acute{e}t}(\bk(X),\G_m)=0$, there exists $D(x)\in \GL_N(\bk(X))$ whose image in $\PGL_N(\bk(X))$ is $C(x).$
Consider the rational morphism 
$f: X\times (\A^N\setminus \{0\})\dashrightarrow X\times (\A^N\setminus \{0\})$ sending $(x,y)$ to $(x, D(x)y).$
Denote $\theta:(\A^N\setminus \{0\})\to \P^{N-1}$ the morphism $(x_1,\dots,x_N)\mapsto [x_1:\dots,:x_N].$
Set $$\phi:=\id\times \theta: X\times (\A^N\setminus \{0\})\to X\times \P^{N-1}.$$
Then we have $\phi\circ f=h\circ \phi.$
Denote by $I$  the set of all $f$-invariant subvarieties in $X\times \A^{N}$ for which every irreducible component of $V$ dominates $X$ under the projection map $\pi_1: X\times \A^N\to X$. 
For every $W\in J$, define $\hat{W}:=\overline{\phi^{-1}(W)}.$ We have $\hat{W}\in I$. 

By \cite[Theorem 2.1]{Ghioca2018a},  at the expense of  replacing $h,f$ by a suitable iterate, 
 there exists an irreducible variety $Y$ endowed with a dominant rational self-map 
$g':Y\dashrightarrow  Y$, a generically finite dominant map $\tau: Y\dashrightarrow X$ satisfying $\tau\circ g'=g\circ \tau$, 
and a birational map $\hat{\beta}$ on $Y\times \A^N=Y\times_{X} X\times\A^N$ of the form  $(x,y)\mapsto (x,\hat{T}(x)y)$ where $\hat{T}(x)\in \GL_N(\bk(Y))$, such that the following holds:
For every subvariety $V\in I$, we have 
$\hat{\beta}^{-1}((\tau\times_X\id)^{\#} (V))=Y\times \alpha(V)\subseteq Y\times \A^N,$  where $\alpha(V)$ is a subvariety of $\A^{N}$ and $(\tau\times_X\id)^{\#} (V)$ is the strict transform of $V$ by the rational map $\tau\times_X\id$. 
Let $T(x)$ be the image of $\hat{T}(x)$ in $\PGL_{N}(\bk(X))$. Let $\beta: Y\times \P^{N-1}\to Y\times \P^{N-1}$ be the morphism $(x,y)\mapsto (x,T(x)y).$
For every $V\in J$, define $\gamma(V):=\theta(\alpha(\hat{V})).$ Then $\beta^{-1}((\tau\times_X\id)^{\#} (V))=Y\times \gamma(V)\subseteq Y\times \P^{N-1},$ which concludes the proof.
\endproof

\subsection{Diophantine condition}
We say that  $\la_1, \la_2\in \C_p\setminus \{0\}$ satisfy the  \emph{Diophantine condition} \cite{Herman1983}, if $|\la_1|=|\la_2|=1$ and
there exists $C, \beta>0$ such that for every $n_1,n_2\in \Z_{\geq 0}, n_1+n_2\geq 2$ and $i=1,2$, we have
$$|\la_1^{n_1}\la_2^{n_2}-\la_i|\geq C|n_1+n_2|_{\R}^{-\beta},$$
here $|\cdot|_{\R}$ is the absolute value on $\R.$

\begin{pro}\label{proembdiop}Let $\la_1,\la_2\in \bk\setminus\{0\}$ be two multiplicatively independent\footnote{We say that $\la_1,\la_2$ are \emph{multiplicatively independent} if for every $(m_1,m_2)\in \Z^2\setminus \{(0,0)\}$,
$$\la_1^{m_1}\la_2^{m_2}\neq 1.$$} elements. Then there exists a prime $p$, a positive integer $m\in \Z_{> 0}$ and an embedding $\tau:\Q(\la_1,\la_2)\hookrightarrow \C_p$ such that $\tau(\la_1)^m,\tau(\la_2)^m$ satisfy the Diophantine condition.
\end{pro}
\proof
If $\la_1,\la_2\in \overline{\Q}$, we conclude the proof by \cite[Theorem 1]{Yu1990}.

If $\trd_{\Q}\Q(\la_1,\la_2)=2$, we only need to show that there are elements $\mu_1,\mu_2\in 1+p\Z_p$ which are algebraically independent over $\Q$ and satisfy the Diophantine condition. This follows from \cite[Proposition 3]{Herman1983}.

Now assume $\trd_{\Q}\Q(\la_1,\la_2)=1$.
We may assume that $\la_1\not\in \overline{\Q}$ and $\la_2\in \overline{\Q(\la_1)}.$
Denote by $P(\la_1,x)=x^d+a_{d-1}x^{d-1}+\dots+a_0,$ where $a_i\in \overline{\Q}(\la_1), i=1,\dots, d,$
the minimal polynomial of $\la_2$ over $\overline{\Q}(\la_1)$.
Observe that $a_0\neq 0.$

There exists a finite set $S\subseteq \overline{\Q}$ such that $a_i\in \sO(\A^1\setminus S)$ for $i=1,\dots, d$.
Denote by $Y$ the curve $\{(t,y)\in (\A^1\setminus S)\times (\A^1\setminus \{0\})|\,\,y^d+a_{d-1}(t)y^{d-1}+\dots+a_0(t)=0\}$ and $\pi: Y\to \A^1\setminus S$ the projection to the first coordinate.
After enlarging $S$, we may assume that $\pi$ is \'etale.
Pick a root of unity $\mu\in \overline{\Q}\setminus S$ and a point $(\mu,c)\in \pi^{-1}(\mu).$ Observe that $c\in \overline{\Q}.$
Let $K$ be a number field who contains $\mu,c$ and all coefficients of $a_i, i=0,\dots,d-1.$
By Lemma \ref{lem:Lech-Bell}, there exists a prime $p>2$ which does not divide the order of $\mu$, an embedding $\tau_1: K\to \Q_p$ such that 
$a_i(\mu)\in \Z_p$ and $(d\pi)|_{(\mu,c)}$ is invertible modulo $p.$
Then there exists 
$\phi(t)=\sum_{i\geq 1}c_it^i,$
where $c_i\in \Z_p, c_i\to 0$ as $i\to \infty$ such that
$P(\mu+t, c+\phi(t))=0.$
There exists $m\in \Z_{>0}$, such that $p\not| m,$
$\mu^m=1$ and $c^m=1 \mod p.$
Set $\alpha(t):=(\mu+t)^m-1.$
We have $\alpha(0)=0$ and $\alpha|_{p\Z_p}: p\C_p\to p\C_p$ is an analytic automorphism.
Set $\beta(t):=(c+\phi(\alpha^{-1}(t)))^m-1.$ Then $\beta$ converges on $p\C_p$ and $\beta(p\C_p)\subseteq p\C_p.$
Observe that the coefficients of $\beta$ are contained in $\Z_p\cap \overline{\Q}.$
For every $u\in p\C_p\setminus \overline{\Q}$, there exists an embedding 
$\tau_u: \Q(\la_1,\la_2)\hookrightarrow \C_p$ sending $\la_1$ to $1+u$ and $\la_2$ to $1+\beta(u).$
Then we only need to show that there exists $u\in p\C_p\setminus \overline{\Q}$, such that 
$1+u$ and $1+\beta(u)$ satisfy the Diophantine condition.

We note that $\exp(x)$ and $\log(1+x)$ are well-defined analytic function on $p\C_p$ and satisfy:
$|\log(1+x)|=|\exp(x)-1|=|x|, x\in p\C_p.$
Set $\delta(t):=\log(1+\beta(\exp(t))),$ which converges on $p\C_p.$
Write $\delta(t)=\sum_{i\geq 0}b_it^i.$ 
\begin{lem}\label{lemexlogdi}There exists $r\in \Q_{>0}\cap [1,+\infty)$, and $C, \beta>0$ such that for every $v\in \C_p$ with norm 
$|v|=p^{-r}$,  for every $m,n\in \Z_{\geq -1}, m+n\geq 1$ and $i=1,2$, we have
$|mv +n\delta(v)|\geq C|m+n+1|_{\R}^{-\beta}.$
\end{lem}
Pick $r$ as in Lemma \ref{lemexlogdi}.
Pick $u\in p\C_p\setminus \overline{\Q}$ with norm $|u|=p^{-r}.$
Then $v:=\log(1+u)$ has norm $|v|=p^{-r}.$ We conclude the proof by Lemma \ref{lemexlogdi}.
\endproof
\proof[Proof of Lemma \ref{lemexlogdi}]
We first have the following observations:
\begin{points}
\item for $m,n\in \Z_{\geq -1}, m+n\geq 1$, we have 
$$\max\{|m|_{\R},|n|_{\R}\}\leq |m+n+1|_{\R};$$
\item for $n\in \Z\setminus \{0\}$, we have $|n|\geq |n|_{\R}^{-1}.$
\end{points}

\medskip

Write $\delta(t)=\sum_{i\geq 0}b_it^i.$  We note that $b_i\in \overline{\Q}, i\geq 0.$

\medskip

We first treat the case where $b_i=0$ for all $i\neq 1.$ Then we have $\delta(t)=b_1t.$ Since $\exp(t)$ and $\exp(\delta(t))$ is multiplicatively independent, we have 
$b_1\not\in \Q.$
Pick  any $r=1$ and $\beta:=3$. 
By the $p$-adic Thue-Siegel-Roth theorem\cite{Ridout1958}, there exists $C_1>0$ such that for $m,n\in \Z, m+n\geq 1$, we have 
$$|m+nb_1|\geq C_1\max\{|m|_{\R},|n|_{\R}\}^{-3}.$$
Set $C:=C_1p^{-1}.$
Then for every $v\in p\C_p$ with norm $|v|=p^{-1},$ $m,n\in \Z_{\geq -1}, m+n\geq 1$, 
 we have 
$$|mv +n\delta(v)|=|m+nb_1|p^{-1}\geq C_1p^{-1}\max\{|m|_{\R},|n|_{\R}\}^{-3}\geq C|m+n+1|_{\R}^{-3}.$$
This concludes the proof.

Now, we may assume that the set $\{i\in \Z_{\geq 0}\setminus \{1\}|\,\, b_i\neq 0\}$ is nonempty and let  
$s$ be the smallest integer in this set.
There exists $l\in \Z_{\geq 1}$ such that $|b_sp^{ls}|>|b_ip^{li}|$ for all $i\neq 1,s.$
Pick $r:=l+1/(s+2)$. Let $v$ be any element with norm $|v|=p^{-r}.$ 
If $n=0,m\geq 1$, we get 
$$|mv +n\delta(v)|=|mv|\geq p^{-r}|m|_{\R}^{-1}\geq  p^{-r}|m+n+1|_{\R}^{-1}.$$
For every $n\in \Z\setminus \{0\}$, We have 
\begin{points}
\item[$\d$] $|nb_sv^s|>|n(\sum_{i\geq 0, i\neq 1,s}b_iv^i)|;$
\item[$\d$] $|nb_sv^s|\neq |(m+nb_1)v|.$
\end{points}
It follows that for $m,n\in \Z_{\geq -1}, m+n\geq 1,n\neq 0$, 
we have 
$$|mv +n\delta(v)|=|(m+nb_1)v+nb_sv^s+n(\sum_{i\geq 0, i\neq 1,s}b_iv^i)|$$
$$=\max\{|(m+nb_1)v|,|nb_sv^s|\}\geq |nb_s|p^{-sr}$$$$\geq |b_s|p^{-sr}|n|_{\R}^{-1}\geq   |b_s|p^{-sr}|n+m+1|_{\R}^{-1}$$
We conclude the proof by setting $\beta:=1$ and $C:=\min\{|b_0|p^{-sr}, p^{-r}\}.$
\endproof

The proofs of \cite[Proposition 2.3]{E.Amerik2011}, \cite[Lemma 2.6]{E.Amerik2011} and \cite[Corollary 2.7]{E.Amerik2011} show that:
\begin{pro}\label{promulindfixam}Let $p$ be a prime. Let $X_{\C_p}$ be an irreducible surface over $\C_p$
and let $f:X_{\C_p}\dashrightarrow X_{\C_p}$ be a dominant rational self-map.
Let $o$ be a smooth point in $X_{\C_p}(\C_p)\setminus I(f)$ satisfying $f(o)=o.$
If the eigenvalues of $df|_o$ are nonzero and 
satisfy Diophantine condition,
then for every $p$-adic neighborhood $V$ of $o$, there exists a nonempty $p$-adic open set $U\subseteq V$ such that for every point $y\in U$, the orbit of $y$ is well-defined and Zariski dense in $X_{\C_p}.$
\end{pro}

Combining Proposition \ref{proembdiop} with Proposition \ref{promulindfixam}, we get the following result.
\begin{cor}\label{cormulindfix}Let $X$ be an irreducible surface defined over $\bk$. 
and $f:X\dashrightarrow X$ be a dominant rational self-map.
Let $o$ be a smooth point in $X(\bk)\setminus I(f)$ satisfying $f(o)=o.$
If the eigenvalues of $df|_o$ are nonzero and multiplicatively independent, then the pair $(X,f)$ satisfies the SAZD-property.
\end{cor}

\section{Applications of the adelic topology}\label{sectionappadelic}
In this section, we assume $\trd_{\Q}\bk<\infty.$
\subsection{Product by endomorphisms of $\P^1$.}
Let $X$ be an irreducible projective variety over $\bk$. Let $g: X\dashrightarrow X$ be a dominant rational self-map.
Let $h:\P^1\to \P^1$ be a dominant endomorphism.
Denote by $f: X\times \P^1\dashrightarrow X\times \P^1$ the rational self-map defined by $(x,y)\mapsto (g(x),h(y)).$

\begin{thm}\label{thmprodendosup}Assume that $(X,g)$ satisfies the AZD-property (resp. SAZD-property). If $h$ has a superattracting fixed point, then $(X\times \P^1, f)$ satisfies the AZD-property (resp. SAZD-property).
\end{thm}
%
%
%
%
%

\begin{rem}
Combining the proof of Theorem \ref{thmprodendosup} and \cite[Lemma 14.3.4.1]{Bell2016} ( see Lemma \ref{lemnotpcfmanpat} also),  we can replace the assumption that $h$ has a superattracting fixed point by assuming $h$ is not post-critially finite. 
\end{rem}

\proof[Proof of Theorem \ref{thmprodendosup}] 
Denote by $\pi: X\times \P^1\to X$ the first projection.
If $H\in \bk(X)^g\setminus \bk$, then $\pi^*H\in \bk(X\times \P^1)^f\setminus \bk.$
So we only need to do the proof for SAZD-property.

We may assume that there exists a nonempty adelic open subset $A$ of X, such that for every $x\in A$, the orbit of $x$ is well-defined and Zariski dense in $X$.
Let $o$ be a  superattracting fixed point of $h.$
Let $K$ be a subfield of $\bk$ which is finitely generated over $\Q$, such that $\overline{K}=\bk$ and $f,X, h,o$ are defined over $K.$

By Proposition \ref{proinvpoly} , after replacing $f$ by a positive iterate, we may assume that  there exists a nonempty adelic open subset $B$ of $(X\times \P^1)(\bk)$ such that for every point $z\in B$, the orbit of $z$ is well-defined and its closure is irreducible.
By Proposition \ref{proinvpoly} again, after replacing $f$ by a positive iterate, we may assume that 
there exists a prime $p\geq 3$ and an embedding $i:K\hookrightarrow \C_p$, such that  there exists an open subset $V\simeq (\C_p^{\circ})^d$ of $X_K(\C_p)$,  which is invariant by $f$, the orbit of the points in $V$ are well-defined and $f|_V=\id \mod p.$
Moreover, there exists an analytic action $\Phi: \C_p^{\circ}\times V\to V$ of
$(\C_p^{\circ},+)$ on $V$ such that for every $n\in \Z_{\geq 0}$,  $\Phi(n,\cdot)=f^{mn}|_{V}(\cdot).$

Let $U$ be an invariant neighborhood of $o$ in $\P^1_K(\C_p)$ such that for every $y\in U$, $h^n(y)\to o$ when $n\to \infty.$
For every $z:=(x,y)\in V\times U$, we have $f^{p^n}(z)\to (x,o)=(\pi(z),o)$ when $n \to \infty.$
Denote by $Z_z$ the Zariski closure of the orbit of $z$. 
Then we have $(\pi(z),o)\in Z_z.$ It follows that $Z_{(\pi(z),o)}\subseteq Z_z.$
Then for every $z\in (X\times \P^1)_K(i, (V\times U)\setminus X\times\{o\} )\cap \pi^{-1}(A)\cap B$, we have
\begin{points}
\item the orbits of $z$ and $\pi(z)$ are well-defined;
\item the Zariski closure $Z_{\pi(z)}$ of the orbit of $\pi(z)$ is $X$;
\item the Zariski closure $Z_{z}$ of the orbit of $z$ is irreducible;
\item $Z_{\pi(z)}\times\{o\}\subseteq Z_z;$
\item $z\in Z_z\setminus Z_{\pi(z)}\times\{o\}.$ 
\end{points}
It follows that $Z_z=X\times \P^1$, which concludes the proof.
\endproof


\proof[Proof of Theorem \ref{thmspendopoly}]
Extend $f$ to an endomorphism of $(\P^1)^N.$
By Theorem \ref{thmzardenseratadelic}, we may assume that $\deg(f_i)\geq 2$ for all $i=1,\dots,N.$
Using Theorem \ref{thmprodendosup}, we concludes the proof by induction on the number of factors $N\geq 1$.
\endproof

\subsection{Endomorphisms of abelian varieties}\label{subsectionendoabelian}
Let $A$ be an abelian variety defined over $\bk$. Let $f: A\to A$ be a dominant endomorphism.

\begin{lemma}\label{lemfixdimone}If there is an irreducible subvariety $V$ of $A$ such that $\dim V\geq 1$ and $f|_V=\id$, then $\bk(A)^f\neq \bk.$
\end{lemma}
\proof[Proof of Lemma \ref{lemfixdimone}]
We may assume that $0\in V.$ Then $f$ is an isogeny.  We have $V\subseteq \ker(f-\id)$. So $\dim \ker(f-\id)\geq 1.$
Write the minimal polynomial of $f$ as $(1-t)^rP(t)$ where $P(1)\neq 0$. We have $r\geq 1.$ Set $N:=(\id-f)^{r-1}P(f).$
Then $\dim N(A)\geq 1$ and $N(A)\subseteq \ker(f-\id).$
Pick a nonconstant rational function $F$ on $N(A).$
Set $H:=F\circ N,$ which is a nonconstant rational function on $A$.
We have $$f^*H=F\circ N\circ f=F\circ f\circ N=F\circ \id \circ N=H,$$
which concludes the proof.
\endproof

\medskip

For every irreducible subvariety $V$ of $A$, define $$S_V:=\{a\in A|\,\, a+V= V\}.$$
Then $S_V$ is a group subvariety of $A$. Denote by $S_V^0$ the identity component of $S_V.$ Then $S_V^0$ is an abelian subvariety.

\begin{pro}\cite[Theorem 1.2]{Krieger2017}\label{proinvabel}
Let $V$ be an irreducible and $f$-invariant subvariety of $A.$
Then there is an irreducible subvariety $W\subseteq V$ with $\kappa(W)=\dim(W)=\kappa(V)$, and some iterate $f^m$, such that $V=W+S_V^{0}$ and $f^m(S_V^0+w)=S_V^0+w$ for every $w\in W.$
\end{pro}

Let $V$ be an irreducible and $f$-invariant subvariety of $A.$
Set $B:=A/S_V^{0}$ and denote by $\pi: A\to B$ the quotient morphism. There exists a unique endomorphism $f_B: B\to B$ such that  $f_B\circ \pi=\pi\circ f.$
Since $f$ is dominant, $f_B$ is dominant.
Set $V_B:=\pi(V)=\pi(W).$ 
We have $(f_B^m)|_{V_B}=\id.$
If $\dim(V_B)=0$, $V$ takes form $a+S_V^0$ where $a\in V.$

\begin{lemma}\label{lemwbpodim}Assume that $\dim(V_B)\geq 1$, then $\bk(A)^f\neq \bk.$
\end{lemma}
\proof[Proof of Lemma \ref{lemwbpodim}]
After replacing $f$ by $f^m$, we may assume that $(f_B)|_{V_B}=\id.$
Since $\pi$ is surjective, we only need to show that $\bk(B)^{f_B}\neq \bk.$
Apply Lemma \ref{lemfixdimone} for $f_B, B$ and $V_B$, we conclude the proof.
\endproof

\begin{lemma}\label{lemnoinvrtransub}If $\bk(A)^f= \bk,$
then every irreducible $f$-invariant subvariety $V$ takes form 
$a+A_0$ where $a\in A$ and $A_0$ is an abelian subvariety of $A$.
Moreover, if $f$ is an isogeny,  then $\Fix(f)$ is finite and $V\cap \Fix(f)\neq \emptyset.$ 
\end{lemma}

\proof[Proof of Lemma \ref{lemnoinvrtransub}]
By Lemma \ref{lemwbpodim},  $V$ takes form 
$a+A_0$ where $a\in A$ and $A_0$ is an abelian subvariety of $A$.

Now assume that $f$ is an isogeny. Since $f(V)=V$, $f(a)+f(A_0)=a+A_0$. It follows that 
$f(A_0)=A_0$ and $f(a)-a\in A_0.$
By Lemma \ref{lemfixdimone}, $\Fix(f)$ is finite and $(f-\id)|_{A_0}$ is an isogeny. So there is $x\in A_0$ such that 
$f(x)-x=a-f(a).$ Then $f(a+x)=a+x$ and $a+x\in a+A_0=V.$ This concludes the proof.
\endproof

\proof[Proof of Theorem \ref{thmendoabadelic}]
We may assume that $\bk(A)^f\neq \bk.$
Let $K$ be a subfield of $\bk$ which is finitely generated over $\Q$ satisfying $\bk=\overline{K}$, such $A$ and $f$ are defined over $K.$

\medskip

We first treat the case where $f$ is an isogeny. 
Denote by $A[2]$ the finite subgroup of the $2$-torsion points in $A$. 
By Lemma \ref{lemfixdimone}, $\Fix(f)$ is finite. 
After replacing $K$ by a finite field extension, we may assume that
every point in $\Fix(f)\cup A[2]$ is defined over $K.$

For every $l\in \Z$, denote by $[l]:A\to A$ the morphism $x\to lx.$ Since $[3]\circ f=f\circ [3]$, $[3](\Fix(f))\subseteq \Fix(f)$. There is $m\geq 1$ such that for every $x\in \Fix(f)$, $[3^{2m}](x)=[3^m](x).$

By Proposition \ref{proinvpoly}, after replacing $f$ by a positive iterate, we may assume that  there exists a nonempty adelic open subset $P$ of $A(\bk)$ such that for every point $z\in A$, $Z_f(z):=\overline{O_f(z)}$ is irreducible. Then we have $f(Z_f(x))=Z_f(x).$ By Lemma \ref{lemnoinvrtransub}, $Z_f(z)$ takes form $a+H$ where $H$ is an abelian subvariety of $A$ and $a\in \Fix(f)$.

We have $|A[2]|=2^{2\dim A}$.  Moreover, for every abelian subvariety $H'$ of $A$, we have 
$$|A[2]\cap H'|=|H'[2]|=2^{2\dim H'}.$$ In particular, if $A[2]\subseteq  H'$, then $H'=A.$

Pick an embedding $\tau: K\hookrightarrow \C_3.$ 
We note that $0\in A_K(\C_3)$ is an attracting fixed point for $[3].$
There exists an open neighborhood $U\subseteq A_K(\C_3)$ of $0$ such that for every $x\in U,$ $\lim_{n\to \infty}[3^n]x=0.$
Set $C:=P\cap (\cap_{y\in A[2]}A_K(\tau, y+U)),$ which is a nonempty basic adelic subset of $A(\bk).$ 

We only need to show that for every $x\in C$, $Z_f(x)=A.$
For $j\in \sI_{\tau}$, denote by $\phi_j: A(\bk)\hookrightarrow A_K(\C_3)$ the embedding induced by $j:\bk\hookrightarrow \C_3.$
For every $y\in A[2]$, there exists $j_y\in \sI_{\tau}$ such that $a_y:=\phi_{j_y}(x)\in U+y.$ Write $Z_x=a+H$ where $H$ is a semiablian subvariety of $A$ and $a\in \Fix(f)$.
Then for $n\geq 1$, $[3^{nm}](Z_x)=[3^{nm}](a+H)=[3^m](a)+H$

We note that $a_y=\phi_{j_y}(x)\in \phi_{j_y}(Z_x)$ for every $y\in A[2].$
Then for every $n\geq 0, y\in A[2]$, we have 

$$y+[3^{nm}](a_y-y)=[3^{nm}](a_y)\in \phi_{j_y}([3^m](a)+H).$$
Since $a_y-y\in U$, let $n\to \infty$, we get $y\in \phi_{j_y}([3^m](a)+H).$
Since $y$ is defined over $K$, $y\in [3^m](a)+H.$
 It follows that $A[2]\subseteq [3^m](a)+H.$ 
 In particular $0\in [3^m](a)+H$, so $[3^m](a)+H=H$. Then $H=A.$ It follows that $Z_f(x)=A$.

 Now we treat the general case. 
 Let $V$ be an irreducible subvariety of $A$ which has minimal dimension in all $f$-periodic subvarieties.
 By  Lemma \ref{lemnoinvrtransub}, $V$ is a translate of an abelian subvariety of $A.$
 After changing the origin of $A$ and replacing $f$ by a suitable iterate, we may assume that $V$ itself is an abelian subvariety of $A$ and $f(V)=V.$
 Set $B:=A/V$ and denote by $\pi: A\to B$ the quotient morphism. 
 There is endomorphism $f_B: B\to B$ such that $f_B\circ \pi=\pi\circ f.$  Since $f(V)=V$ and $f$ is dominant, $f_B$ is an isogeny.
 Then there is a nonempty adelic subset $D$ of $B(\bk)$ such that for every $x\in D$, $O_{f_B}(x)$ is Zariski dense in $B.$
 By Proposition \ref{proinvpoly}, after replacing $f$ by a positive iterate, we may assume that  there exists a nonempty adelic open subset $P$ of $A(\bk)$ such that for every point $z\in A$, $Z_f(x)$ is irreducible. 
 
 We claim that for every $x\in \pi^{-1}(D)\cap P$, $Z_f(x)=A.$ Since $O_{f_B}(\pi(x))$ is Zariski dense in $B$, $\pi(Z_f(x))=B.$ So $Z_f(x)\cap V\neq \emptyset.$
Since $f(Z_f(x)\cap V)\subseteq Z_f(x)\cap V$, there is an $f$-periodic subvariety contained in $Z_f(x)\cap V.$ The minimality of $\dim V$ implies that $V\subseteq Z_f(x).$ 
 By Lemma \ref{lemnoinvrtransub}, $Z_f(x)$ is an abelian subvariety of $A.$ Then $Z_f(x)=A.$ 
\endproof

\section{General facts of endomorphisms of projective surfaces}\label{secgen}
Let $X$ be an irreducible projective surface over $\bk$ and $f:X\dashrightarrow X$ be a dominant rational self-map.
We mainly interest in the case when $f$ is an endomorphism. 
When $f$ is an endomorphism, 
by \cite[Lemma 5.6]{fa}, $f$ is finite. 

\subsection{Amplified endomorphisms}
Assume that $f$ is an endomorphism.
Recall that $f: X\to X$ is said to be \emph{amplified} \cite{Krieger2017}, if there exists a line bundle $L$ on $X$ such that $f^*L\otimes L^{-1}$ is ample.
In particular, a polarized endomorphism is amplified. 

\begin{lem}\label{lemampit}
Let $n$ be a positive integer.
Then $f$ is amplified if and only if $f^n$ is amplified.
\end{lem}
\proof
If $f$ is amplified, then there exists a line bundle $L$ on $X$ such that $H:=f^*L\otimes L^{-1}$ is ample.
Since $f$ is finite, for every $i\geq 0$, $(f^i)^*H=(f^{i+1})^*L\otimes (f^{i})^*L^{-1}$ is ample.
It follows that 
$(f^n)^*L\otimes L^{-1}=\otimes_{i=0}^{n-1}H$ is ample.
Then $f^n$ is amplified.
\medskip

If $f^n$ is amplified, then there exists a line bundle $L$ on $X$ such that $(f^n)^*L\otimes L^{-1}$ is ample.
Set $M:=\otimes_{i=0}^{n-1}(f^i)^*L.$ Then we have $f^*M\otimes M^{-1}=(f^n)^*L\otimes L^{-1}$ is ample.
Then $f$ is amplified.
\endproof

Denote by $\Fix(f)$ the set of fixed points of $f$.  The proof of  \cite[Theorem 5.1]{fa} shows that when $f$ is amplified,  the set of periodic points of $f$ is Zariski dense and for all $n\geq 1$, $\Fix(f^n)$ is finite.

\begin{lem}\label{leminvcurvedegatltwo}Assume that $f$ is amplified. Let $C$ be an irreducible curve in $X$ satisfying $f(C)=C.$ Then the degree of $f|_C$ is at least $2$ and at most $d_f$.
In particular the normalization of $C$ is either $\P^1$ or an elliptic curve.
\end{lem}

\proof
Since $f$ is amplified, there exists a line bundle $L$ on $X$ such that $H:=f^*L\otimes L^{-1}$ is ample.
Since $f$ is finite, $\deg(f|_C)\geq 1.$ If $\deg(f|_C)=1,$ then $(f|_C)^*L|_C$ is numerically equivalent to $L|_C.$
So $(f^*L\otimes L^{-1})|_C=(f|_C)^*L|_C\otimes L|_C^{-1}$ is both ample and numerically trivial, which is a contradiction. So $\deg(f|_C)\geq 2.$ 
Then the normalization of $C$ is either $\P^1$ or an elliptic curve.

Let $x$ be a general point in $C(\bk)$, we have $\deg(f|_C)=|f|_C^{-1}(x)|\leq |f^{-1}(x)|\leq d_f,$
which concludes the proof.
\endproof

For an irreducible curve $C$ in $X$ satisfying $f(C)=C,$ denote by $\pi_C:\overline{C}\to C$ the normalization of $C$ and $f_{\overline{C}}:\overline{C}\to \overline{C}$ the  endomorphism induced by $f|_C.$
For a point $o\in \Fix(f)$ and an irreducible curve $C$ of $X$, denote by $m_C(o)$ the number of branches of $C$ centered at $o$ which is invariant by $f.$ 
We claim that if $f$ is amplified, then we have 
\begin{equation}\label{eqmco}
m_C(o)\leq [d_f+2d_f^{1/2}+1]+1.
\end{equation}
Indeed, if $C$ is not $f$-invariant, we have $m_C(o)=0$ for every $o\in \Fix(f).$ If $C$ is $f$-invariant, $\overline{C}$ is either $\P^1$ or an elliptic curve. When $\overline{C}\simeq \P^1$, we have $$|{\rm Fix}(f_{\overline{C}})|\leq \deg(f_{\overline{C}})+1\leq d_f+1\leq [d_f+2d_f^{1/2}+1]+1.$$
When $\overline{C}$ is an elliptic curve, we have 
$$|\Fix(f_{\overline{C}})|=|(f_{\overline{C}}-\id)^{-1}(0)|=|\alpha-1|^2$$ where $\alpha$  is a complex number satisfying 
$|\alpha |^2=\deg(f_{\overline{C}})\geq 2.$
It follows that 
$$|{\rm Fix}(f_{\overline{C}})|=|\alpha|^2-2 {\rm Re} (\alpha)+1\leq \deg(f_{\overline{C}})+2\deg(f_{\overline{C}})^{1/2}+1$$$$\leq d_f+2d_f^{1/2}+1\leq [d_f+2d_f^{1/2}+1]+1.$$
Since every invariant branch of $C$ corresponds to a fixed point of $f_{\overline{C}}$ in $\overline{C}$, we get 
$$m_C(o)\leq |{\rm Fix}(f_{\overline{C}})|\leq [d_f+2d_f^{1/2}+1]+1.$$

\begin{lem}\label{lemseqpinc}Assume that $f$ is amplified. Let $C$ be an irreducible curve in $X$ satisfying $f(C)=C.$ Then there exists a sequence of distinct points $o_i\in C(\bk), i\geq 0$ such that
\begin{points}
\item $o_0\in \Fix(f)\cap C$;
\item $f(o_i)=o_{i-1}$ for $i\geq 1.$
\end{points}
\end{lem}
\proof
By Lemma \ref{leminvcurvedegatltwo}, we have $\deg(f_{\overline{C}})\geq 2.$
Denote by ${\Exc}(f_{\overline{C}})$ the set of exceptional points i.e.
the points $x\in \overline{C}$ whose inverse orbit $\cup_{i\geq 0}f_{\overline{C}}^{-i}(x)$ is finite. 
We claim that $$\Fix(f_{\overline{C}})\setminus \Exc(f_{\overline{C}})\neq \emptyset.$$

Recall that $\overline{C}$ is a either $\P^1$ or an elliptic curve. 

When $\overline{C}\simeq \P^1$, it is well know that $|{\Exc}(f_{\overline{C}})|\leq 2.$
If there exists $x\in \Fix(f_{\overline{C}})$ with multiplicity at least $2$, then $x\not\in {\Exc}(f_{\overline{C}}).$ Otherwise if all fixed points of $f_{\overline{C}}$ are of multiplicity $1$, then we have 
$$|\Fix(f_{\overline{C}})|=\deg(f|_{\overline{C}})+1\geq 3> |{\Exc}(f_{\overline{C}})|.$$ We concludes the claim.

When $\overline{C}$ is an elliptic curve, $f_{\overline{C}}$ is \'etale. So we have ${\Exc}(f_{\overline{C}})=\emptyset.$
On the other hand 
$$|\Fix(f_{\overline{C}})|=|(f_{\overline{C}}-\id)^{-1}(0)|=|\alpha-1|^2$$ where $\alpha$  is some complex number satisfying 
$|\alpha |^2=\deg(f_{\overline{C}})\geq 2.$ Since $\alpha\neq 1$, we get $|\Fix(f_{\overline{C}})|>0,$ which concludes the claim.

Pick $q_0\in \Fix(f_{\overline{C}})\setminus \Exc(f_{\overline{C}})$. There exists a sequence $q_i\in \overline{C},i\geq 1$ such that $f_{\overline{C}}(q_i)=q_{i-1}.$
Then $q_i, i\geq 0$ are distinct.
Set $o_i':=\pi_C(q_i)\in C(\bk), i\geq 0$. Since $\pi_{C}$ is finite, the sequence $o_i',i\geq 0$ is infinite. We have $f(o_i')=o_{i-1}', i\geq 1$ and $o_0'\in \Fix(f).$ 
There is a maximal $t\geq 0$ such that $o_t'=o_0'.$ Set $o_i:=o_{i+t}$, then 
$o_i, i\geq 0$ are distinct, which concludes the proof.
\endproof

\subsection{Definition field of a subvariety} 
Let $K$ be a subfield of $\bk$ such that $X, f$ are defined over $K$.
\rem
There exists always such a field $K$  which is finitely generated over $\Q$.
\endrem

Set $G:=\Gal(\bk/K).$ 
It naturally acts on $X(\bk).$ For every $x\in X(\bk)$, we denote by $G_x$ the stabilizer of $x$ under this action.
For every sub-extension $K'/K$ of $\bk/K$, we write $X(K')$ for the set of points in $X(\bk)$ defined over $K'.$ We particularly interest in the case $K'=\overline{K}.$

\medskip

For a subvariety $S$ of $X$, 
define the subgroup $G_S:=\{g\in G|\,\, g(S)=S\}$ of $G$. Define $K_S:=\bk^{G_S}$, which is the smallest field extension of $K$, over which $S$ is defined.
In particular, if $S$ is $G$-invariant, then we have $K_S=K.$

Define $$G^S:=\cap_{x\in S(\bk)}G_x$$ which is a subgroup of $G_S.$ Define $K^S:=\bk^{G^S}$ which is the 
the smallest field extension of $K$ such that all points in $S(\bk)$ are defined over $K^S$. Observe that  $K^S$ is a Galois extension of $K_S$ whose Galois group $G_S/G^S$ is the image of $G_S$ in the permutation group of $S(\bk)$.
When $S$ is finite, $[K^S:K_S]$ divides $|S|!.$

\medskip

\begin{lem}\label{lemseqdefdeg}Assume that $f$ is an endomorphism. Let $p_0,\dots,p_n$ be a sequence of points in $X(\bk)$ satisfying $f(p_i)=p_{i-1}, i=1,\dots, n.$
Then we have $$[K^{\{p_0,\dots,p_n\}}: K^{\{p_0\}}]|\,\, (d_f!)^n.$$
\end{lem}
\proof[Proof of Lemma \ref{lemseqdefdeg}]
We have a filtration of fields
$$K^{\{p_0\}}\subseteq K^{\{p_0,p_1\}}\subseteq \dots\subseteq K^{\{p_0,\dots,p_n\}}.$$
We only need to show that 
$$[K^{\{p_0,\dots,p_{i+1}\}}: K^{\{p_0,\dots,p_{i}\}}]|\,\, d_f!, i=0,\dots,n-1.$$
After replacing $K$ by $K^{\{p_0,\dots,p_{i}\}}$, we only need to prove this lemma in the case $n=1$ and $K=K^{\{p_0\}}.$

Now assume $n=1$ and $K=K^{\{p_0\}}.$  Since $f^{-1}(p_0)$ is $G$-invariant, we have $K_{f^{-1}(p_0)}=K.$
Then we have 
$K=K_{f^{-1}(p_0)}\subseteq K^{\{p_0,p_1\}}\subseteq K^{f^{-1(p_0)}}.$
It follows that 
$$[K^{\{p_0,p_1\}}:K]|\,\, [K^{f^{-1(p_0)}}:K], [K^{f^{-1(p_0)}}:K] | \,\, |f^{-1}(p_0)|!$$
and $|f^{-1}(p_0)|! \,|\,\, d_f!,$
 which concludes the proof.
\endproof

Then we get the following constraint on the field of definition of  invariant curves. 
\begin{cor}\label{cordefinvcur}
Assume that $f$ is an amplified endomorphism.
Let $C$ be an irreducible curve in $X$ satisfying $f(C)=C$. Then $K^{\Fix(f)}$ is a finite field extension of $K$.
Assume that $K^{\Fix(f)}=K$, then
there exists $n\geq 1$ such that 
$[K_{C}:K]|\,\, (d_f!)^n.$
\end{cor}
\proof[Proof of Corollary \ref{cordefinvcur}]
Since $f$ is amplified, $\Fix(f)$ is finite. Then all points in $\Fix(f)$ are defined over $\overline{K}.$ Then $K^{\Fix(f)}$ is a finite field extension of $K.$

Now assume that $K^{\Fix(f)}=K.$
By Lemma \ref{lemseqpinc},  there exists a sequence of distinct points $o_i\in C(\bk), i\geq 0$ such that
\begin{points}
\item $o_0\in \Fix(f)\cap C$;
\item $f(o_i)=o_{i-1}$ for $i\geq 1.$
\end{points}

Let $M$ be an ample line bundle on $X$ defined over $K$. 
Denote by $Y$ the space of curves $D$ in $X$ satisfying $M\cdot D\leq M\cdot C$, which is a quasi-projective variety over $\bk$.
Moreover, it is defined over $K.$ So $G$ naturally acts on $Y.$

For every $i\geq 0$, denote by $H_i$ the closed subset of $Y$ consisting of curves $D\in Y$ satisfying $o_j\in D,$ for all $j=0,\dots, i.$
Then $H_i, i\geq 0$ is decreasing and $\cap_{i\geq 0}H_i=\{C\}.$ There exists $n\geq 1$ such that $\cap_{i=0}^nH_i=\{C\}.$
For every $g\in G^{\{o_0,\dots,o_n\}}$, we have $g(C)\in Y$ and $o_i\in g(C)$ for every $i=0,\dots, n.$ It follows that $g(C)\in \cap_{i=0}^nH_i=\{C\}.$
Then $G^{\{o_0,\dots,o_n\}}<G_C.$ It follows that $K_C\subseteq K^{\{o_0,\dots,o_n\}}.$
By Lemma \ref{lemseqdefdeg}, we have 
$[K^{\{o_0,\dots,o_n\}}: K^{\{o_0\}}]|\,\, (d_f!)^n.$
Since $K\subseteq K^{\{o_0\}}\subseteq K^{\Fix(f)}=K,$
we get $[K_{C}:K]|\,\, (d_f!)^n.$
\endproof

\section{Local dynamics}\label{seclocdy}
Assume $\trd_{\Q}\bk<\infty.$
Let $X$ be a smooth irreducible projective surface over $\bk$ and $f:X\dashrightarrow X$ be a dominant rational self-map.


\subsection{Fixed points}

Let $o$ be a fixed point of $f$. 
Let $\la_1,\la_2$ be the eigenvalues of the tangent map $df|_o:T_{X,o}\to T_{X,o}.$ 


\medskip

If we blow up $o$, we get a new surface $X_1$. Denote by $E$ the exceptional curve. Then $f$ induces a rational self-map $f_1$ on $X_1$. 
Assume that $df|_o$ is invertible. Then $f_1$ is regular along $E.$

If $\la_1\neq \la_2$, then there are exact two fixed points $o_1,o_2$ of $f_1$ in $E$. At $o_i,i=1,2$, $df_1|_{o_i}$ is semi-simple and the tangent vectors in $E$ is an eigenvector of $df_1|_{o_i}$. We may assume that the eigenvalue for this vector at $o_1$ is $\la_2/\la_1$ and the other eigenvalue is $\la_1$. Then the eigenvalues of $df|_{o_2}$ are $\la_1/\la_2, \la_2.$

If $\la_1=\la_2$ and $df|_o$ is semi-simple, then every point in $E$ is fixed by $f_1.$ At a point $q$ in $E$, $df|_{q}$ is semi-simple and the eigenvalues of $df_1|_{q}$ are $1, \la_1=\la_2.$
If $\la_1=\la_2$ and $df|_o$ is not semi-simple, then there exists a unique point $q$ in $E$ fixed by $f_1.$ The eigenvalues of $df|_{q}$ are $1, \la_1=\la_2.$

\medskip

If $C$ is a branch of curve centered at $o$ and invariant under $f$. Then the strict transform of $C$ in $X'$ is a branch of curve passing through a fixed point in $E$ and it is invariant by $f'.$
After a finite sequence of blowups at the center of the strict transform of $C$, we may get a strict transform $\overline{C}$ of $C$ where the composition $\pi_C:\overline{C}\to C$ of these blowups is the normalization
of $C$. The induced morphism $f_{\overline{C}}: \overline{C}\to \overline{C}$ from the blowups coincides the one induced by the normalization. Denote by $\overline{o}$ the center of $\overline{C}.$ The above computation shows that 
\begin{equation}\label{equationeigencure}df_{\overline{C}}|_{\overline{o}}=\la_1^s\la_2^t
\end{equation}
for some $s,t\in \Z.$
\medskip
\begin{lem}\label{lemmultiblowupeighenone}Assume that $df|_o$ is invertible and semi-simple. Assume that $\la_1=\mu^{m_1}$ and $\la_2=\mu^{m_2}$, where $\mu\in \bk$ and $m_1,m_2\in \Z_{>0}$ satisfying $(m_1,m_2)=1$. 
Then there exists a sequence of birational morphisms
$\pi_i:X_i\to X_{i-1}, i=1,\dots,l$ with a point $o_i\in X_i, i=0,\dots,l$ 
such that 
\begin{points}
\item $X_0=X, o_0=o;$
\item $\pi_i$ is the blowup at $o_{i-1}$;
\item $o_i$ is a fixed point of the rational map $f_i: X_i\dashrightarrow X_i$ induced by $f$;
\item $o_i$ is in the exceptional curve $E_i$ of $\pi_i$;
\item the eigenvalues of $df_i|_{o_i}, i=0,\dots,l-1$ take form $\mu^s$, for some $s\geq 1;$
\item the two eigenvalues of $df_{l-1}|_{o_{l-1}}$ are $\mu,\mu;$
\item $f_l|_{E_l}=\id$.
\end{points}
Moreover, if $K$ is a subfield of $\bk$ such that $X,f,o$ and $\mu$ are defined over $K$,
then we may ask $o_i$ to be defined over $K$ for $i=0,\dots,l.$
\end{lem}

\proof[Proof of Lemma \ref{lemmultiblowupeighenone}]
We prove the lemma by induction on $\max\{m_1,m_2\}.$
When $\max\{m_1,m_2\}=1$, we have $m_1=m_2=1.$
Define $\pi_1:X_1\to X_1$ the blowup of $o.$
Then $f_l|_{E_l}=\id.$
Let $o_1$ be any point in $E_1$ ( if $\mu\in K$, then pick $o_1\in E_1(K)$), we conclude the proof.

Assume that the lemma holds for $\max\{m_1,m_2\}\leq N$, where $N\geq 1.$
Assume that $\max\{m_1,m_2\}= N+1\geq 2.$ Since $(m_1,m_2)=1$, we have $m_1\neq m_2.$
Assume that $m_1<m_2$.
Define $\pi_1:X_1\to X_1$ to be the blowup at $o.$ 
If $\mu\in K$, the two fixed points in $E_1$ are defined over $K.$
In $E_1$, there exists a fixed point $o_1$ of $f_1$ such that 
the eigenvalues of $df_1|_{o_1}$ is $\mu^{m_1}, \mu^{m_2-m_1}.$
Since $m_2-m_1\geq 1$, $(m_1,m_2-m_1)=1$ and $\max\{m_1,m_2-m_1\}\leq m_2-1\leq N$,
we may apply the induction hypothesis the $(f_1, X_1, o_1)$ to conclude the proof.
\endproof

\medskip

\begin{defi}The fixed point $o\in X(\bk)$ is said to be \emph{good} if $df|_o$ is invertible and one of the following holds:
\begin{points}
\item[(1)] $\la_1$ and $\la_2$ are multiplicatively independent;
\item[(2)]  there exists a prime $p$ and an embedding $\tau: \bk\hookrightarrow \C_p$ such that 
$$|\tau(\la_1)+\tau(\la_2)|\leq1 \text{ and } |\tau(\la_1)||\tau(\la_2)|<1$$ where $|\cdot|$ is the p-adic norm on $\C_p.$
\end{points}
\end{defi}
\begin{rem} We note the that the condition (2) just means that both $|\tau(\la_1)|$ and $|\tau(\la_1)|$ are at most one and there exists $i=1,2$ much that $|\tau(\la_i)|<1.$
\end{rem}

\begin{defi} We say that $f$ has R-property if 
there exists a fixed point $o$ of $f$ and
an embedding $\sigma: \bk\hookrightarrow \C$ such that 
both $|\sigma(\la_1)|$ and $|\sigma(\la_2)|$ are strictly great than $1,$
where  $\la_1,\la_2$ are the eigenvalues of the tangent map $df|_o:T_{X,o}\to T_{X,o}.$ 
\end{defi}

\subsection{The existence of good fixed points}
In this section, assume that $f$ is an amplified endomorphism on $X$. Let $L$ be a line bundle on $X$ such that 
$f^*L\otimes L^{-1}$ is ample.

The aim of this section is to prove the following result.
\begin{lem}\label{lemextgoodfix}Assume that $f$ is an amplified endomorphism which has R-property.
Then either $(X,f)$ satisfied the SAZD-property or there exists $n\geq 1$, such that $f^n$ has a good fixed point. 
\end{lem}

Let $R$ be a finitely generated $\overline{\Q}$ sub-algebra of $\bk$, such that $\bk$ is the algebraically closure of $\Frac R$ and $X,f,L$ are defined over $\Frac R.$
There exists a variety $X_{\Frac R}$ over $\Frac R$ and an endomorphism $f_{\Frac R}: X_{\Frac R}\to X_{\Frac R}$, such that $X=X_{\Frac R}\times_{\Spec {\Frac R} }\Spec \bk$ and
and $f=f_{\Frac R}\times_{\Spec {\Frac R}} \id.$

After shrinking $W:=\Spec R$, we may assume that $W$ is smooth, there exists a smooth projective $R$-scheme $\pi: X_R\to W$ whose generic fiber is $X_{\Frac R}$,  $f_{\Frac R}$ extends to a finite endomorphism $f_R$ on $X_{R}$ and there exists a line bundle $L_R$ on $X_R$ such that $f_R^*L_R\otimes L_R^{-1}$ is $\pi$-ample.
For every point $t\in W(\overline{\Q})$, denote by $X_t$ the special fiber $X_R\times_{W}\Spec \overline{\Q}$ of $X_R$  over $t$. Let $L_t,f_t$ be the restriction of $L_R,f_R$ on $X_t.$

\begin{lem}\label{lemliftgoodf}Assume that there exists $t\in W(\overline{\Q})$ such that $f_t$ has a good fixed point in $X_t$. Then $f$ has a good fixed point in $X$.
\end{lem}
\proof[Proof of Lemma \ref{lemliftgoodf}]
Denote by $\Fix(f_R)$ the subscheme of $X_R$ of the fixed points of $f_R.$ It is isomorphic to the intersection of the graph of $f_R$ and the diagonal of $X_R\times_{W} X_R.$
Let $o$ be a good fixed point of $f_t\in X_t\subseteq  X_R.$ We have $o\in \Fix(f_R).$ Since $o$ is smooth, every irreducible component of $\Fix(f_R)$ passing through $o$ has absolute dimension at least $2(\dim W+2)-(\dim W+2+2)=\dim W.$  Pick $S$ an irreducible component of $\Fix(f_R)$ passing through $o$. For every $s\in W(\overline{\Q})$, $X_s\cap S\subseteq \Fix(f_s)$, which is finite.
It follows that $\pi|_S:S\to R$ is projective, quasi-finite and dominant. So it is
finite and surjective. 

Let $\la_1,\la_2$ be the eigenvalues of the tangent map $d(f_t)|_o:T_{X_t,o}\to T_{X_t,o}.$ 
Since $o$ is a good fixed point, then $d(f_t)|_o$ is invertible and one of the following holds:
\begin{points}
\item[(1)] $\la_1,\la_2$ are multiplicatively independent;
\item[(2)]  there exists a prime $p$ and an embedding $\tau: \overline{\Q}\hookrightarrow \C_p$ such that 
$$|\tau(\la_1)|,|\tau(\la_2)|\leq1 \text{ and } |\tau(\la_1)||\tau(\la_2)|<1$$ where $|\cdot|$ is the p-adic norm on $\C_p.$
\end{points}
Identify $X_{\bk}$ as the geometric generic fiber of $\pi$. There exists a point $o_{\bk}$ of $X_{\bk}$, whose Zariski closure in $X_R$ is $S$. 
Denote by $(\la_1)_{\bk},(\la_2)_{\bk}$ be the eigenvalues of the tangent map $d(f)|_{o_{\bk}}:T_{X,o_{\bk}}\to T_{X,o_{\bk}}.$
Observe that $\la_1,\la_2$ are the specializations of $(\la_1)_{\bk},(\la_2)_{\bk}$ (up to some permutation). If $\la_1,\la_2$ are multiplicatively independent, then $(\la_1)_{\bk},(\la_2)_{\bk}$ are multiplicatively independent.
Now we may assume that there exists a prime $p$ and an embedding $\tau: \overline{\Q}\hookrightarrow \C_p$ such that 
$$|\tau(\la_1)+\tau(\la_2)|\leq1 \text{ and } |\tau(\la_1)||\tau(\la_2)|<1$$ where $|\cdot|$ is the p-adic norm on $\C_p.$

The embedding $\tau$ induces embeddings $W(\overline{\Q})\hookrightarrow W(\C_p)$ and $X_R(\overline{\Q})\hookrightarrow X_R(\C_p).$
For every $x\in S(\C_p)$, we denote  by $\la_1(x),\la_2(x)$ the eigenvalues of the tangent map $d(f|_{X_{\pi(x)}})|_{x}:T_{X_{\pi(x)},x}\to T_{X_{\pi(x)},x}.$
Since $\la_1(x)+\la_2(x)$ and $\la_1(x)\la_2(x)$ are continuous functions on $S(\C_p)$, there exists a neighborhood $U\subseteq S(\C_p)$ of $o$ such that for every $x\in U$, 
$|\la_1(x)+\la_2(x)|\leq 1$ and $0<|\la_1(x)\la_2(x)|<1.$ We note that $\pi(U)$ is a nonempty open subset of $W(\C_p).$

For every $P\in R\setminus \{0\},$ denote by $Z_P$ the set $\{z\in W(\C_p)|\,\, P(z)=0\}.$ It is a nowhere dense closed subset of $W(\C_p).$ Observe that the topology of $W(\C_p)$ can be defined by a complete metric. Since $R\setminus \{0\}$ is countable, by Baire category theorem, $W(\C_p)\setminus (\cup_{R\setminus \{0\}}Z_P)$ is dense in $W(\C_p).$
It follows that $\pi(U)\setminus (\cup_{R\setminus \{0\}}Z_P)$ is not empty.
 Pick any point $z\in \pi(U)\setminus (\cup_{R\setminus \{0\}}Z_P).$ Then $z$ induces an inclusion $\tau_z: R\hookrightarrow \C_p.$ 
 It extends to a inclusion $\tau_z:\Frac(R) \hookrightarrow \C_p.$ Pick $x\in U\cap \pi^{-1}(z)$, we have $|\la_1(x)+\la_2(x)|\leq 1$ and $0<|\la_1(x)\la_2(x)|<1.$
Then $x$ induces an extension $\sigma:=\overline{\tau_z}:\bk=\Frac(R)\hookrightarrow \C_p,$ such that 
$|\sigma((\la_1)_{\bk}+(\la_2)_{\bk})|\leq 1$ and $|\sigma((\la_1)_{\bk}(\la_2)_{\bk})|<1.$ Then $f$ has a good fixed point in $X$.
\endproof

\begin{lem}\label{lemcurveattrgood}Let $o$ be a fixed point of $f$ such that $df|_o$ is invertible. Let $C$ be an irreducible curve in $X$ passing through $o.$ Assume that $f(C)=C,$ and every branch of $C$ at $o$ is invariant under $f$.
Denote by $\pi_C:\overline{C}\to C$ the normalization of $C$ and $f_{\overline{C}}:\overline{C}\to \overline{C}$ the  endomorphism induced by $f|_C.$
Let $q\in \pi_C^{-1}(o)$ and set $\mu:=df_{\overline{C}}|_q\in \bk.$  Assume that there exists an embedding $\alpha:\bk\hookrightarrow \C$ such that $0<|\alpha(\mu)|<1.$
Then there exists $n\geq 0$ such that $f^n$ has a good fixed point.
\end{lem}
\proof[Proof of Lemma \ref{lemcurveattrgood}]
After enlarging $R$, we may assume that $o,C,q$ are defined over $\Frac R$ and $\mu\in R.$
After shrinking $W$, we may assume that there exists an irreducible subscheme $C_R$ of $X_R$ whose generic fiber is $C$ and a section $o_R\in X_R(R)$ whose generic fiber is $o.$
For every point $t\in W$, denote by $C_t$ and $o_t$ the specializations of $C_R$ and $o_R.$ After shrinking $W$, we may assume that $C_t$ is irreducible for every $t\in W.$
There exists a projective morphism $\pi_{C_R}: \overline{C}_R\to C_R$ over $R$ whose generic fiber is $\pi_{C}$ and a $R$-point $q_R\in \overline{C}_R(R)$, whose generic fiber is $q$. After shrinking $W$, we may assume that for all $t\in W,$ the specialization $\pi_{C_t}: \overline{C}_t\to C_t$ of $\pi_{C_R}$ is the normalization of $C_t.$

The embedding $\alpha:R\subseteq \bk\hookrightarrow \C$ defined a point $\eta\in W(\C).$ We view $\mu$ as a function on $W(\C).$ We have $|\mu(\eta)|=|\alpha(\mu)|\in (0,1).$
There exists an euclidean open neighborhood  $U$ of $\eta$, such that $|\mu(\cdot)|\in (0,1)$ on $U.$
Pick $t\in U\cap W(\overline{\Q}),$ we have $0<|\mu(t)|<1.$ By Lemma \ref{lemliftgoodf}, we only need to prove that there exists $n\geq 0$ such that $f_t^n$ has a good fixed point in $X_t$.
Then we reduce to the case $\bk=\overline{\Q}.$

\medskip

Now we may assume that $\bk=\overline{\Q}.$
Assume that $X, f$ are defined over a number field $K$.
There exists a variety $X_K$ over $K$ and an endomorphism $f_K: X_K\to X_K$, such that $X=X_K\times_{\Spec K}\Spec \bk$ and
and $f=f_K\times_{\Spec K} \id.$
%
%
%

\medskip

Let $O_K$ be the ring of integers of $K.$
There exists a projective $O_K$-scheme $X_{O_K}$ which is 
flat over $\Spec O_K$ whose generic fiber is $X_{O_K}.$ Denote by $\pi_{O_K}:X_{O_K}\to \Spec O_K$ the structure morphism.
The endomorphism $f_{K}$ on the generic fiber extends to a rational self-map $f_{O_K}$ on $X_{O_K}.$ 

\medskip

Denote by $\pi^{O_K}_{\Z}:\Spec {O_K}\to \Spec \Z$ the morphism induced by the inclusion $\Z\hookrightarrow {O_K}.$
Let $X_{\Z}$ be the $\Z$-scheme  which is the same as $X_{O_K}$ as an absolute scheme 
with the structure morphism $\pi_{\Z}:=\pi^{O_K}_{\Z}\circ\pi_{O_K}: X_{O_K}\to \Spec \Z.$ Then $X_{\Z}$ is a projective $\Z$-scheme. 
Denote by $f_{\Z}:X_{\Z}\dashrightarrow X_{\Z}$ the rational self-map induced by $f_{O_K}.$

Since the generic fiber of $X_{\Z}$ is smooth and  $f_{\Z}$ is regular on the generic fiber, there exists a finite set $B(f,\Z)$ of primes such that $\pi_{\Z}^{-1}(\Spec \Z\setminus B)$ is smooth and $f_{\Z}$ is regular on $\pi_{\Z}^{-1}(\Spec \Z\setminus B).$ 
Set $B(f,{O_K}):=(\pi^{O_K}_{\Z})^{-1}(B(f,\Z))$, which is a finite subset of $\Spec(O_K).$ Then $\pi_{O_K}^{-1}(\Spec (O_K)\setminus B(f,{O_K}))$ is smooth and
$f_{O_K}$ is regular on $\pi_{O_K}^{-1}(\Spec (O_K)\setminus B(f,{O_K}))$.
If $x$ is a fixed point of $f^m$ for some $m\geq 1$, and $\beta_1,\beta_2$ be the eigenvalues of $df^m|_x$, then 
for every prime $p\not\in B(f,\Z)$, and every embedding $\tau: \overline{K}\hookrightarrow \C_p$, we have $|\tau(\beta_1)|,|\tau(\beta_2)|\leq 1.$

\medskip

By Lemma \ref{leminvcurvedegatltwo}, $\deg(f_{\overline{C}})\geq 2.$
Observe that $\overline{C}$ is either $\P^1$ or an elliptic curve.
Since on a complex elliptic curve, an endomorphism of degree at least $2$ is everywhere repelling, $\overline{C}$ could not be an elliptic curve.
Then we have $\overline{C}\simeq \P^1.$
Since $0<|\alpha(\mu)|<1,$ by \cite[Corollary 11.6]{Milnor1999}, $f_{\overline{C}}$ is not post-critically finite.

We need the following lemma, which is almost the same as \cite[Lemma 14.3.4.1]{Bell2016}.
\begin{lem}\label{lemnotpcfmanpat} Let $g:\P^1\to \P^1$ be an endomorphism over $\overline{\Q}$ of degree at least $2$ which is not post-critically finite. Then for every $N\geq 0$ and a finite subset $Z$ of $\P^1$, there exists a prime $p>N$, a point $x\in \P^1(\overline{\Q})$, $l\geq 1$, and an embedding $\tau:\overline{\Q}\hookrightarrow \C_p$ such that 
$x\not\in Z$, $g^l(x)=x,$ and $|\tau(d(g^l)|_x)|<1.$
\end{lem}

Denote by $J(f)$ the critical locus of $f$.  Since $o\not\in J(f)$ and $o\in C$, we have $C\not\subseteq J(f).$
Then $C\cap J(f)$ is finite. Let $P(f,C)$ be the union of the orbits of all periodic points in $C\cap J(f)$. Then $P(f,C)$ is finite.
Observe the for every $n\geq 1$, $P(f^n,C)=P(f,C).$
By Lemma \ref{lemnotpcfmanpat},  after replacing $f$ by a suitable positive iterate, there exists a prime $p\not\in B(f,\Z)$, an embedding $\tau:\overline{\Q}\hookrightarrow \C_p$ and 
$x\in \Fix(f|_{\overline{C}})\setminus \pi_C^{-1}(P(f,C))$ such that $C$ is smooth at $\pi_C(x)$ and $|\tau(d(f|_{\overline{C}})|_x)|<1.$
Set $q:=\pi_{C}(x).$ 
Since $q\not\in P(f,C)$, $df|_q$ is invertible. 
Since $d(f|_{\overline{C}})|_q$ is an eigenvalue of $df|_q,$ 
$q$ is a good fixed point of $f,$ which concludes the proof.
\endproof

\proof[Proof of Lemma \ref{lemnotpcfmanpat}]
Denote by $J(g)$ the set of critical points of $g.$
Since $g$ is not post-critically finite, there exists $b\in J(g)$ such that the orbit $O_g(b)$ of $b$ is infinite.
There exists $b_1\in \P^1(\overline{\Q})$ such that $g(b_1)=b.$ We have $b\neq b_1.$ Let $W$ be the union of all orbits of periodic points in $J(g)\cup Z.$ 
Then $W$ is finite.


\medskip
After a base change, we may assume that 
$g, b,b_1$, all points of $Z$ and all points of $W$ are defined over a number field $K.$ 
Set $T:=\{b,b_1\}\cup Z\cup W.$
Then $g$ defines a rational map $g_{O_K}:\P^1_{O_K}\dashrightarrow \P^1_{O_K}$ over $O_K.$
There exists a finite subset $B\subseteq \Spec O_K$ such that 
\begin{points}
\item[1)] $g_{O_K}$ is regular over $\Spec O_K\setminus B;$
\item[2)] for every $v\in \Spec O_K\setminus B$, the characteristic of the residue field at $v$ is strictly great than $N;$
\item[3)] for every $v\in \Spec O_K\setminus B$, the specialization of points of  $T$ are distinct.
\end{points}
For every $v\in \Spec O_K\setminus B$,  denote by $\P^1_v$ the special fiber of $\P^1_{O_K}$ at $v$, $g: \P^1_v\to \P^1_v$ the specialization of $g$ at $v$ and for every $x\in \P^1(K)$,
$r_v(x)$ the specialization of $x$ in $\P^1_v.$ By \cite[Lemma 4.1]{Benedetto2012}, there are infinitely may $v\in \Spec O_K\setminus B$, such that there exists $n\geq 1$ satisfying  $g_v^n(r_v(b))=r_v(b_1).$ It follows that $g_v^{n+1}(r_v(b))=r_v(b).$
Denote by $p$ the characteristic of the residue field at $v$. We have $p>N.$
Then $r_v(b)$ is a critical periodic point of $g_v.$  Denote by $K_v$ the completion of $K$ by $v$ and fix an embedding $K\hookrightarrow K_v\subseteq \C_p.$

Then there exists a point in $y\in \P^1(K_v\cap \overline{K})$ whose reduction is $r_v(b)$ and satisfying $g^{n+1}(y)=y.$
Since $b\not\in W$, $r_v(b)\not\in r_v(W).$ It follows that $y\not\in W.$ Since $y$ is  periodic, $y\not\in Z.$
Since the reduction of $dg^{n+1}|_y$ is $dg_v^{n+1}|_{r_v(b)}=0,$ we have $|dg^{n+1}|_y|<1.$
Extend the inclusion $K\subseteq \C_p$ to an embedding $\tau:\overline{K}\hookrightarrow \C_p$, we concludes the proof.
\endproof

\proof[Proof of Lemma \ref{lemextgoodfix}]
Since $f$ has R-property,  there exists a fixed point $o$ of $f$ and
an embedding $\sigma: \bk\hookrightarrow \C$ such that 
both $|\sigma(\la_1)|$ and $|\sigma(\la_2)|$ are strictly greater than $1,$
where  $\la_1,\la_2$ are the eigenvalues of the tangent map $df|_o:T_{X,o}\to T_{X,o}.$ 
It follows that $df|_o$ is invertible, and $\la_1,\la_2$ are not root of unity.


If $\la_1,\la_2$ are multiplicatively independent, then $o$ is a good fixed point of $f.$

Now we may assume that $\la_1,\la_2$ are not multiplicatively independent.
There exists $(m_1,m_2)\in \Z^2\setminus\{(0,0)\}$ such that 
$\la_1^{m_1}\la_2^{m_2}=1.$ 
Since $|\sigma(\la_1)|,|\sigma(\la_2)|>1$, we have $m_1m_2<0$.
We may assume that $m_1>0$ and $m_2<0.$

%

If for every embedding $\alpha: \bk\hookrightarrow \C$ we have
$|\alpha(\la_1)|\geq 1,$  then $\la_1\in \overline{\Q}.$
By product formula, there exists a prime $p$ and an embedding $\tau: \bk\hookrightarrow \C_p$ such that 
$|\tau(\la_1)|< 1.$ Since 
$|\tau(\la_1)|^{m_1}=|\tau(\la_2)|^{-m_2},$
we have $0<|\tau(\la_1)|,|\tau(\la_2)|<1.$  Then $o$ is good for $f.$

\medskip

Now we may assume that there exists an embedding $\alpha: \bk\hookrightarrow \C$ such that 
$|\alpha(\la_1)|<1.$ Since $|\alpha(\la_1)|^{m_1}=|\alpha(\la_2)|^{-m_2},$ we have $|\alpha(\la_2)|<1.$
View $X(\C)$ as a complex surface using the inclusion $\alpha:\bk\hookrightarrow \C.$
Let $\phi_{\alpha}$ be the natural morphism $\phi_{\alpha}:X(\bk)\hookrightarrow X(\C)$ induced by $\alpha.$
We note that $X(\bk)$ is dense in $X(\C)$ in the euclidean topology.
Since $o$ is an attracting fixed point of $f$ in $X(\C),$
there exists an euclidean open set $U$ of $X(\C)$ containing $o$ such that 
$f(\overline{U})\subseteq U$ and $\lim\limits_{n\to \infty}f^n(x)=o$ for every $x\in U.$
\begin{lem}\label{lemexinvcurve}If
$(X,f)$ does not satisfy the SAZD-property, then there exists an irreducible curve $C$ of $X$ over $\bk$ passing through $o$ and $m\geq 1$ such that 
$f^m(C)=C.$
\end{lem}
Assume that $(X,f)$ does not satisfies the SAZD-property. After replacing $f$ by a suitable positive iterate, we may assume that  there exists an irreducible curve $C$ of $X$ passing through $o$ such that $f(C)=C.$ Denote by $\pi_C:\overline{C}\to C$ the normalization of $C$ and $f_{\overline{C}}:\overline{C}\to \overline{C}$ the  endomorphism induced by $f|_C.$
After replacing $f$ by a suitable positive iterate, we may assume that every branch of $C$ at $o$ is invariant under $f$.
Pick $q\in \pi_C^{-1}(o)$.  It is a fixed point of $f_{\overline{C}}.$ Set $\mu:=df_{\overline{C}}|_q\in \bk.$  
By Equation \ref{equationeigencure},  there exists $l_1,l_2\in \Z$ such that $\mu=\la_1^{l_1}\la_2^{l_2}$. Then $\mu\neq 0.$
Since $f$ is attracting at $o\in X(\C)$, we have $|\alpha(\mu)|<1.$  Since $0<|\alpha(\mu)|<1,$ we conclude the proof by Lemma \ref{lemcurveattrgood}.
\endproof

\proof[Proof of Lemma \ref{lemexinvcurve}]
Let $K$ be a subfield of $\bk$ which is finitely generated over $\Q$, such that $\overline{K}=\bk$ and $X,f, o$ are defined over $K.$ 
By Proposition \ref{proinvpoly},
there exists $m\geq 1$ and a nonempty adelic open subset $B$ of $X(\bk)$
such that for every $x\in B$, the Zariski closure of the orbit $O_{f^m}(x)$ in $X$ is irreducible. 
After replacing $f$ by $f^m$, we may assume that for every $x\in B$, the Zariski closure $Z_x$ of the orbit $O_{f}(x)$ in $X$ is irreducible.

Since $f$ is finite, there exists an open neighborhood $V$ of $o$ in $U$ such that $f^{-1}(o)\cap V=\{o\}.$
There exists $l\geq 1$ such that $f^l(\overline{U})\subseteq V.$ Then the set  
$$S:=(\cup_{i\geq 0}f^{-i}(o))\cap U=(\cup_{i= 0}^lf^{-i}(o))\cap U$$ is finite.
For every $x\in U\setminus S$, $O_f(x)$ is infinite. 
Then for every $x\in X_K(\alpha|_K,U\setminus S)\cap A$, $Z_x$ is irreducible and positive dimensional.
Since $(X,f)$ does not satisfy the SAZD-property, 
there exists $x\in X_K(\alpha|_K,U\setminus S)\cap A$, such that  $\dim Z_x=1$.
Since $f^n(\phi_{\alpha}(x))\to o$ in $X_{K}(\C)$ for $n\to \infty$, we have $o\in Z_x$. This concludes the proof.
\endproof

\subsection{Invariant neighborhood}
Let $o$ be a fixed point of $f$ and let $\la_1,\la_2$ be the two eigenvalues of $df|_o$.
Let $K$ be a subfield of $\bk$ which is finitely generated over $\Q$, such that $\overline{K}=\bk$ and $X,f, o, \la_1,\la_2$ are defined over $K.$
 Let $\tau:K\hookrightarrow \overline{\Q_p}\subseteq \C_p$ be an embedding for some prime $p$.
 Assume that $$|\tau(\la_1)|,|\tau(\la_2)|\leq 1.$$ 
Let $K_p$ be the closure of $\tau(K)$ in $\C_p$ which is a finite extension of $\Q_p.$

\medskip

Let $W$ be an affine chart of $X$ containing $o.$ Assume that $W$ is defined over $K.$
Since $o$ is smooth, we may assume that $W$ is a complete intersection. 
Then $W$ can be viewed as a closed subvariety  of $\A^N$ which is defined by the ideal 
$(F_1,\dots, F_{N-2})$ where $F_i, i=1,\dots,N-2$ are contained in $K_p[x_1,\dots,x_N].$
We may assume that $o$ is the origin in $\A^{N}.$
Since $X$ is smooth at $o$,
the matrix $(\partial_{x_j}F_i(0))_{1\leq i\leq N-2,1\leq  j\leq N}$ has rank $N-2.$
Observe that the tangent plane of $W$ at $o$ in $\A^N$ is defined over $K_p.$
After a $K_p$ linear transform, we may assume that tangent plan of $W$ at $o$ in $\A^N$ is spanned by $\partial_{x_{1}}(0)$ and $\partial_{x_{2}}(0)$ and moreover the matrix of $df|_o$ under the base 
$\partial_{x_{1}}(0)$, $\partial_{x_{2}}(0)$ is a Jordan block
$$\left(\begin{array}{cc}\la_1& \epsilon \\ 0 & \la_2\end{array}\right)$$
where $ \epsilon= 0 \text{ or } 1.$
Then the matrix $(\partial_{x_j}F_i(0))_{1\leq i\leq N-2,3\leq  j\leq N}$ is invertible.
Denote by $\pi:W\to \A^2$ the projection $(x_1,\dots, x_N)\mapsto (x_{1},x_2).$
For every $l\geq 0$, denote by $U_l:=\{(x,y)\in \A^2(\C_p)|\,\, x,y\in p^l\C_p^{\circ}\}$, which is a $p$-adic neighborhood of $(0,0)$ in $\A^2(\C_p).$
By implicit function theorem, there exists a $l\in \Z_{>0}$ and an analytic morphism 
$\phi_l: U_l\to W(\C_p)\subseteq \A^N(\C_p)$ such that $\phi_l(U_l)$ is an open neighborhood of $o$ and  $$\pi\circ \phi_l=\id \text{ and } \phi_l\circ \pi|_{\phi_l(U_l)}=\id.$$
Moreover, $\phi_l$ is defined over $K_p$.

For every $n\geq l$, define $V_n:=\phi_l(U_n)$ which is a $p$-adic neighborhood of $o$ in $X_K(\C_p).$
There exists $m\geq l$ such that $f(V_m)\subseteq V_l.$
Then $f$ induces an analytic morphism $F: U_m\to U_l$.
Observe that $(0,0)$ is fixed by $F$ and 
$$dF|_{(0,0)}=\left(\begin{array}{cc}\la_1& \epsilon \\ 0 & \la_2\end{array}\right).$$
We may write $F$ as
$$F:(x_1,x_2)\mapsto (\la_1x_1+\epsilon x_2+\sum_{i,j\geq 0, i+j\geq 2}a_{i,j}x_1^ix_2^j, \la_2x_2+\sum_{i,j\geq 0, i+j\geq 2}b_{i,j}x_1^ix_2^j)$$
where $a_{i,j},b_{i,j}\in K_p.$
There exists $r\in \Z_{>0}$ such that $$\max\{|a_{i,j}|, |b_{i,j}||\,\, i,j\geq 0, i+j\geq 2\}\leq |p|^{-r+1}.$$
Then we have $F(U_r)\subseteq U_r.$
There exists an isomorphism $U:=(\C_p^{\circ})^2\to U_r$ sending $(z_1,z_2)$ to $(p^rz_1,p^rz_2).$
Then $F$ induces a morphism $G: U\to U$ taking form 
$$G:(z_1,z_2)\mapsto (\la_1z_1+\epsilon z_2+\sum_{i,j\geq 0, i+j\geq 2}p^{(i+j-1)r}a_{i,j}z_1^iz_2^j, \la_2z_2+\sum_{i,j\geq 0, i+j\geq 2}p^{(i+j-1)r}b_{i,j}z_1^iz_2^j).$$
Observe that $$|p^{(i+j-1)r}a_{i,j}|,|p^{(i+j-1)r}b_{i,j}|\leq |p|$$ for $i,j\geq 0, i+j\geq 2.$
The reduction $\widetilde{G}: \widetilde{U}=\widetilde{\C_p}^2\to \widetilde{U}$ of $G$ takes form 
$$(z_1,z_2)\mapsto (\widetilde{\la_1}z_1+\widetilde{\epsilon} z_2, \widetilde{\la_2}z_2)$$
where $\widetilde{\la_1}, \widetilde{\la_2}$ are the reduction of $\la_1$ and $\la_2$ in $\widetilde{\C_p}=\C_p^{\circ}/\C_p^{\circ\circ}=\overline{\F_p}.$
\medskip

Summarizing the above, we get the following result.
\begin{pro}\label{proexistinvarnei}
Assume that $|\la_1|,|\la_2|\leq 1.$
Then there exists an analytic diffeomorphism $\phi$ from the unit polydisk $U:=(\C_p^{\circ})^2$ to the open subset $V$ of $X_K(\C_p)$ which is defined over $K_p$ such that,
\begin{points}
\item $\phi((0,0))=o$;
\item the set $V$ is $f$-invariant;
\item the action of $f$ on $V$ is conjugate, via $\phi$, to an analytic endomropshim
on $U=(\C_p^{\circ})^2$ taking form 
$$G:(z_1,z_2)\mapsto (\la_1z_1+\epsilon z_2+\sum_{i,j\geq 0, i+j\geq 2}c_{i,j}z_1^iz_2^j, \la_2z_2+\sum_{i,j\geq 0, i+j\geq 2}d_{i,j}z_1^iz_2^j).$$
where $c_{i,j}, d_{i,j}\in pK_p^{\circ}$, $\epsilon=0$ if $df|_o$ is semi-simple and  $\epsilon=1$ if $df|_o$ is not semi-simple.
\end{points}
In particular, $G$ is defined over $K_p$ and 
the reduction of $G$ modulo $\C^{\circ\circ}$ takes form 
$$\widetilde{G}:(z_1,z_2)\mapsto (\widetilde{\la_1}z_1+\widetilde{\epsilon} z_2, \widetilde{\la_2}z_2).$$
\end{pro}

\begin{lem}\label{lemnotzdamrankone}Assume that $|\la_1|<1$, $|\la_2|=1$ and $f$ is an amplified endomorphism. Then there exists a nonempty open subset $U$ of $X_K(\C_p)$, such that for every point $x\in U$, the orbit $O_f(x)$ is Zariski dense in $X.$ 
\end{lem}

\proof[Proof of Lemma \ref{lemnotzdamrankone}]
Denote by $q$ a uniformizer of $K_p.$
Since  $|\la_1|<1$ and  $|\la_2|=1$, $df|_o$ is semi-simple.
Then the reduction of $\widetilde{G}$ takes form
$$\widetilde{G}:(z_1,z_2)\mapsto (0, \widetilde{\la_2}z_2).$$
 By Section \ref{subsectionzardenan},
 there exists $g\in K_p\{z_1,z_2\}$ taking form $g=z_1+h$ where $h\in qK_p^{\circ}\{z_1,z_2\}$
such that $Y:=\{g=0\}$ is invariant by $f$, $f|_Y$ is an isomorphism, $Y\simeq \C_p^{\circ}$
and $\cap_{n\geq 0} f^n(U)=Y.$
There exists a morphism $\psi:U\to Y$ satisfying $\psi|_Y=\id$ and $$f|_Y\circ \psi=\psi\circ f.$$ 
Since $f$ is finite, $G(U)\not\subseteq Y$.
Since $f$ is amplified, then all periodic points are isolated. It follows that $G|_Y$ is not of finite order.
Then we concludes the proof by Proposition \ref{profordtwozd} and Remark \ref{remcomparezaralan}.
\endproof

\begin{lem}\label{lemampatt}Assume that  $f$ is an amplified endomorphism, $\la_1,\la_2$ and all points in $\Fix(f)$ are defined over $K$.
Assume that $df|_o$ is invertible, $|\la_1|, |\la_2|<1$ and $\la_1,\la_2$ are not multiplicatively independent. 
Then there exists a nonempty open subset $U$ of $X_K(\C_p)$, such that for every point $x\in U$, the orbit $O_f(x)$ is Zariski dense in $X.$ 
\end{lem}
\proof[Proof of Lemma \ref{lemampatt}]
Since $\la_1,\la_2$ are not multiplicatively independent, there exists $(l_1,l_2)\in \Z^2\setminus \{(0,0)\}$ such that 
$$\la_1^{l_1}\la_2^{-l_2}=1.$$
Since $|\la_1|, |\la_2|<1$, we may assume that $l_1, l_2>0.$
After replacing $f$ by $f^{(l_1,l_2)}$ we may assume that $(l_1,l_2)=1.$
There exists $s,t\in \Z$ such that $sl_1+tl_2=1.$
Set $\mu:=\la_1^t\la_2^s.$ Since $\la_1,\la_2\in K$,we have $\mu\in K.$
Observe that $$\mu^{l_2}=\la_1^{tl_2}\la_2^{sl_2}=\la_1^{1-sl_1}\la_2^{sl_2}=\la_1(\la_1^{l_1}\la_2^{-l_2})^{-s}=\la_1.$$
The same, we have $\mu^{l_1}=\la_2.$

\medskip

We first treat the case where $df|_o$ is not semi-simple.
In this case, $\la_1=\la_2=\mu.$
Denote by $\pi_1:X_{1}\to X$ the blowup of $o.$ Denote by $E_1$ the exceptional curve and
 $f_1$ the rational self-map of $X_1$ induced by $f$. Observe that $f_1$ is regular along $E_1.$
 In $E_1$, there exists a unique fixed point $o_1$ of $f_1$ which is defined over $K.$
The two eigenvalues of $df_1|_{o_1}$ is $1, \mu.$ 
By Proposition \ref{proexistinvarnei}, 
there exists an analytic diffeomorphism $\phi$ from the unit polydisk $U:=(\C_p^{\circ})^2$ to the open subset $V$ of $X_{1,K}(\C_p)$ which is defined over $K_p$ such that,
\begin{points}
\item $\phi((0,0))=o$;
\item the set $V$ is $f_{1}$-invariant;
\item the action of $f_{1}$ on $V$ is conjugate, via $\phi$, to an analytic endomropshim
on $U=(\C_p^{\circ})^2$ taking form 
$$G:(z_1,z_2)\mapsto (z_1+\sum_{i,j\geq 0, i+j\geq 2}c_{i,j}z_1^iz_2^j, \mu z_2+\sum_{i,j\geq 0, i+j\geq 2}d_{i,j}z_1^iz_2^j).$$
where $c_{i,j}, d_{i,j}\in pK_p^{\circ}$.\end{points}
In particular,  
the reduction of $G$ takes form 
$$\widetilde{G}:(z_1,z_2)\mapsto (z_1, 0).$$
We may assume further that $Y:=\phi^{-1}(E_1)=\{z_2=0\}.$
Since $E_1$ is fixed by $f_1$, we have $d_{i,0}=0$ for $i\geq 2.$
By Section \ref{subsectionzardenan} and Remark \ref{remattractordete},
 $Y$ is invariant by $G$, $G|_Y$ is an isomorphism, $Y\simeq \C_p^{\circ}$
and $\cap_{n\geq 0} f^n(U)=Y.$
There exists a morphism $\psi:U\to Y$ satisfying $\psi|_Y=\id$ and $$f|_Y\circ \psi=\psi\circ f.$$ 
Since $f$ is finite, $G(U)\not\subseteq Y$. Since $f_1|_{E_1}$ is not of finite order, 
$G|_Y$ is not of finite order.
Then we concludes the proof by Proposition \ref{profordtwozd} and Remark \ref{remcomparezaralan}.

\medskip

Now we may assume that $df|_o$ is semi-simple.
By Lemma \ref{lemmultiblowupeighenone},
there exists a sequence of birational map
$\pi: X'\to X$ defined over $K$, an irreducible component $E$ of $\pi^{-1}(o)$ defined over $K$, and a point $o'\in E(K)$
such that 
\begin{points}
\item $\pi$ is an isomorphism above $X\setminus \{o\};$
\item $o'$ is a smooth point of $\pi^{-1}(o)$;
\item the induced rational map $f': X'\dashrightarrow X'$ is regular along $\pi^{-1}(o);$
\item the eigenvalues of $df'|_{o'}$ are $1,\mu;$
\item $f'|_{E}=\id.$
\end{points}
By Proposition \ref{proexistinvarnei} and the fact that $f'|_E=\id,$
there exists an analytic diffeomorphism $\phi$ from the unit polydisk $U:=(\C_p^{\circ})^2$ to the open subset $V$ of $X_K'(\C_p)$ which is defined over $K_p$ such that,
\begin{points}
\item $\phi((0,0))=o'$;
\item the set $V$ is $f'$-invariant;
\item the action of $f'$ on $V$ is conjugate, via $\phi$, to an analytic endomropshim
on $U=(\C_p^{\circ})^2$ taking form 
$$G:(z_1,z_2)\mapsto (z_1+z_2(\sum_{i,j\geq 0, i+j\geq 1}c_{i,j}z_1^iz_2^j), \mu z_2+z_2(\sum_{i,j\geq 0, i+j\geq 1}d_{i,j}z_1^iz_2^j)).$$
where $c_{i,j}, d_{i,j}\in p\C_p^{\circ}$. \end{points}
In particular,  
the reduction of $G$ takes form 
$$\widetilde{G}:(z_1,z_2)\mapsto (z_1, 0).$$
We have that $Y:=\{z_2=0\}=\phi^{-1}(E)$. After shrinking $U$ and replacing $f$ by a suitable iterate, we may assume that the morphism $\beta:=\pi\circ\phi:U\setminus Y\to X_K(\C_p)$ is an homeomorphism onto an open subset of $X_K(\C_p).$
Observe that $Y\simeq \C_p^{\circ},$ and  $G|_Y=\id$.
Let $q$ be a uniformizer of $K_p.$ Let $r$ be a positive integer which is prime to $([d_f+2d_f^{1/2}+1]+1)!.$ Set $W_Y:=\{(z_1,0)\in Y|\,\, |z_1|=|q|^{1/r}\}.$
Then $W_Y$ is a nonempty open subset of $Y.$
\begin{lem}\label{lemdivextnotin}
Let $L$ be any finite extension of $K$ satisfying $$[L:K]|(([d_f+2d_f^{1/2}+1]+1)!)^n$$ for some $n\geq 1.$ 
Let $\tau_L:L\hookrightarrow \C_p$ be any extension of $\tau$.
Denote by $\phi_{\tau_L}:X_K(L)\hookrightarrow X_K(\C_p)$ the induced inclusion.
Then $W_Y\cap \phi^{-1}(\phi_{\tau_L}(E_K(L)))=\emptyset.$
\end{lem}

By Section \ref{subsectionzardenan} and Remark \ref{remattractordete},
there exists a morphism $\psi:U\to Y$ such that $\psi|_Y=\id$, $\psi=\psi\circ f$
and for every point $x\in U$, $f^n(x)$ tends to $\psi(x).$
Since $\pi|_{U\setminus E}$ is a homeomorphism to its image and $f$ is injective in a neighborhood of $o$, after shrinking 
$U$, we may assume that $f(U\setminus E)\subseteq U\setminus E$.

Set $W:=\psi^{-1}(W_Y)\setminus Y$, which is a nonempty open in $U.$ 
Then $\beta(W)$ is a nonempty open subset of $X_K(\C_p)$.
Then we conclude the proof by the following lemma.
\begin{lem}\label{lemxwbarkorbdense}
For every $y\in \beta(W)$, the orbit $O_f(y)$ is Zariski dense in $X$.
\end{lem}
\endproof

\proof[Proof of Lemma \ref{lemdivextnotin}]
Denote by $L_p$ the closure of $\tau_L(L)$ in $\C_p$. 
Then we have $[L_p:K_p] |[L:K]$.
Since $[L:K] |(([d_f+2d_f^{1/2}+1]+1)!)^n$, we have 
$$[L_p:K_p] |(([d_f+2d_f^{1/2}+1]+1)!)^n.$$

Let $t$ be a uniformizer of $L_p$, we have $|t|^e=|q|$ for some positive integer $e|[L_p:K_p].$
We have $|L_p|=\{0\}\cup |q|^{e^{-1}\Z}.$ We note that $|q|^{1/r}\not\in |L_p|$.

Since $\phi$ is defined over $K_p\subseteq L_p,$ we have $\phi^{-1}(E(L_p))\subseteq Y(L_p):=\{(z_1,0)| z_1\in L_p^{\circ}\}\subseteq Y.$
In particular, for every $x=(z_1,0)\in \phi^{-1}(E(L_p))$, we have $|z_1|\neq |q|^{1/r}.$ It follows that $x\not\in W_Y.$
We conclude the proof.
\endproof

\proof[Proof of Lemma \ref{lemxwbarkorbdense}]
Extension $\tau$ to an embedding $\bk\hookrightarrow \C_p$. Using this embedding, we may view $X(\bk)$ as a subset of $X_K(\C_p).$

Assume that the orbit $O_f(y)$ is not Zariski dense in $X$.
Denote by $C$ the Zariski closure of $O_f(y)$ in $X$.  We have $\dim C=1.$

Set $x:=\beta^{-1}(y)$ which is contained in $W.$
Set $c:= \psi(x)\in W_Y.$
By Example \ref{exefixedcase},  $D_c:=\psi^{-1}(c)\simeq \C_p^{\circ}$ and it contains the orbit of $x$.
It follows that $\beta(D_c)$ is an irreducible analytic curve in $X_K(\C_p)$ which contains $O_f(y)$ and $o$.
Then $C$ is exactly the Zariski closure of $\beta(D_c)$ in $X.$ It follows that $C$ is an irreducible curve containing $o$. 
Moreover, since $f(\beta(D_c))\subseteq \beta(D_c)$, we have $f(C)=C.$
By Corollary \ref{cordefinvcur}, there exists a finite field extension $H$ over $K$ satisfying $[H:K]| (d_f!)^l$ for some $l\geq 0$ such that 
$C$ is defined over $H.$ 

Denote by $C'$ the strict transform of $C$ in $X'$. 
Then $\phi^{-1}(C')$ is a Zariski closed subset of $U.$
Observe that $Y\not\subseteq \phi^{-1}(C')$.
Then $\phi^{-1}(C')$ takes form $\sqcup_{i=1}^sD_{c_i}$ where $c_i\in Y$ and $D_{c_i}:=\psi^{-1}(c_i)\simeq \C_p^{\circ}$.
We may assume that $c_0=c.$
 Then $D_c$ is the unique irreducible component of $\phi^{-1}(C')$ which meets $c.$
It follows that $C'$ has only one branch passing through $\phi(c).$

Denote by $\pi_C:\overline{C}\to C$ the normalization of $C$ and $f_{\overline{C}}$ the endomorphism induced by $f|_C.$
Set $F_o:= \Fix(f_{\overline{C}})\cap \pi_C^{-1}(o).$ We have $|F_o|\leq m_C(o)\leq [d_f+2d_f^{1/2}+1]+1.$
Since $C$ and $o$ are defined over $H$, $F_o$ is defined over $H.$ 
Then there exists a finite field extension $I$ over $H$ satisfying $[I:H]|([d_f+2d_f^{1/2}+1]+1)!$ such that 
every point in $F_o$ is defined over $I.$
We note that $$[I:L]=[I:H][H:L]|([d_f+2d_f^{1/2}+1]+1)!(d_f)^l.$$

The rational map $\pi^{-1}\circ\pi_{C}:\overline{C}\dashrightarrow C'$  extends to a morphism $\pi_{C'}:\overline{C}\to C'.$
It is the normalization of $C'.$ The morphism $\pi_{C'}$ is defined over $H.$
It follows that the image of every point of $F_o$ in $X'$ is defined over $I.$ 
Then we have $\phi(c)\in \pi_{C'}(F_o)\subseteq E_K(I).$ 
Then we  have $c\in W_Y\cap \phi^{-1}(E_K(I)).$ Since $[I:L]$ divides some power of $([d_f+2d_f^{1/2}+1]+1)!$, this contradicts 
Lemma \ref{lemdivextnotin}.
We concludes the proof.
\endproof

\subsection{Amplified endomorphisms of smooth surfaces}
\begin{pro}\label{prozaridenseorbitsurfendoamp}Let $X$ be a smooth projective surface over $\bk$.
Let $f:X\to X$ be an amplified endomorphism. Assume that $f$ satisfies the R-property.
Then the pair $(X,f)$ satisfies the SAZD-property.
\end{pro}

\proof[Proof of Proposition \ref{prozaridenseorbitsurfendoamp}]
By Lemma \ref{lemextgoodfix}, after replacing $f$ by a positive iterate, we may assume that $f$ has a good fixed point $o\in X(\bk)$.
Let $\la_1,\la_2$ be the two eigenvalue of $df|_o$. 
By Corollary \ref{cormulindfix}, we may assume that $\la_1,\la_2$ are not multiplicatively independent and there exists an embedding $\tau:\bk\hookrightarrow \C_p$ for some prime $p$ such that 
$|\tau(\la_1)|,|\tau(\la_2)|\leq 1$ and $0<|\tau(\la_1)||\tau(\la_2)|<1.$
Let $K$ be a subfield of $\bk$ which is finitely generated over $\Q$ such that $\overline{K}=\bk$ and $X,f, o,\la_1,\la_2$ are defined over $K$.
Lemma \ref{lemnotzdamrankone} and Lemma \ref{lemampatt} show that there exists a nonempty open subset $U\subseteq X_K(\C_p)$ such that for every $x\in X_K(\tau|_K, U)$, the orbit of $x$ is Zariski dense in $X.$ This concludes the proof
\endproof
Let $X$ be a smooth projective variety over $\bk$.
Let $f:X\to X$ be a dominant endomorphism. We denote by $\la_1(f)$ the first dynamical degree i.e. 
$$\la_1(f):=\lim_{n\to \infty}((f^*)^nL\cdot L)^{1/n}$$ where $L$ is an ample line bundle on $X.$
The limit always exists and does not depend on the choice of the ample line bundle $L$.


\begin{cor}\label{cortopdomampz}Let $X$ be a smooth projective surface over $\bk$.
Let $f:X\to X$ be an endomorphism with $d_f>\la_1(f)$.
Then the pair $(X,f)$ satisfies the SAZD-property.

In particular, when $X=\P^2$ and $f$ is an endomorphism of $\P^2$ of degree at least $2$, then the pair $(X,f)$ satisfies the SAZD-property.
\end{cor}

\proof[Proof of Corollary \ref{cortopdomampz}]
The following Lemma shows that $f$ is amplified.
\begin{lem}\label{lemindamplify}Let $X$ be a smooth projective surface over $\bk$ and $f: X\to X$ be a dominant endomorphism. If $d_f>\la_1(f)$, then there is an ample line bundle $L$ of $X$ such that $f^*L\otimes L^{-1}$ is ample. In particular, $f$ is amplified.
\end{lem}

Pick any embedding $\sigma:\bk\hookrightarrow \C$. 
We view $X(\C)$ as a complex surface induced by $\sigma.$
Since $d_f>\lambda_1(f)$,
by \cite[Theorem 3.1, iv)]{Guedj2005} (see also \cite[Theorem 1.1]{Dinh2015}),  there exists a repelling periodic point of $f.$
It implies that $f^n$ has R-property for some $n\geq 1.$
We conclude the proof by Proposition \ref{prozaridenseorbitsurfendoamp}.
\endproof

\proof[Proof of Lemma \ref{lemindamplify}]
We note that the pull-back $f^*: N^1(X)\to N^1(X)$ and the push-forward $f_*: N^1(X)\to N^1(X)$ are dual to each other via the prefect paring on $N^1(X)$ defined by the intersection number.
So the set of eigenvalues of $f^*$ and $f_*$ are the same. We denote it by $S.$ On the other hand $f_*f^*=d_f\id.$ This shows that 
$S=d_f/S:=\{d_f/\mu|\,\, \mu\in S\}.$
Then $\la_1(f)=\max S=d_f/\min S.$ Since $d_f>\la_1(f)$, $\min S>1.$
We conclude the proof by \cite[Theorem 1.1]{Meng2020}.
\endproof

\section{Proof of Theorem \ref{thmzaridenseorbitsurfendoadelic}}\label{secproof}
\proof Assume $\trd_{\Q}\bk<\infty.$
Let $X$ be a smooth projective surface over $\bk$.
Let $f:X\to X$ be a dominant endomorphism.

If $f$ is an automorphism, we conclude the proof by Corollary \ref{corbirsurzd}.
Now we may assume that $d_f\geq 2.$

\medskip

If $\kappa(X)=2$, by \cite[Proposition 2.6]{Fujimoto2007}, $f$ is an automorphism, which concludes the proof.

\medskip

Recall the following result \cite[Lemma 2.3 and Proposition 3.1]{Fujimoto2007}.

\begin{lem}\label{lemkodeta}If the Kodaira dimension of $X$ is nonnegative and $f$ is not an automorphism, then $X$ is minimal and $f$ is \'etale.
\end{lem}
If $\kappa(X)=1$,  by \cite[Section 8]{Matsuzawa2018}, there exists a projective curve $B$, surjective morphism $\pi: X\to B$ and $m\geq 1$ such that 
$\pi\circ f^m=\pi.$ Then we concludes the proof by Lemma \ref{leminvratfunite}.

\medskip

Now we assume that $\kappa(X)=0$.
Then $X$ is either an abelian surface, a hyperelliptic surface, a K3 surface, or an Enriques surface. Since $f$ is \'etale, by \cite[Corollary 2.3]{Fujimoto2007}, we have 
$\chi(X,\sO_X)=d_f\chi(X,\sO_X).$
Since $d_f\geq 2$, we have $\chi(X,\sO_X)=0.$ Then $X$ is either an abelian surface or a hyperelliptic surface, because K3 surfaces and Enriques surfaces have nonzero Euler characteristics. 
When $X$ is an abelian surface, we concludes the proof by Theorem \ref{thmendoabadelic}.
Now we may assume that $X$ is a hyperelliptic surface.

Let $\pi:X\to E$ be the Albanese map of $X$. Then $E$ is a genus one curve,  $\pi$ is a surjective morphism with connected fibers. 
There exists a morphism $g: E\to E $ satisfying $g\circ \pi=\pi\circ f$.
Moreover there is an \'etale cover $\phi: E'\to E$ such that
$X':=X\times_E E'=F\times E'$,where $F$ is a genus one curve. Denote by $\pi_1: X'=X\times_E E'\to X$ the first projection, which is a finite \'etale morphism.
By \cite[Lemma 6.3]{Matsuzawa2018}, after a further \'etale base change, we may assume that there exists an endomorphism $g':E'\to E'$ such that 
$\phi\circ g'=g\circ \phi.$
Define $f':=f\times_E g':X'\to X'$ the induced endomorphism on $X'.$ Then we have $\pi_1\circ f'=f\circ \pi_1.$
Since $X'$ is an abelian surface, Theorem \ref{thmzaridenseorbitsurfendoadelic} holds for $f'$ by Theorem \ref{thmendoabadelic}.
Then we conclude the proof by  Lemma \ref{lemdesfiniteadelic}.

\medskip

Now assume that  $\kappa(X)=-\infty$.
Recall the following result \cite[Proposition 10]{nakayama}.
\begin{lem}\label{lemnegminf}
Assume that $\kappa(X)=-\infty$ and $f$ is not an automorphism.
Then there is a positive integer $m$ such that for every irreducible curve $E$ on $X$ with negative self-intersection, we have $f^m(E) = E.$
\end{lem}
By Lemma \ref{lemnegminf}, after replacing $f$ by $f^m$, we may assume that $f$ fixes all $(-1)$-curves.
If we contract a $(-1)$-curve of $X$ to get a new surface $X'$, $f$ induces an endomorphism $f'$ on $X'.$
By Lemma \ref{lemdesfiniteadelic}, we only need to prove Theorem \ref{thmzaridenseorbitsurfendoadelic} for $X'.$
Continue this process until there is no $(-1)$-curve, we may assume that $X$ is minimal.
Then $X$ is either $\P^2$ or a $\P^1$-bundle over a smooth projective curve $B.$
If $X=\P^2$, we conclude the proof by Corollary \ref{cortopdomampz}.

\medskip

Now we may assume that $X$ is a $\P^1$-bundle $\pi: X\to B$ over a smooth projective curve $B.$
By \cite[Lemma 5.4]{Matsuzawa2018}, after replacing $f$ by $f^2$, we may assume that there exists an endomorphism $f_B: B\to B$
such that $\pi\circ f=f_B\circ \pi.$ Denote by $d_B$ the degree of $f_B.$
For $b\in B$, set $F_b:=\pi^{-1}(b).$ Denote by $d_F$ the degree of the morphism $f|_{F_b}: F_b\to F_{f_B(b)}.$ We have $d_F\times d_B=d_f.$
Since $d_f\geq 2$, either $d_B\geq 2$ or $d_F\geq 2.$

We have $\dim N^1(X)=2.$ Denote by $A$ the nef cone of $X$ in $N^1(X).$
Denote by $F$ the class of a fiber of $\pi$ in $N^1(X).$
There exists $E\in N^1(X)$ such that the boundary of $A$ is the union of $\R_{\geq 0}F$ and $\R_{\geq 0}E.$ 
We note that for every $s,t>0$, $sE+tF$ is ample.
Since $f_*, f^*$ preserve the nef cone and $f_*f^*=d_f\id,$ we have $f^*(A)=A.$
Since $f$ preserves $\pi$, we have $f^*(F)=d_BF$. Then there exists $c>0$ such that $f^*(E)=cE.$
By Hodge index theorem, we have $(E\cdot F)>0.$ Since 
$d_f(E\cdot F)=(f^*E\cdot f^*F)=(d_BF\cdot c E),$ we get $c=d_f/d_B=d_F.$

If $d_F,d_B\geq 2,$ then 
$\la_1(f)=\max\{d_B,d_F\}< d_B\times d_F=d_f.$
We conclude the proof by  Corollary \ref{cortopdomampz}.

\medskip

Now we may assume that there is exactly one of $d_B,d_F$ equals to $1.$
In particular, $d_B\neq d_F.$
Then we have 
$$E\cdot E=d_F^{-1}(f^*E\cdot E)=d_F^{-1}(E\cdot f_*E)=d_B/d_F(E\cdot E).$$
It follows that $E\cdot E=0.$

 If $f_B$ is of finite order, then there exists
$m\geq 1$ such that 
$\pi\circ f^m=\pi.$
We concludes the proof by Lemma \ref{leminvratfunite}.
Now we assume that $f|_B$ is not of finite order.

For a curve $C$ in $X$, we denote by $[C]$ its class in $N^1(X).$
Write $C=aF+bE$. For every $m\geq 0$, we have 
$$f_*^mC=af_*^mF+bf_*^mE=ad_F^mF+bd_B^mE.$$
If $C$ is an irreducible periodic curve, then $[C]\in \R_{+}E\cup \R_+F.$
Moreover, if  $[C]\in  \R_+F,$ we have $[C]\cdot F=0$. It follows that $C$ is a fiber of $\pi.$

By Proposition \ref{proinfinvcurve}, after replacing $f$ by a suitable iterate, we may assume that 
there are infinite distinct irreducible curves $C_i, i\geq 1$ of $f$.
Since $f_B\neq \id,$ there are at most finite fixed points of $f_B.$
Then there are at most finite $C_i$ which are fibers of $\pi.$
After replacing $C_i, i\geq 1$ by a subsequence, we may assume that $[C_i]\in \R_+E$ for all $i\geq 1.$

Let $\overline{C_i}, i=1,2,3$ be the normalization of $C_i$ and let $f_{\overline{C_i}}$ be the endomorphism of $\overline{C_i}$ induced by $f|_{C_i}.$ 
Taking base changes by $\overline{C_i} \to B$ consecutively and by Lemma \ref{lemdesfiniteadelic}, we may assume that $C_1,C_2,C_3$ are sections of $\pi.$
Since $C_i\cdot C_j=0, i,j\in\{1,2,3\}$, $C_i\cap C_j=\emptyset$ when $i,j\in \{1,2,3\}$ and $i\neq j.$
For every $b\in B(\bk)$, denote by $C_{i,b}\in F_b(\bk)$ the fiber of $C_i, i=1,2,3$ over $b.$
Fix $3$ distinct points $o_i,i=1,2,3$ in $\P^1(\bk).$
There exists a unique morphism $\psi_b: \P^1\to F_b$ sending $o_i$ to $C_{i,b}$. 
The morphism $\psi: B\times \P^1\to X$ sending $(b,x)$ to $\psi_b(x)\in F_b\subseteq X$
is an isomorphism.
Then we may identify $X$ with $B\times \P^1.$ Denote by $\pi':X=B\times \P^1\to \P^1$ the projection to the second factor.
Then we may assume that $E$ is the class of a fiber of $\pi'.$
Since $f_*E \cdot E=d_BE\cdot E=0,$
$f$ preserves $\pi'.$
Then $f: B\times \P^1\to B\times \P^1$ takes form $(x,y)\to (f_B(x),g(y))$
where $g$ is an endomorphism of $\P^1$ of degree $d_F.$ For every $i\geq 0$, $[C_i]\in \R_+E$ and  $C_i$ is irreducible. 
Since $E\cdot C_i=0, i\geq 0$, $C_i$ is a fiber of $\pi'$ for every $i\geq 0$. Since $C_i$ is $f$-invariant, $\pi'(C_i), i\geq 0$ are distinct fixed points of $g$.
Then $g$ has infinitely many fixed points. It follows that $g=\id$, which concludes the proof.
%
%
%
%
\endproof


\section{Appendix A: Endomorphisms on the $\bk$-affinoid spaces}
In this appendix, we use the terminology of Berkovich spaces. See \cite{Berkovich1990,Berkovich1993} for the general theory of Berkovich spaces.
Our aim is to show that for certain endomorphism $f$ on a $\bk$-affinoid space $X$, the attractor $Y$ of $f$ is a Zariski closed subset and 
the dynamics of $f$ is  semi-conjugate to the its restriction to $Y.$  A special case of this result is used in the proof of the main theorem.
In the sequels to the papers \cite{Xied}, we will generate this result to the global setting.

\medskip

Denote by $\bk$ a complete valued field with a nontrivial non-archimedean norm $|\cdot |.$ Denote by $\bk^{\circ}:=\{f\in \bk|\,\,|f|\leq 1\}$ the valuation ring of $\bk$ and $\bk^{\circ\circ}:=\{f\in \bk|\,\,|f|< 1\}$ its maximal ideal. Denote by $\widetilde{\bk}:=\bk^{\circ}/\bk^{\circ\circ}$ the residue field.

Let $A$ be a strict and reduced $\bk$-affinoid algebra. Let $\rho_A(\cdot)$ be the
spectral norm on $A$.  
Set $X:=\sM(A).$
Denote by $\widetilde{X}:=\Spec(A^{\circ}/A^{\circ\circ})$ the reduction of $X$ and $\pi:X\to \widetilde{X}$ the reduction map. Let $f: X\to X$ be an endomorphism. Denote by $\widetilde{f}: \widetilde{X}\to \widetilde{X}$ the reduction of $f$. For any subspace $V$ of $A$, write $V^{\circ}:=V\cap A^{\circ}$, $V^{\circ\circ}:=V\cap A^{\circ\circ}$ and $\widetilde{V}:=V^{\circ}/V^{\circ\circ}.$

\medskip

For every $h\in A$,  the sequence $\rho_A((f^n)^*h), n\geq 0$ is decreasing, so the limit 
$$\rho_f(h):=\lim_{n\to \infty}\rho_A((f^n)^*h)$$ exists.
It is easy to see that $\rho_f(\cdot): A\to [0,+\infty)$ is a power multiplicative semi-norm on $A$, which is bounded by $\rho_A.$
Define $J^f$ to be the ideal of $A$ consisting of $h\in A$ satisfying $\rho_f(h)=0.$
The following result shows that for $h\in J^f$, $(f^n)^*h$ converges to $0$ uniformly.

\begin{pro}\label{proatleastlinear}There exists $b\in (0,1)$ and $m\geq 1$ such that for all $g\in J^f$, $\rho_A((f^*)^m(g))\leq b\rho_A(g).$
\end{pro}

\proof[Proof of Proposition \ref{proatleastlinear}]
%
%
Write $J^f=(g_1,\dots,g_s)$ where $\rho_A(g_i)=1, i=1,\dots,s.$ 
There exists $C>0$ such that for every $g\in J^f,$ we may write 
$g=\sum_{i=1}^sg_ih_i$ where $\rho_A(h_i)\leq C\rho_A(g).$
There exists $m\geq 1,$ such that 
$$\rho_A((f^*)^m(g_i))<(1+C)^{-1}, i=1,\dots, s.$$
Then $(f^*)^{m}(g)=\sum_{i=1}^s(f^*)^{m}(g_i)(f^*)^{m}(h_i).$
For $i=1,\dots,s,$ we have 
$$\rho_A((f^*)^{m}(g_i)(f^*)^{m}(h_i))\leq \rho_A((f^*)^{m}(g_i))\rho_A((f^*)^{m}(h_i))< (1+C)^{-1}C\rho_A(g).$$
It follows that 
$$\rho_A((f^*)^{m}(g))\leq (1+C)^{-1}C\rho_A(g)$$ for all $g\in J^f.$
Set $b:=(1+C)^{-1}C$. We conclude the proof.
\endproof

\begin{rem}
For simplifying the notation,  we write $\rho_{f^*}:=\rho_f$ and $J^{f^*}:=J^f.$
\end{rem}

\bigskip

The main result of Appendix A is the following theorem.

\begin{thm}\label{thmstattrainv}Assume that $X$ is distinguished. 
Assume that there exists a subvariety $Z\subseteq \widetilde{X}$ such that  $\widetilde{f}(\widetilde{X})=Z$ and $\widetilde{f}|_Z$ is an automorphism of 
$Z.$
Denote by $\widetilde{I}$ the ideal of $\widetilde{A}$ defined by $Z.$
Let $Y$ be the Zariski closed subset of $X$ defined by $J^f$.

Then we have 
\begin{points}
\item[1)] $\widetilde{J^f}=\widetilde{I};$ 
\item[2)] both the residue norm on $A/J^f$ w.r.t. $\rho_A$ and the spectral norm on $A/J^f$ are equal to the norm on $A/J^f$ induced by $\rho_f(\cdot);$
 \item[3)] $\widetilde{Y}\simeq \pi(Y)=Z$ where the first isomorphism is induced by the inclusion of $Y$ in $X$;
\item[4)] $f(Y)\subseteq Y;$
\item[5)] $Y$ is distinguished and $f|_Y$ is an automorphism of $Y;$
\item[6)] there exists a unique morphism $\psi:X\to Y$ satisfying $\psi|_Y=\id$ and $$f|_Y\circ \psi=\psi\circ f;$$ 
moreover there exists $C>0, \beta\in (0,1)$ such that 
for every $h\in A,$
$x\in X$ and $n\geq 0$, we have $$||h(f^n(x))|-|h(f^n(\psi(x)))||\leq C\beta^n\rho(h).$$
\end{points}
\end{thm}

\begin{rem}
Since $\widetilde{f}(\widetilde{X})\subseteq \widetilde{Z}$ and $\psi|_Y=\id$, we have $\widetilde{\psi}=\widetilde{f}|_Z^{-1}\circ \widetilde{f}.$
\end{rem}


\begin{rem}
By \cite[6.4.3, Theorem 1]{Bosch1984}, when $\bk$ is stable, $X$ is distinguished.
We note that, if $\bk$ is algebraically closed, or if the valuation on $\bk$ is discrete, then $\bk$ is stable.
\end{rem}

\rem
Under the assumption of Theorem \ref{thmstattrainv}, we have $Y=\cap_{n\geq 0}f^n(X).$
\endrem

\proof[Proof of Theorem \ref{thmstattrainv}]
For simplifying the notations, we first state an equivalent affinoid algebra statement of Theorem \ref{thmstattrainv}.
\begin{algstate}
Let $A$ be a strict and reduced $\bk$-affinoid algebra which is 
distinguished. Let $f:A\to A$ be an endomorphism.
Let $\widetilde{I}$ be a reduced ideal of $\widetilde{A}$ such that $\widetilde{f}(\tilde{I})=0$ and the induced morphism $\widetilde{f}|_{\widetilde{A}/\widetilde{I}}: \widetilde{A}/\widetilde{I}\to \widetilde{A}/\widetilde{I}$ is an isomorphism. Denote by $\tau$ the quotient morphism $\tau:A\to A/J^f.$

Then we have 
\begin{points}
\item[1)] $\widetilde{J^f}=\widetilde{I};$ 
\item[2)] both the residue norm on $A/J^f$ w.r.t. $\rho_A$ and the spectral norm on $A/J^f$ are equal to the norm on $A/J^f$ induced by $\rho_f(\cdot);$
 \item[3)] $\widetilde{A/J^f}=\widetilde{A}/\widetilde{J^f}=\widetilde{A}/\widetilde{I};$
\item[4)] $f(J_f)\subseteq J_f;$
\item[5)] $A/J^f$ is distinguished and the morphism $f|_{A/J^f}: A/J^f\to A/J^f$ induced by $f$ is an isomorphism;
\item[6)] there exists a unique morphism $\psi:A/J^f\to A$ such that 
$\tau\circ \psi=\id$ and $$\psi\circ f|_{A/J^f}=f\circ \psi;$$ 
moreover there exists $C>0, \beta\in (0,1)$ such that 
for every $h\in A,$
$x\in \sM(A)$ and $n\geq 0$, we have  
$$||f^n(h)(x)|-|(\psi\circ\tau\circ f^n)(h)(x)||\leq C\beta^n\rho(h).$$
\end{points}
\end{algstate}

\medskip

We only need to prove the algebraic form.
There exists a distinguished epimorphism  
$$\phi: \sT:=\bk\{T_1,\dots,T_r\}\twoheadrightarrow A.$$
There exists a morphism $F: \bk\{T_1,\dots,T_r\}\to \bk\{T_1,\dots,T_r\}$
such that $$\phi\circ F=f\circ \phi.$$

Denote by $K$ the kernel of $\phi.$ 
We have $F(K)\subseteq K.$
Since $\phi$ is distinguished, $\widetilde{K}$ is exactly the kernel of $\widetilde{\phi}.$
Set $\widetilde{I}_1:=\widetilde{\phi}^{-1}(\widetilde{I}).$ Since $\widetilde{f}(\tilde{I})=0$, 
$\widetilde{F}(\widetilde{I}_1)\subseteq \widetilde{K}.$
Moreover, since $\widetilde{f}|_{\widetilde{A}/\widetilde{I}}$ is an automorphism of $\widetilde{A}/\widetilde{I}$, for every $h\in \widetilde{\sT},j\geq 1$, there exists $h'\in \widetilde{\sT}$
 such that $h-(\widetilde{F})^j(h')\in \widetilde{I}_1.$ In other words, we have $\widetilde{\sT}=\widetilde{I}_1+(\widetilde{F})^j(\widetilde{\sT}).$
 
 Write $\widetilde{I}_1=(\widetilde{G}_1,\dots, \widetilde{G}_s).$
Set $\widetilde{g}_i=\widetilde{\phi}(\widetilde{G_i}), i=1,\dots,s,$
then $\widetilde{I}=(\widetilde{g}_1,\dots, \widetilde{g}_s).$ 
There are $G_i\in T^{\circ}, i=1,\dots,s$, such that $\widetilde{G}_i, i=1,\dots,s$ is the reduction of $G_i$.
Since $\widetilde{F}(\widetilde{G}_i)\in \widetilde{K}, i=1,\dots,s$, there exists $c\in (0,1)$ such that for all $i=1,\dots,s$, 
$\rho_A(f(\phi(G_i)))\leq c.$

By \cite[Chapter 2.3, Corollary 7]{Bosch2014}, we may write $K=(K_1,\dots,K_m)$ where $\rho_{\sT}(K_i)=1, i=1,\dots,m$ such that $\widetilde{K}=(\widetilde{K}_1,\dots,\widetilde{K}_m)$ and $K^{\circ}=\sum_{i=1}^mK_i\sT^{\circ}.$
For a multi-index $I=(i_1,\dots,i_r)$, denote by $T^I$ the monomial $T_1^{r_1}\dots T_r^{i_r}$.


\begin{lem}\label{lemdecompline}There are three disjoint sets $S_1,S_2,S_3$, elements
 $E_i, i\in S_1\sqcup S_2$ and $E_i^j, i\in S_3, j\geq 1$ of $\sT$
such that 
\begin{points}
\item[$\d$] for every $j\geq 1$, $E_i, i\in S:=S_1\sqcup S_2, E_i^j, i\in S_3$ is an orthonormal basis of $\sT$;
\item[$\d$] $E_i, i\in S_1$ is an orthonormal basis of $K$;
\item[$\d$] $\widetilde{E}_i, i\in S_1\sqcup S_2$ is a base of $\widetilde{I}_1$;
\item[$\d$] for every $i\in S_2$, $E_i$ takes form $G_lT^I$ for some $l\in \{1,\dots, s\}$ and some multi-index $I;$
\item[$\d$] for every $i\in S_3, j\geq 1$, $E_i^j$ takes form $(F)^j(P^j_i)$ where $P^j_i=T^I$ for some multi-index $I.$
\end{points}
\end{lem}

For every $i\in S_2$, we have $\rho_{A}(\phi(F(E_i)))\leq c.$
Set $W:=\oplus_{i\in S_2}\bk E_i$.
We note that $\widetilde{W}$ is generated by $\widetilde{E_i}, i\in S_2.$
So $\widetilde{\phi}(\widetilde{W})=\widetilde{I}.$
For every $H\in W$, we may write $H=\sum_{i\in S_2}a_iE_i$ with $\rho_{\sT}(H)=\max_{i\in S_2}|a_i|.$
Then $\rho_{A}(\phi(F(H)))\leq c\rho_{\sT}(H)$ for all $H\in W.$

For the convenience, we set $E_i^j:=E_i$ for $i\in S_1\sqcup S_2, j\geq 1.$ 
For every $H\in \sT,j\geq 1$, we may write
$H=\sum_{i\in S} a_i^jE_i^j, a_i^j\in \bk$ and $\rho_{\sT}(H)=\max_{i\in S}|a_i^j|.$
Write $$K^j(H):=\sum_{i\in S_1}a_i^jE_i\in K, 
W^j(H):=\sum_{i\in S_2}a_i^jE_i\in W, \text{ and } Q^j(H):=\sum_{i\in S_3}a_i^jP_i^j.$$
We have $$H=K^j(H)+W^j(H)+F^j(Q^j(H))$$ and $$\rho_{\sT}(H)=\max\{\rho_{\sT}(K^j(H)),\rho_{\sT}(W^j(H)),\rho_{\sT}(Q^j(H))\}.$$
Moreover since $\phi$ is distinguished, for every $j\geq 0$
$$\rho_A(\phi(H))=\max\{\rho_{\sT}(W^j(H)),\rho_{\sT}(Q^j(H))\}.$$

\bigskip

For every $H\in (K\oplus W)\cap \sT^{\circ}$, we will define sequences $H_i\in \sT^{\circ},i\geq 0, a_i\in \bk,i\geq 1$ such that $|a_i|\leq c$, 
 \begin{points}
\item[$\d$]  $F^i(H_i)\in (W\oplus K)\cap \sT^{\circ}$;
\item[$\d$] $\rho_{\sT}(W^i(F^i(H_i)))\leq \prod_{j=1}^i|a_j|;$
\item[$\d$] $\rho_{\sT}(H_j-H_i)\leq c^{i+1}$ for $j> i\geq 0$.
\end{points}
In particular, $\widetilde{H_i}=\widetilde{H}$ for $i\geq 0$ and the sequence $H_i$ converges when $i\to \infty.$

 Now we do the construction by recurrence. 
Set $H_0:=H.$ For $i\geq 0$, we have
$$F^i(H_i)=K^i(F^i(H_i))+W^i(F^i(H_i)).$$
There is $V_i\in W\cap \sT^{\circ}$ such that $$W^i(F^i(H_i))=(\prod_{j=1}^ia_j)V_i.$$
Then
$$F^{i+1}(H_i)=(\prod_{j=1}^ia_j)F(V_i)\mod K.$$
Since $V_i\in T^{\circ}\cap W$, we have 
$\rho_{A}(\phi(F(V_i)))\leq c.$
Pick $a_{i+1}\in \bk$ with $|a_{i+1}|=\rho_A(\phi(F(V_i))).$
Since $\phi$ is distinguished, we have $$F(V_i)=a_{i+1}U_{i} \mod K$$ for some $U_{i}\in T^{\circ}.$
We have $$U_i=W^{i+1}(U_i)+F^{i+1}(Q^{i+1}(U_i)) \mod K.$$

It follows that 
$$F^{i+1}(H_i)= (\prod_{j=1}^{i+1}a_j)W^{i+1}(U_i)+(\prod_{j=1}^{i+1}a_j)F^{i+1}(Q^{i+1}(U_i))\mod K,$$
thus
$$F^{i+1}\Big(H_i-(\prod_{j=1}^{i+1}a_j)Q^{i+1}(U_i)\Big)= (\prod_{j=1}^{i+1}a_j)W^{i+1}(U_i)\mod K.$$
Set  $$H_{i+1}:=H_i-(\prod_{j=1}^{i+1}a_j)Q^{i+1}(U_i).$$
We note that $\rho_{\sT}((\prod_{j=1}^{i+1}a_j)Q^{i+1}(U_i))\leq c^{i+1}.$
The sequences $H_i\in T^{\circ},i\geq 0, a_i\in \bk,i\geq 1$ are what we need.

\medskip

We claim that for every $\widetilde{g}\in \widetilde{I}$, there exists $g\in J^f\cap A^{\circ}$ such that $\widetilde{g}$ is the reduction of  $g.$
Now we prove the claim. For $\widetilde{g}\in \widetilde{I}$, there exists $\widetilde{G}\in \widetilde{I}_1$ such that $\widetilde{g}=\widetilde{\phi}(\widetilde{G}).$
Write $$\widetilde{G}=\sum_{i\in S_1\sqcup S_2}\widetilde{a_i}\widetilde{E}_i$$ where $a_i\in \bk^{\circ}.$
Set $H:=\sum_{i\in S_1\sqcup S_2}a_iE_i\in (K\oplus W)\cap T^{\circ}$. The reduction of $H$ is $\widetilde{G}.$
By the construction in the previous paragraph, we have a sequence $H_i, i\geq 0.$
Set $H_{\infty}:=\lim\limits_{i\to \infty}H_i$ and $g:=\phi(H_{\infty}).$ We have $\rho_{\sT}(H_{\infty}-H_i)\leq c^{i+1}$ for $i\geq 0$.
In particular, $\widetilde{H}_{\infty}=\widetilde{G}.$ Then the reduction of $g$ is $\widetilde{g}.$
For every $i\geq 1$, we have 
$$\rho_{\sT}(F^i(H_{\infty})-F^i(H_i))\leq c^{i+1} \text{ and }
F^i(H_i)\in K+(\prod_{j=1}^{i}a_i)W^{\circ}.$$
Then
$\rho_A(f^i(g))=\rho_A(f^{i}(\phi(H_{\infty})))\leq \max\{c^i,c^{i+1}\}=c^{i}.$
So $g\in J^f$. 

\medskip

The above argument shows that $\widetilde{I}\subseteq \widetilde{J^f}.$
Since $\widetilde{f}_Z$ is an automorphism of $Z$, we have $\widetilde{I}=\ker (\widetilde{f}^{i})$ for all $i\geq 1.$
For every $g\in J^f\cap A^{\circ}$ there exists $n\geq 1$ such that  $f^n(g)\in A^{\circ\circ}.$ 
Then we have $\widetilde{g}\in \ker(\widetilde{f}^{n})=\widetilde{I}.$
Then we get $\widetilde{I}=\widetilde{J^f},$
which proves 1).

Let $Y$  be the Zariski closed subset of $X$ defined by the ideal $J^f.$
We have $Y=\sM(A/J^f).$
Denote by $\|\cdot\|_{A/J^f}$ the residue norm on $A/J^f$.
Denote by $\rho_{A/J^f}(\cdot)$ the spectral norm on $A/J^f.$ 
We still denote by  $\rho_f(\cdot)$ the norm on $A/J^f$ induced by $\rho_f(\cdot)$ on $A$.
Since $\rho_f(\cdot)\leq \rho_A(\cdot)$ on $A$ and it is power-multiplicative,
we have that, on  $A/J^f$, $$\|\cdot\|_{A/J^f}\geq \rho_{A/J^f}(\cdot)\geq \rho_f(\cdot).$$
To prove 2), we only need to show that for every $\overline{g}\in A/J^f$, 
$\rho_f(\overline{g})\geq \|\overline{g}\|_{A/J^f}.$
\begin{lem}\label{lemdistquotien}Let $A$ be a distinguished $\bk$-affinoid algebra. Let $I$ be a reduced ideal of $A$. Denote by $\pi: A\twoheadrightarrow B:=A/I$ the quotient map.
Denote by $\|\cdot\|_B$ the residue norm on $B$ w.r.t. the spectral norm $\rho_A$ on $A$. Then for every $g\in B$, there exists $G\in \pi^{-1}(g)$ such that $\rho_A(G)=\|g\|_B.$
\end{lem}
By Lemma \ref{lemdistquotien}, there exists $g\in A$ whose image in $A/J^f$ is $\overline{g}$
such that $\rho_A(g)=\|\overline{g}\|_{A/J^f}.$
We may assume that $\rho_A(g)=\|\overline{g}\|_{A/J^f}=1$. We only need to show that $\rho_f(g)=1.$
Otherwise $\rho_f(g)<1$, then there exists $n\geq 1$ such that $\rho_A(f^n(g))<1.$
In other words, $\widetilde{g}\in \ker (\widetilde{f}^{n})=\widetilde{I}.$
Since $\widetilde{I}=\widetilde{J^f},$ there exists $w\in J^f$ and $h\in A^{\circ\circ}$ such that 
$g=w+h.$ Then $\|\overline{g}\|_{A/J^f}\leq \rho_A(h)<1,$ which is a contradiction.
Then we proved 2).

By $2)$, $\widetilde{A}\twoheadrightarrow \widetilde{A}/\widetilde{J^f}=\widetilde{A/J^f}.$ By 1), we have $\widetilde{I}=\widetilde{J^f}.$ This proves 3).
The definition of $J^f$ shows that $f(J^f)\subseteq J^f$. So we get 4). 
Since $\widetilde{f}|_{\widetilde{A}/\widetilde{I}}$ is an automorphism of $\widetilde{A}/\widetilde{I}$, $\widetilde{f|_{A/J^f}}=\widetilde{f}|_{\widetilde{A}/\widetilde{I}}$ is an automorphism of $\widetilde{A/J^f}=\widetilde{A}/\widetilde{I}.$
By 2) and the assumption that $A$ is distinguished, $A/J^f$ is distinguished.

\begin{pro}\label{proreducesursur}Let $A, B$ be two distinguished $\bk$-affinoid algebras. Let $g:A\to B$ be a morphism.
If the reduction $\widetilde{g}:\widetilde{A}\twoheadrightarrow \widetilde{B}$ is surjective. Then the morphism $g:A\to B$ is surjective.
\end{pro}

By Proposition \ref{proreducesursur}, $f|_{A/J^f}$ is surjective.
By \cite[ 6.3.1, Theorem 6]{Bosch1984}, $f|_{A/J^f}$ is injective. Then $f|_{A/J^f}$ is an isomorphism of $A/J^f$.
Then we get 5).

\medskip

We now construct the morphism $\psi: A/J^f\to A.$
Pick a bounded $\bk$-linear map $\chi:A/J^f\to A$ satisfying 
$\tau\circ \chi=\id.$ There exists $C>0$ such that $\rho_A(\chi(h))\leq C\rho_{A/J^f}(h).$
By Proposition \ref{proatleastlinear}, there exists $b\in (0,1)$ and $m\geq 1$ such that for all $g\in J^f$, $$\rho(f^m(g))\leq b\rho(g).$$
Then for every $n\geq m,g\in J^f$ we have $$\rho_A(f^n(g))\leq \rho_A(f^{[n/m]m}(g))\leq b^{[n/m]}\rho_A(g)\leq b^{n/2m}\rho_A(g).$$

For $n\geq 0$, we define the bounded $\bk$-linear maps 
$$\psi_n:=f^{n}\circ \chi\circ (f|_{A/J^f})^{-n}: A/J^f\to A.$$
Denote by $\|\cdot\|$ the operator norm of $\Hom_{\bk}(A/J^f, A)$.
We have $\|\psi_n\|\leq C.$
For every $h\in A/J^f,j\geq i\geq m$,
we have $$\psi_j(h)-\psi_i(h)=f^{j}\circ \chi\circ (f|_{A/J^f})^{-j}(h)-f^{i}\circ \chi\circ (f|_{A/J^f})^{-i}(h)$$
$$=f^{i}\circ\Big(f^{(j-i)}\circ \chi\circ (f|_{A/J^f})^{-j}(h)-\chi\circ (f|_{A/J^f})^{-i}(h)\Big).$$

Observe that $$\tau\Big(f^{(j-i)}\circ \chi\circ (f|_{A/J^f})^{-j}(h)-\chi\circ (f|_{A/J^f})^{-i}(h)\Big)$$
$$=\tau\circ f^{(j-i)}\circ \chi\circ (f|_{A/J^f})^{-j}(h)-\tau\circ\chi\circ (f|_{A/J^f})^{-i}(h)$$
$$=(f|_{A/J^f})^{(j-i)}\circ(\tau\circ \chi)\circ (f|_{A/J^f})^{-j}(h)-(\tau\circ\chi)\circ (f|_{A/J^f})^{-i}(h)$$
$$=0.$$

We have $f^{(j-i)}\circ \chi\circ (f|_{A/J^f})^{-j}(h)-\chi\circ (f|_{A/J^f})^{-i}(h)\in J^f$ and 
$$\rho_A\Big(f^{(j-i)}\circ \chi\circ (f|_{A/J^f})^{-j}(h)-\chi\circ (f|_{A/J^f})^{-i}(h)\Big)\leq C\rho_{A/J^f}(h).$$

It follows that $$\rho_A\Big(\psi_j(h)-\psi_i(h)\Big)=\rho_A\Big(f^{i}\big(f^{(j-i)}\circ \chi\circ (f|_{A/J^f})^{-j}(h)-\chi\circ (f|_{A/J^f})^{-i}(h)\big)\Big)$$
$$\leq b^{i/2m}\rho_{A/J^f}(h).$$
Then the sequence of operators $\psi_i,i\geq 0$ converges to a bounded $\bk$-linear map $\psi:A/J^f\to A$ with $\|\psi\|\leq C.$

For $g,h\in A/J^f,n\geq m$,
we have 
$$\psi_n(gh)-\psi_n(g)\psi_n(h)=f^{n}\circ (\chi\circ (f|_{A/J^f})^{-n}(gh)-\chi\circ (f|_{A/J^f})^{-n}(g)\chi\circ (f|_{A/J^f})^{-n}(h))$$
Observe that $$\tau(\chi\circ (f|_{A/J^f})^{-n}(gh)-\chi\circ (f|_{A/J^f})^{-n}(g)\chi\circ (f|_{A/J^f})^{-n}(h))$$$$=(f|_{A/J^f})^{-n}(gh)-(f|_{A/J^f})^{-n}(g)(f|_{A/J^f})^{-n}(h)=0.$$
We have $\chi\circ (f|_{A/J^f})^{-n}(gh)-\chi\circ (f|_{A/J^f})^{-n}(g)\chi\circ (f|_{A/J^f})^{-n}(h)\in J^f$ and of norm at most $C\rho_{A/J^f}(g)\rho_{A/J^f}(h).$
So we have 
$$\rho\Big(\psi_n(gh)-\psi_n(g)\psi_n(h)\Big)=\rho\Big(f^{n}\big(\chi\circ (f|_{A/J^f})^{-n}(gh)-\chi\circ (f|_{A/J^f})^{-n}(g)\chi\circ (f|_{A/J^f})^{-n}(h)\big)\Big)$$$$\leq Cb^{n/2m}\rho_{A/J^f}(g)\rho_{A/J^f}(h).$$
Let $n\to \infty$, we get $\psi(gh)=\psi(g)\psi(h).$ Then $\psi$ is a morphism of $\bk$-algebra.
Since
$$\tau\circ \psi_n=\tau\circ f^{n}\circ \chi\circ (f|_{A/J^f})^{-n}=\tau\circ(f|_{A/J^f})^{n}\circ (\tau\circ\chi)\circ (f|_{A/J^f})^{-n}=\id$$ for all $n\geq 0,$
$\tau\circ \psi=\id.$ 
we have $\psi_n\circ f|_{A/J^f}=f\circ \psi_{n-1}.$
Let $n\to \infty$, we get $\psi\circ f|_{A/J^f}=f\circ \psi.$ 
Then we get (i).

\medskip

For every $x\in \sM(A)$, $h\in A^{\circ},n\geq m$
we have 
$$\rho_A\Big(f^n\circ\psi\circ \tau(h)-f^n\circ\psi_n\circ \tau(h)\Big)\leq \|\psi-\psi_n\|\rho_A(h)\leq Cb^{n/2m}.$$
We note that $\tau(\psi_n\circ \tau(h)-h)=\tau(h)-\tau(h)=0$. Then we $\psi_n\circ \tau(h)-h\in J^f$ and its norm is at most $C.$
It follows that 
$$\rho_A\Big(f^n\circ\psi_n\circ \tau(h)-f^n(h)\Big)\leq b^{n/2m}C.$$
Then for every $x\in \sM(A)$, we have
$$||f^n(h)(x)|-|(\psi\circ\tau\circ f^n)(h)(x)||\leq |\Big(f^n(h-\psi(\tau(h)))\Big)(x)|$$$$\leq \rho_A\Big(f^n\circ\psi\circ \tau(h)-f^n(h)\Big)$$
$$\leq \max\{\rho_A(f^n\circ\psi\circ \tau(h)-f^n\circ\psi_n\circ \tau(h)),  \rho_A(f^n\circ\psi_n\circ \tau(h)-f^n(h)) \}$$$$\leq Cb^{n/2m}.$$

Now we only need to prove the uniqueness of $\psi.$
Assume that there is another morphism $\psi_1: A/J^f\to A$ satisfying $\psi_1|_{A/J^f}:=\tau\circ \psi_1=\id$ and  $\psi_1\circ f|_{A/J}=f\circ \psi_1$. We want to show $\psi=\psi_1.$
Since  $f\circ \psi=\psi\circ f|_{A/J^f}$, for every $n\geq 0$, we have $f^n\circ \psi=\psi\circ f|_{A/J^f}^n.$
Then we have $$f\circ \psi\circ f|_{A/J^f}^{-n}=\psi.$$
The same, we get $$f^n\circ \psi_1\circ f|_{A/J^f}^{-n}=\psi_1.$$

For every $h\in A/J^f, n\geq 0$, we have $$\tau(\psi\circ f|_{A/J^f}^{-n}(h)-\psi_1\circ f|_{A/J^f}^{-n}(h))=0.$$
Then we have $\psi\circ f|_{A/J^f}^{-n}(h)-\psi_1\circ f|_{A/J^f}^{-n}(h)\in J_f.$
Then for every $n\geq m,$  
$$\rho_A(\psi(h)-\psi_1(h))=\rho_A\Big(f^n\big(\psi\circ f|_{A/J^f}^{-n}(h)-\psi_1\circ f|_{A/J^f}^{-n}(h)\big)\Big)$$
$$\leq b^{n/2m}\rho_A(\psi\circ f|_{A/J^f}^{-n}(h)-\psi_1\circ f|_{A/J^f}^{-n}(h))\leq b^{n/2m}\rho_A(h).$$
Let $n\to \infty$, we get $\psi(h)-\psi_1(h)=0$. It implies that $\psi=\psi_1.$ We get 6).
\endproof

\proof[Proof of Lemma \ref{lemdecompline}]
By \cite[Chapter 2.3, Proposition 3]{Bosch2014}, there exists a Bald subring $R$ of $\bk^{\circ}$ such that all coefficients of $F, G_i, i=1,\dots,s$ and $K_i, i=1,\dots,m$ are contained in $R.$
After localizing  $R$ by all elements of norm $1$, we may assume that $R$ is a B-ring. Moreover, after taking completion, we may assume that $R$ is complete. 
Then $\widetilde{R}=R^{\circ}/R^{\circ\circ}$ is a subfield of $\widetilde{\bk}$.

\medskip
 
We have $\widetilde{K}\subseteq \widetilde{I}_1$ and $\widetilde{\sT}=\widetilde{I}_1+\widetilde{F}^j(\widetilde{\sT}), j\geq 1.$
There is a basis $\widetilde{E}_i, i\in S_1$ of $\widetilde{K}$ such that 
for all $i\in S_1$, $\widetilde{E}_i$ takes form $\widetilde{K}_{n_i}\widetilde{T}^{I_i}$ for some $n_i\in \{1,\dots, m\}$ and some multi-index $I_i.$
Since $\widetilde{I}_1$ is spanned by $\widetilde{G_l}\widetilde{T}^I, l=1,\dots s, I\in \Z_{\geq 0}^s$, there exist $\widetilde{E_i}, i\in S_2$ such that 
$\widetilde{E}_i, i\in S_1\sqcup S_2$ is  a basis of $\widetilde{I}_1$ and for all $i\in S_2$, $\widetilde{E}_i$ takes form $\widetilde{G}_{n_i}\widetilde{T}^{I_i}$ for some $n_i\in \{1,\dots, s\}$ and some multi-index $I_i.$

For every $j\geq 1$, 
since $\widetilde{\sT}=\widetilde{I}_1+\widetilde{F}^j(\widetilde{\sT}),$ and $\widetilde{F}^j(\widetilde{T}^I), I\in \Z_{\geq 2}^r$ spans $\widetilde{F}^j(\widetilde{\sT})$,
there exist $\widetilde{E_i^j}, i\in S_3$ such that 
$\widetilde{E}_i, i\in S_1\sqcup S_2, \widetilde{E}_i^j\in S_3$ is  a basis of $\widetilde{\sT}$ and for all $i\in S_3$, $\widetilde{E}_i^j$ takes form $\widetilde{F}^j(\widetilde{T}^{I_i^j})$ for some multi-index $I_i^j.$

\medskip

Define $E_i:=K_{n_i}{T}^{I_i}$ for $i\in S_1$,   $E_i:=G_{n_i}{T}^{I_i}$ for $i\in S_2$, $E_i^j:=F^j({T}^{I_i^j})$ for $i\in S_3,j\geq 1$ and $S:=S_1\sqcup S_2\sqcup S_3$.
We note that $E_i, E_n^j\in R\{T_1,\dots,T_r\}$ for $i\in S_1\sqcup S_2, n\in S_3, j\geq 1$.
Now \cite[Theorem 6]{Bosch2014} implies that for every $j\geq 1$, $E_i, i\in S:=S_1\sqcup S_2, E_n^j, n\in S_3$ forms an orthonormal basis of $T$, which concludes the proof.
\endproof

\proof[Proof of Lemma \ref{lemdistquotien}]
Pick a distinguished epimorphism $$\phi:\sT:=\bk\{T_1,\dots,T_r\}\twoheadrightarrow A.$$
Then the spectral norm $\rho_A$ equals to the residue norm w.r.t. $\rho_{\sT}.$
It implies that the norm $\|\cdot\|_B$ on $B$ is the residue norm w.r.t. $\rho_{\sT}.$
So we may assume that $A=\sT$ and conclude the proof by \cite[Chapter 2.3, Corollary 7]{Bosch2014}.
\endproof

\proof[Proof of Proposition \ref{proreducesursur}]
Let $\phi_A:\sT_A:=\bk\{T_1,\dots, T_r\}\twoheadrightarrow A$ be a distinguished epimorphism. 
Let $\phi_B:\sT_B:=\bk\{U_1,\dots, U_s\}\twoheadrightarrow B$ be a distinguished epimorphism. 
There exists a morphism $F:\sT_A\to \sT_B$ such that $g\circ\phi_A=\phi_B\circ F.$
Denote by $K$ the kernel of $\phi_B.$
By \cite[Corollary 7]{Bosch2014}, we may write $K=(K_1,\dots,K_m)$ where $\rho_B(K_i)=1, i=1,\dots,m$ such that $\widetilde{K}=(\widetilde{K}_1,\dots,\widetilde{K}_m)$ and $K^{\circ}=\sum_{i=1}^mK_i\sT_B^{\circ}.$

By \cite[Chapter 2.3, Proposition 3]{Bosch2014}, there exists a Bald subring $R$ of $\bk^{\circ}$ such that all coefficients of $F$ and $K_i, i=1,\dots,m$ are contained in $R.$
After localizing  $R$ by all elements of norm $1$, we may assume that $R$ is a B-ring. Moreover, after taking completion, we may assume that $R$ is complete. 
Then $\widetilde{R}=R^{\circ}/R^{\circ\circ}$ is a subfield of $\widetilde{\bk}$.

\medskip

We have a base $\widetilde{E}_i, i\in S_1$ of $\widetilde{K}$ such that 
for all $i\in S_1$, $\widetilde{E}_i$ takes form $\widetilde{K}_{j_i}\widetilde{U}^{I_i}$ for some $j_i\in \{1,\dots, m\}$ and some multi-index $I_i.$

Since $\widetilde{f}$ and $\widetilde{\phi}_A$ are surjective, $\widetilde{\sT_B}=\widetilde{K}+\widetilde{F}(\widetilde{\sT_A}).$
We note that $\widetilde{F}(\widetilde{T}^I), I\in \Z_{\geq 2}^r$ spans $\widetilde{F}(\widetilde{\sT_A}).$
There exist $\widetilde{E_i}, i\in S_2$ such that 
$\widetilde{E}_i, i\in S_1\sqcup S_2$ is  a basis of $\widetilde{\sT_A}$ and for every $i\in S_2$, $\widetilde{E}_i$ takes form $\widetilde{F}(\widetilde{T}^{I_i})$ for some multi-index $I_i.$

\medskip

Define $E_i:=K_{j_i}{U}^{I_i}$ for $i\in S_1$,   $E_i:=F({T}^{I_i})$ for $i\in S_2$.
We note that $E_i\in R\{U_1,\dots,U_s\}$ for $i\in S_1\sqcup S_2$.
Now \cite[Chapter 2.3, Theorem 6]{Bosch2014} implies that $E_i, i\in S:=S_1\sqcup S_2$ forms an orthonormal basis of $\sT_B$.
Since $E_i\in K$ for $i\in S_1$ and $E_i\in F(\sT_A)$ for $i\in S_2$, we get 
$\sT_B=K+F(\sT_A).$
Then $B=F(\sT_A)/(K\cap F(\sT_A))$. It implies that $g$ is surjective.
\endproof

%

\medskip

\begin{pro}\label{proideadisofdisred}Assume that $A$ is distinguished. Let $J$ be a reduced ideal of $A$. Assume that
the residue norm $\|\cdot\|_{A/J}$ on $A/J$ w.r.t. $\rho_A$ equals to $\rho_{A/J}$.
Let $g_1,\dots, g_m$ be elements in $J^{\circ}$ such that their reductions $\widetilde{g_1}, \dots,\widetilde{g_m}$ generate 
$\widetilde{J}.$ 
Then $g_1,\dots,g_m$ generate $J.$ 
\end{pro}

\proof[Proof of Proposition \ref{proideadisofdisred}]
Since $\|\cdot\|_{A/J}=\rho_{A/J}(\cdot)$, we have 
$\widetilde{A/J}=\widetilde{A}/\widetilde{J}.$

Pick a distinguished epimorphism $\phi:\sT:=\bk\{T_1,\dots,T_r\}\twoheadrightarrow A.$
Denote by $\rho_{\sT}$ the spectral norm on $\sT$.
Then the spectral norm $\rho_A$ on $A$ is the residue norm on $A$ w.r.t. $\rho_{\sT}.$ Set $I:=\ker(\phi).$ 
Pick $F_1,\dots,F_s$ in $I^{\circ}$ such that  their reductions $\widetilde{F_1}, \dots,\widetilde{F_s}$ generate $\widetilde{I}.$ 
Since $\phi$ is distinguished, by \cite[Chapter 2.3, Corollary 7]{Bosch2014}, for every $i=1,\dots,m$, there exists $G_1,\dots, G_m\in T$ such that $\phi(G_i)=g_i, i=1,\dots,m$ and
$\rho_{\sT}(G_i)=\rho_A(g_i)\leq 1, i=1,\dots,m.$
We note that  $F_1,\dots, F_s, G_1,\dots,G_m\in \phi^{-1}(J).$
We only need to show that $F_i,i=1,\dots s, G_i,i=1,\dots, m$  generate $\phi^{-1}(J)$.

Denote by $\psi:\sT\twoheadrightarrow A/J=\sT/\phi^{-1}(J)$ the composition of $\phi$ and the quotient morphism $A\to A/J.$
Since $\phi$ is distinguished and  $\|\cdot\|_{A/J}=\rho_{A/J}(\cdot)$, $\psi$ is distinguished.
Then
$\widetilde{A}/\widetilde{J}=\widetilde{A/J}=\widetilde{\sT}/\widetilde{\phi^{-1}(J)}.$
Since $\phi$ is distinguished, we have 
$\widetilde{A}=\widetilde{\sT/I}=\widetilde{\sT}/\widetilde{I}.$
Then we get 
$\widetilde{T}/\widetilde{\phi^{-1}(J)}=\widetilde{T}/\widetilde{\phi}^{-1}(\widetilde{J}),$
which implies that $\widetilde{\phi^{-1}(J)}=\widetilde{\phi}^{-1}(\widetilde{J}).$
Observe that 
$\widetilde{F}_1,\dots, \widetilde{F}_s, \widetilde{G}_1,\dots,\widetilde{G}_m$  generate $\widetilde{\phi}^{-1}(\widetilde{J})$.
Then the proof of  \cite[Chapter 2.3, Corollary 7]{Bosch2014},  shows that 
$F_1,\dots, F_s, G_1,\dots,G_m$  generate $\phi^{-1}(J)$, which concludes the proof.
\endproof

\subsection{The Zariski density of orbits}\label{subsectionzardenan}
In this section, we assume that $\bk=\C_p.$ 
Let $K_p\subseteq \C_p$ be a finite field extension of $\Q_p.$

Let $f: \D^2:=\sM(\bk\{x,y\})\to \D^2$ be an endomorphism defined over $K_p$ whose reduction $\widetilde{f}:\A^2_{\widetilde{\bk}}\to \A^2_{\widetilde{\bk}}$ takes form 
 $$\widetilde{f}:(x,y)\mapsto (\widetilde{a}x+\widetilde{b},0)$$ where $\widetilde{a}\in \widetilde{\bk}\setminus \{0\}, \widetilde{b}\in \widetilde{\bk}.$

By Theorem \ref{thmstattrainv} and Proposition \ref{proideadisofdisred},
there exists $g\in K_p \{x,y\}$ taking form $g=y+h$ where $h\in K_p\{x,y\}^{\circ\circ}$
such that $J^f=(g).$ Set $Y:=\sM(\bk\{x,y\})/(g)$. We have $Y\simeq \D^1.$
There exists a unique morphism $\psi:\D^2\to Y$ satisfying $\psi|_Y=\id$ and $f|_Y\circ \psi=\psi\circ f.$
There exists $C>0, \beta\in (0,1)$ such that 
for every $F\in \bk\{x,y\},$
$x\in \D^2$ and $n\geq 0$, we have $$||F(f^n(x))|-|F(f^n(\psi(x)))||\leq C\beta^n\rho(F).$$

\begin{rem}\label{remattractordete}
Assume that $f$ takes form 
$(x,y)\mapsto (ax+b+P, yQ)$ where $|a|=1, P,Q\in \bk\{x,y\}, \rho(P),\rho(Q)<1.$
Then for $n\geq 1$, we have 
$$(f^*)^n(y)=yQf^*(Q)\cdots (f^*)^n(Q).$$
It follows that $\rho((f^*)^n(y))\leq \rho(Q)^n.$ In particular, we have $y\in J^f.$
By Proposition \ref{proideadisofdisred}, we have $J^f=(y).$
\end{rem}
\begin{exe}\label{exefixedcase}
Assume that $f$ takes form 
$(x,y)\mapsto (x+yP, yQ)$ where $P,Q\in \bk\{x,y\}, \rho(P),\rho(Q)<1.$
By Proposition \ref{proideadisofdisred}, we have $J^f=(y).$ So we have $Y:=\{y=0\}.$
In this case we may compute the morphism $\psi: \D^2\to Y$ explicitly.
Follows the proof Theorem \ref{thmstattrainv}, $\psi$ equal to $\lim\limits_{n\to\infty} f^n$, which is defined by 
$$(x,y)\mapsto (\lim_{n\to \infty}(f^*)^n(x), \lim_{n\to \infty}(f^*)^n(y))=(x+\sum_{i\geq 1}(f^*)^i(y)(f^*)^i(P), 0).$$
We note that $\rho((f^*)^i(y)(f^*)^i(P))\leq \rho((f^*)^i(y))\leq \rho(Q)^i.$
In particular, for every $c\in \bk^{\circ}$, $\psi^{-1}((c,0))=\{x+\sum_{i\geq 1}(f^*)^i(y)(f^*)^i(P)=c\}.$
By implicit function theorem, $\psi^{-1}((c,0))\simeq \D.$
\end{exe}

\begin{pro}\label{profordtwozd}
Assume that $f^n|_Y\neq \id$ for every $n\geq 1$ and $f^{-1}(Y)\neq \D^2$. Then there exists a nonempty strict affinoid subdomain $V$ of $\D^2$ such that for every $o\in V(\bk)$, the orbit $O_f(o)$ of $o$ is Zariski dense in $\D^2.$
\end{pro}

\rem\label{remcomparezaralan}
Assume that $X$ is an projective surface over $\bk$. Denote by $X^{\an}$ the analytification of $X$.  Then we have a natural morphism $\pi_X:X^{\an}\to X.$
We note that $\pi_X$ gives a bijection between $X^{\an}(\bk)$ and $X(\bk).$

Assume that there exists a strict affinoid subdomain $U$ of $X^{\an}$. 
Then the Zariski topology of $U$ is finer than the pullback by $\pi_X|_{U}$ of the Zariski topology of $X.$
So if $X$ is irreducible and a set $S$ of $U(\bk)$ is Zariski dense in $U(\bk)$, then $\pi_X(U(\bk))$ is Zariski dense in $X.$

\endrem

\proof[Proof of Proposition \ref{profordtwozd}]
Fix an identification $\D^1=\sM(\bk\{T\})\simeq Y.$
The reduction of $f|_Y$ takes form $\widetilde{f|_Y}: T\mapsto \widetilde{a}T+\widetilde{b}.$
There exists $m\geq 0$ such that $\widetilde{a}^m=1.$ After replacing $f$ by $f^{mp}$, we may assume that $\widetilde{a}=1, \widetilde{b}=0.$ Then we may assume that $\widetilde{f|_Y}=\id$.

Denote by $\Delta_{f|_Y}:=f|_Y^*-\id: \bk\{T\}\to \bk\{T\}$ the difference operator which is a bounded linear operator on the Banach space $\bk\{T\}.$
Denote by $\|\Delta_{f|_Y}\|$ th operator norm of $\Delta_{f|_Y}.$
Since $\widetilde{{f|_Y}}=\id$, and $Y,f|_Y$ are defined over a discrete valuation field $K_p$, we have $\|\Delta_{f|_Y}\|< 1$.
By \cite[Remark 4]{Poonen2014}, there exists $r\geq 1$ such that $\|\Delta_{({f|_Y})^r}\|< p^{-2}.$
After replacing $f$ by $f^r$, we may assume that $\|\Delta_{{f|_Y}}\|< p^{-2}.$
Then \cite[Theorem 1]{Poonen2014} shows that the set of preperiodic  points of ${f|_Y}$ in $\D^1(\bk)$ is the 
set of fixed points $\Fix({f|_Y})$ of ${f|_Y}$ in $\D^1(\bk).$ Since $f|_Y^n\neq \id$ for every $n\geq 1$, $\Fix({f|_Y})$ is finite. 
Since $f^{-1}(Y)\neq \D^2$, $f^{-1}(Y)$ is a union of finitely many irreducible curves. Let $Y_1$ to be the union of all irreducible components of $f^{-1}(Y)$ except $Y.$
Then $Y\cap Y_1$ is a finite union of closed points.

\medskip

Pick $a\in Y(\bk)\setminus (\Fix(f|_Y)\cup (Y\cap Y_1)).$ There exists $s\geq 1$ such that 
the ball $B:=\{t\in Y=\D^1|\,\, |(T-a)(t)|\leq p^{-s}\}$ does not meet $\Fix(f|_Y)\cup (Y\cap Y_1).$
Observe that $B$ is a strict affinoid subdomain of $Y$.
By \cite[Remark 4]{Poonen2014}, after replacing $f$ by some positive iterate, we may assume that $\|\Delta_{f|_Y}\|< p^{-s}.$
It follows that the ball $B:=\{t\in \D^1|\,\, |(T-a)(t)|\leq p^{-s}\}$ is invariant under $f|_Y.$

\medskip

We note that $Y\cap \psi^{-1}(B)\cap Y_1=B\cap Y_1=\emptyset.$ 
There exists $l\geq 1$ such that 
$$Y^l\cap \psi^{-1}(B)\cap Y_1=\emptyset$$ where $Y^l$ is the affinoid subdomain $\{t\in \D^2|\,\,|g(x)|\leq p^l \}.$
Observe that $f(Y^l)\subseteq Y^l.$
Then $Y^l\cap  \psi^{-1}(B)$ is an analytic subdomain of $\D^2$ which is invariant by $f$.
Moreover $(Y^l\cap  \psi^{-1}(B))\setminus Y$ is also invariant by $f.$
Since $W:=(Y^l\cap  \psi^{-1}(B))\setminus Y$ contains a strict affinoid subdomain of $X$, we only need to show that for every $o\in W(\bk)$,
the orbit of $o$ is Zariski dense. 
Otherwise, there exists $o\in W(\bk)$ such that the Zariski closure $Z$ of $O_f(o)$ satisfies $\dim(Z)\leq 1.$ Since $f^n(o)\not\in Y, n\geq 0$ and  tends to $Y$, we have $\dim(Z)=1.$
After replacing $f$ by a positive iterate and $o$ by $f^s(o)$ for some $s\geq 0$, we may assume that $Z$ is an irreducible curve.
The intersection $Z\cap Y$ is a finite set of closed points. Since $f|_Y$ is an automorphism and $f(Z)\subseteq Z$, 
every point in $Z\cap Y$ is periodic. It follows that $(Z\cap Y)(\bk)\subseteq \Fix(f|_Y).$
By \cite[Remark 4]{Poonen2014}, we have $f^{p^n}(\psi(o))\to \psi(o)$ when $n\to \infty.$ 
Assume that $Z$ is defined by $g_1,\dots,g_r.$ For every $i=1,\dots,r$,  $n\geq 0$
we have 
$$|g_i(f^{p^n}(\psi(o)))|=||g_i(f^{p^n}(\psi(o)))|-|g_i(f^{p^n}(o))||\leq C\beta^{p^n}\rho(g_i).$$
Let $n\to \infty$, we have $g_i(\psi(o))=0.$ It follows that $\psi(o)\in (Z\cap Y)(\bk)\subseteq \Fix(f|_Y)$, which is a contradiction.
Then we concludes the proof.
\endproof


\section{Appendix B (Joint work with Thomas Tucker): The Zariski dense orbit conjecture for endomorphisms of $(\P^1)^N$}
Let $\bk$ be an algebraically closed field of characteristic zero.

The aim of this appendix is to prove Theorem \ref{thmzaridenseorbitsurfendop1n}. By Corollary \ref{coradelicstrong}, we only need to prove its adelic version as follows.
\begin{thm}\label{thmzaridenseorbitsurfendop1nadelic}Assume $\trd_{\Q}\bk<\infty.$
Let $f:(\P^1)^N\to (\P^1)^N$ be a dominant endomorphism of $(\P^1)^N, N\geq 1$.
Then the pair $((\P^1)^N,f)$ satisfies the adelic ZD-propety.
\end{thm}

Now assume $\trd_{\Q}\bk<\infty.$

\medskip

For $i=1,\dots,N$, denote by $\pi_i: (\P^1)^N\to \P^1$ the projection to the $i$-th coordinate. 
Denote by $H_i, i=1,\dots, N$ the class in $N^1((\P^1)^N)$ represented by $\pi_i^{*}(\sO_{\P^1}(1)).$
The Nef cone $\sC$ of $(\P^1)^N$ in $N^1((\P^1)^N)$ is the convex cone spanned by $H_1,\dots, H_N.$
Since $f^*, f_*$  preserve the Nef cone, after replacing $f$ by some positive iteration, we may assume that $f^*H_i$ is some multiple to $H_i$.
Then $f$ preserves all $\pi_i, i=1,\dots,N.$
We may assume that $f$ takes form $(x_1,\dots, x_N)\mapsto (f_1(x_1),\dots, f_N(x_N)),$ where $f_i$ is an endomorphism of $\P^1$ of degree at least $1.$

\subsection{Exceptional endomorphisms of curves}
Let $C$ be an irreducible projective curve.
Let $g$ be an endomorphism $g: C\to C$ with $\deg g\geq 2$. Then $C$ is either $\P^1$ or an elliptic curve.

We say that $g$ is \emph{of Latt\'es type}, if it semi-conjugates to an endomorphism of an elliptic curve i.e.  there exists an endomorphism of an elliptic curve $h: E\to E$ and a finite morphism
$\pi: E\to C$ such that $f\circ\pi=\pi\circ h.$

We say that $g$ is \emph{of monomial type}, if it semi-conjugates to an endomorphism of a monomial map i.e.  there exists a monomial endomorphism $h: \P^1\to \P^1$ taking form $x\mapsto x^{\pm d}, d\geq 2$ and a finite morphism
$\pi: \P^1\to C$ such that $f\circ\pi=\pi\circ h.$ We note that in this case $C\simeq \P^1.$

We say that $g$ is \emph{exceptional} if it is of  Latt\'es type or monomial type.
Otherwise, it is said to be \emph{nonexceptional}.

For every endomorphism $g: \P^1\to \P^1$ of $\deg g\geq 2$, it has
exactly one type in Latt\'es, monomial, and nonexceptional.  Moreover,
the types of $g^n,n\geq 1$ are the same.

The following facts are well known.
\begin{points}
\item If two endomorphisms of curves are semi-conjugacy, then they have the same type.
\item If there is a nonzero rational differential form $\omega$, such that $g^*\omega=\mu \omega$ for some $\mu\in \bk^*$, then $g$ is exceptional.
\item  If $g$ has an exceptional point i.e. a point in $C$ whose inverse orbits is finite, then $g^2$ is polynomial. In particular, when $g$ is of monomial type, $g^2$ is polynomial.
\end{points}

\subsection{Invariant subvarieties}
For every $l=1,\dots,N$, let $S_l$ be the set of subsets of $\{1,\dots, N\}$ of $l$ elements.
For every $I\subseteq \{1,\dots, N\}$, the ordering in $I$ is induced by the ordering in $\{1,\dots, N\}$.
For every subset $I$ of $\{1,\dots, N\}$, denote by $\pi_I: (\P^1)^N\to (\P^1)^{|I|}$ the projection $(x_i)_{i=1,\dots,N}\mapsto (x_i)_{i\in I}.$
Let $f_I:(\P^1)^{|I|}\to (\P^1)^{|I|}$ be the endomorphism $(x_i)_{i\in I}\mapsto (f_i(x_i))_{i\in I}.$
We have $\pi_I\circ f=f_I\circ \pi_I.$

\medskip

Every hypersurface $V$ of $(\P^1)^N$, as a divisor,  is linearly equivalent  to $\sum_{i=1}^N d_iH_i$ for some $d_i\geq 0, i=1,\dots,N.$ 
It is ample if and only if $d_i>0$ for all $i=1,\dots,m.$ Moreover, $d_i=0$ if and only if $\pi_{\{1,\dots,N\}\setminus \{i\}} |_V: V\to (\P^1)^{N-1}$ is not dominant.

\medskip

The following results on the invariant subvarieties was obtained in \cite{Medvdev} using model theory. 
When $\bk=\overline{\Q}$, it was also obtained by Ghioca, Nguyen and Ye in \cite[Theorem 1.2]{Ghioca2018c}, as a consequence of their solution of the Dynamical 
Manin-Mumford Conjecture in this case. Here we give a new proof which is purely geometric.
\begin{pro}\label{proinvsubvspl}Assume that $N\geq 2$, $\deg f_i\geq 2, i=1,\dots,N$ and all $f_i,i=1,\dots,N$ are nonexceptional. Let $V$ be a proper irreducible subvariety of $(\P^1)^N$ which is invariant under  $f$.
Then there exists $I\in S_2$ such that $V\subseteq \pi_I^{-1}(C)$ where $C$ is an $f_I$-invariant curve in $(\P^1)^2.$
\end{pro}

To prove this, we need the following lemmas.
\begin{lem}\label{lemirr}Let $V$ be an irreducible hypersurface of $(\P^1)^N$.
Assume that $V$ is ample as a divisor.
Then for every $J\in S_l, l\leq N-2$ and a general point $z\in (\P^1)^{|J|}$, $(\pi_J|_V)^{-1}(z)$ is irreducible.
\end{lem}
\proof[Proof of Lemma \ref{lemirr}]
We may assume that $N\geq 3$ and $J=\{1,\dots,l\}.$
Let $\mu:Y\to V$ be a desingularization of $V$. Set $\mu_i:=\pi_i\circ \mu, i=1,\dots,N.$

By Theorem of Bertini, for general $a_1\in \P^1$, $\mu_1^{-1}(a_1)$ is smooth.
Since $V|_{\pi_1^{-1}(a)}$ is ample and $\dim \pi_1^{-1}(a)\geq 2$, 
by \cite[Proposition, Page 67]{Hartshorne1970}, $V|_{\pi_1^{-1}(a)}$ is connected. Since $\mu$ is birational, by Zariski's main theorem, every fiber of $\mu$ is connected.
This implies that $\mu_1^{-1}(a_1)=\mu^{-1}(\pi_1^{-1}(a))$ is smooth and connected. So $\mu_1^{-1}(a_1)$ is smooth and irreducible
of dimension $N-2$. Repeat this argument, we get that,
for general $a_2\in \P^1$, $\mu_1^{-1}(a_1)\cap \mu_2^{-1}(a_2)$ is irreducible and smooth of dimension $N-3$;\dots;
for general $a_{l}\in \P^1$, $\mu_1^{-1}(a_1)\cap\dots\cap  \mu_l^{-1}(a_l)$ is irreducible and smooth of dimension $N-l-1$.
Then the geometric generic fiber of $\pi_J\circ \mu$ is irreducible and smooth.  So for a general point $z\in (\P^1)^{|J|}$, $(\pi_J\circ \mu)^{-1}(z)$ is irreducible and smooth.
Then, $(\pi_J|_V)^{-1}(z)=\mu((\pi_J\circ \mu)^{-1}(z))$ is irreducible.
\endproof

\medskip

\begin{lem}\label{lemsurn3}Assume that $N=3$, $\deg f_i\geq 2, i=1,2,3$ and $f_2$ is nonexceptional. Let $V$ be a proper irreducible hypersurface of $(\P^1)^3$ which is invariant under  $f$.
Assume that $\pi_{\{1,3\}}(V)=(\P^1)^2,\pi_{\{2,3\}}(V)=(\P^1)^2$, then $\pi_{\{1,2\}}(V)\neq (\P^1)^2.$
\end{lem}

\proof[Proof of Lemma \ref{lemsurn3}] 
Assume that $\pi_{\{1,2\}}(V)=(\P^1)^2.$ Then $V$ is ample.

Set $g:=f|_V.$ Set $P_i:=\pi_i|_V, i=1,2,3.$
Then we have three nonzero rational differential forms $dP_i,i=1,2,3$ on $V.$
We have $$g^*dP_i=d(P_i\circ g)=d(f_i\circ P_i)=P_i^*(f_i')dP_i.$$

Since $\pi_{\{1,2\}}(V)=(\P^1)^2$, $dP_1,dP_2$ are linearly independent at a general point in $V.$
So there are rational functions $G_1,G_2\in \bk(V)$ such that $dP_3=G_1dP_1+G_2dP_2.$
Since $\pi_{\{1,3\}}(V)=(\P^1)^2,\pi_{\{2,3\}}(V)=(\P^1)^2$, $G_1,G_2$ are nonzero.

For every $a\in \P^1$, define  $V_a:=\pi^{-1}_1(a)\cap V.$ 
By Lemma \ref{lemirr}, there exists a nonempty Zariski open subset $U$ of $\P^1$, such that for every $a\in U$, $V_a$ is irreducible.
After shrinking $U$, we may assume that for every $a\in U$, $G_1|_{V_a}$ and $G_2|_{V_a}$ are nonzero.

Since the set of critical periodic points of $f_1$ is finite, after shrinking $U$, we may assume that for all $a\in U\cap \Per(f_1),$ $(f_1^n)'(a)\neq 0$ for every $n\geq 0.$
We note that $\Per(f_1)$ is Zariski dense in $\P^1$. 
Pick $a\in U\cap \Per(f_1).$
There exists $s\geq 1$, such that $f_1^s(a)=a$. Then $V_a$ is $g^s$-invariant.

We have 
$$(P_3|_{V_a}^*(f_3^s)')dP_3|_{V_a}=(g^s|_{V_a})^*dP_3|_{V_a}=(g^s|_{V_a})^*(G_1|_{V_a}dP_1|_{V_a}+G_2|_{V_a}dP_2|_{V_a})$$
$$=((g^s|_{V_a})^*G_1|_{V_a})((f_1^s)'(a))dP_1|_{V_a}+((g|_{V_a}^s)^*G_2|_{V_a})(P_2|_{V_a}^*(f_2^s)')dP_2|_{V_a}.$$

It follows that 
$$\frac{((g^s|_{V_a})^*G_1|_{V_a})((f_1^s)'(a))}{((g|_{V_a}^s)^*G_2|_{V_a})(P_2|_{V_a}^*(f_2^s)')}=\frac{G_1|_{V_a}}{G_2|_{V_a}}.$$
Then
$$(g^s|_{V_a})^*(\frac{G_2|_{V_a}}{G_1|_{V_a}}dP_2|_{V_a})=(g^s|_{V_a})^*(\frac{G_2|_{V_a}}{G_1|_{V_a}})(P_2|_{V_a}^*(f_2^s)')dP_2|_{V_a}$$
$$=((f_1^s)'(a))\frac{G_2|_{V_a}}{G_1|_{V_a}}dP_2|_{V_a}.$$
Then $g^s|_{V_a}$ is exceptional.
Since $g^s$ semi-conjugates to $f_2^s$ and $f_2$ is nonexceptional, we get a contradiction. 
\endproof

\proof[Proof of Proposition \ref{proinvsubvspl}]
We do the proof by induction on $N\geq 2.$
When $N=2$, there is nothing to prove.

Assume that Proposition \ref{proinvsubvspl} is known for $N=2,\dots,l.$ We need to show it when $N=l+1\geq 3.$

We first assume that there exists $J\in S_l$ such that $\pi_J(V)\neq (\P^1)^{|J|}$. Then it is a proper irreducible subvariety of $(\P^1)^{|J|}$ which is invariant under  $f^J$.
We conclude the proof by the induction hypothesis.

Now we may assume that for every $J\in S_l$, $\pi_J(V)=(\P^1)^{|J|}$. Then $\dim V= l,$ and for every $J\in S_t, 1\leq t\leq l$, $\pi_J(V)=(\P^1)^{|J|}$.
So $V$ is ample.

When $l=2,N=3$, we conclude the proof by Lemma \ref{lemsurn3}.
Now we may assume that $l\geq 3.$
Set $K:=\{1,\dots,N\}$,$J:=\{4,\dots,N\}$ and $I:=\{1,2,3\}.$
By Lemma \ref{lemirr}, there exists a nonempty Zariski open subset $U$ of $(\P^1)^{N-3}$ such that  for every $a\in U$,
$V_a:=\pi_J^{-1}(a)\cap V$ is an irreducible surface.

For every $i\in J$,  $\Per(f_i)$ is Zariski dense in $\P^1$. Since the set $C_i$ of critical $f_i$-periodic points of is finite, $P_i:=\Per(f_i)\setminus C_i$ is Zariski dense in $\P^1$.
Then $(\prod_{i\in J}P_i)\cap U$ is Zariski dense in $(\P^1)^{|J|}.$
Pick $a\in (\prod_{i\in J}P_i)\cap U$, there exists $s\geq 1$ such that $f_J^s(a)=a.$ Then $V_a$ is invariant under $f_{I}^s.$
By Lemma \ref{lemsurn3}, there exists $i\in I$ such that $\pi_{I\setminus \{i\}}(V_a)\neq (\P^1)^2.$  
Pick $o\in  (\P^1)^2\setminus \pi_{I\setminus \{i\}}(V_a)$.  Then we have $(o,a)\subseteq (\P^1)^{N-1}\setminus \pi_{K\setminus \{i\}}(V).$ It follows that 
$\pi_{K\setminus \{i\}}(V)\neq (\P^1)^{N-1},$ which contradicts our assumption. We conclude the proof.
\endproof

\subsection{Proof of Theorem \ref{thmzaridenseorbitsurfendop1nadelic}}
\proof
By Theorem \ref{thmzardenseratadelic}, we may assume that $\deg f_i\geq 2, i=1,\dots,N.$
By Theorem \ref{thmprodendosup}, we may assume that $f_i^2, i=1,\dots,N$ are not polynomial. In particular, none of $f_i, i=1,\dots, N$ is of monomial type.

\medskip

By Proposition \ref{proinvpoly}, after replacing $f$ by some positive iterate,
there exists a nonempty adelic open subset $A$ of $(\P^1)^N(\bk)$ such that 
 for every $x\in A$, the Zariski closure $Z_x$ of the orbit $O_{f}(x)$ is irreducible.

\medskip

We may assume that there exists $0\leq s\leq N$ such that $f_i$ is nonexceptional for $i\leq s$ and it is of type Latt\'es for $i\geq s+1$.
Define $l(f):=\min\{s, N-s\}\geq 0.$

\medskip

We first treat the case $l(f)=0.$

If $s=0$, then all $f_i$ are of type Latt\'es. Then there exists an abelian variety $A'$, a dominant endomorphism $g: A'\to A'$ and a finite morphism $\pi: A'\to (\P^1)^N$ such that 
$f\circ \pi=\pi\circ g$. By Theorem \ref{thmendoabadelic}, the pair $(A',g)$ satisfies the AZD-property.
Then we conclude the proof by Lemma \ref{lemdesfiniteadelic}.

If $s=N,$ then all $f_i$ are nonexceptional.  For every $J\in S_2$, $f_J$ is an amplified endomorphism on $(\P^1)^2$, whose topological degree is strictly larger than $\la_1(f_J).$
By Corollary \ref{cortopdomampz}, there exists a nonempty adelic open subset $A_J$ of $(\P^1)^2$ such that the $f_J$-orbits of every point in $A_J$ are Zariski dense in $(\P^1)^2$.

Then for every point $x$ in the nonempty adelic open subset $A\cap (\cap_{J\in S_2}\pi_J^{-1}(A_J))$, the Zariski closure $Z_x$ of the orbit  of $x$ is irreducible, $f$-invariant and $\pi_J(Z_x)=(\P^1)^2$ for $J\in S_2.$
By Proposition \ref{proinvsubvspl}, $Z_x=(\P^1)^N,$ which concludes the proof.

\medskip

Now we do the proof by induction on  $l(f)N\geq 0.$
Assume that Theorem \ref{thmzaridenseorbitsurfendop1nadelic} holds when $0\leq l(f)N\leq m.$ We need to prove it when $l(f)N=m+1\geq 1.$

By the induction hypothesis, for every $J\in S_{N-1}$, 
there exists a nonempty adelic open subset $A_J$ of $(\P^1)^{N-1}(\bk)$, such that the $f_J$-orbit of every $x\in A_J$ is Zariski dense in $(\P^1)^{N-1}(\bk).$

Set $B:=A\cap(\cap_{J\in S_{N-1}}\pi_J^{-1}(A_J)).$ For every $x\in B$, the Zariski closure $Z_x$ of the orbit $O_{f}(x)$ is irreducible, invariant by $f$ and for every $J\in S_{N-1}$, we have $\pi_J(Z_x)=(\P^1)^{N-1}.$ Then we have $\dim Z_x\geq N-1,$ and for every $J\in S_t, 1\leq t\leq N-1$, $\pi_J(Z_x)=(\P^1)^{|J|}$.  In particular, if $Z_x\neq (\P^1)^N$, $Z_x$ is ample.

Assume that $\dim Z_x=N-1.$
We note that $f_1$ is nonexceptional and $f_N$ is of Latt\'es type.
Set $I:=\{2,\dots, N-1\}.$
Lemma \ref{lemirr} shows that there exists a nonempty Zariski open subset $U\subseteq (\P^1)^{N-2}$, such that 
for every $a\in U$, $(\pi_I|_{Z_x})^{-1}(a)$ is an irreducible curve.
Denote by $\Per(f_I)$ the set of periodic points of $f_I.$ It is Zariski dense in $(\P^1)^{N-2}.$
Pick $a\in \Per(f_I)\cap U.$ Then $C_a:=(\pi_I|_{Z_x})^{-1}(a)$ is an irreducible curve in $(\P^1)^2$, which is invariant under $f_{\{1,N\}}^s$ for some $s\geq 1.$
Then $f_{\{1,N\}}^s|_{C_a}: C_a\to C_a$ is an endomorphism of degree at least $2$. Denote by $p_i,i=1,N: (\P^1)^2\to \P^1$ the projection $(x_1,x_N)\mapsto x_i$. 
If $p_i(C_a)=\P^1$ for $i=1,2$, then $f_{\{1,N\}}^s|_{C_a}$ semiconjugates to both $f_1^s$ and $f_N^s.$ In particular, $f_1$ and $f_N$ has the same type, which is a contradiction.
It follows that $C_a$ is a fiber of $p_i$, $i=1$ or $N$. 
This implies that for all $b\in (\P^1)^{N-2},$ $C_b:=(\pi_I|_{Z_x})^{-1}(b)$ is a fiber of $p_i.$
Then, for $J:=\{1,\dots,N\}\setminus \{N+1-i\},$ we have $\pi_J(Z_x)\neq (\P^1)^{N-1},$ which contradicts our assumption.
It follows that $\dim Z_x=N$, which concludes the proof.
\endproof

\bibliography{dd}
\end{document}